\documentclass[12pt]{elsarticle}

\usepackage{graphics,graphicx}
\usepackage[usenames,dvipsnames]{xcolor}
\usepackage{amsmath}
\usepackage{latexsym}
\usepackage{amsfonts}
\usepackage{mathrsfs}
\usepackage{bm}
\usepackage[margin=1.5cm]{geometry}
\usepackage{multirow}
\usepackage{amsthm}
\usepackage{amssymb}
\usepackage{caption}
\usepackage{subcaption}
\usepackage{epstopdf}
\usepackage{float}
\usepackage[colorlinks,citecolor=blue,urlcolor=blue,bookmarks=false,hypertexnames=true]{hyperref}
\usepackage[title]{appendix}
\usepackage{tikz}
\usetikzlibrary{shapes.geometric, arrows}
\tikzstyle{startstop} = [rectangle, rounded corners, minimum width=3cm, minimum height=1cm,text centered, draw=black]
\tikzstyle{arrow} = [thick,->,>=stealth]
\graphicspath{ {Fig/} }

\usepackage[labelformat=simple]{subcaption}

\biboptions{sort&compress}
\newtheorem{note}{Note}
 \numberwithin{equation}{section}

\newtheorem{theorem}{Theorem}[section]
\newtheorem{corollary}{Corollary}[theorem]
\newtheorem{lemma}[theorem]{Lemma}
\newtheorem{remark}{Remark}[section]

\usepackage{tikz}

\usetikzlibrary{shapes.geometric, arrows, shapes}
\tikzstyle{startstop} = [rectangle, rounded corners, minimum width=3cm, minimum height=1cm,text centered, draw=black]
\tikzstyle{arrow} = [thick,->,>=stealth]

\newcommand{\tikzcircle}[2][red,fill=red]{\tikz[baseline=-0.5ex]\draw[#1,radius=#2] (0,0) circle ;}%
\newcommand\solidrule[1][1cm]{\rule[0.5mm]{8mm}{3pt}}
\newcommand\dashedrule{\mbox{%
\rule[0.5mm]{2.5mm}{2.5pt}\hspace{2mm}\rule[0.5mm]{2mm}{2.5pt}\hspace{2mm}\rule[0.5mm]{2mm}{2.5pt}}}

\makeatletter

\newcommand{\mathleft}{\@fleqntrue\@mathmargin0pt}
\newcommand{\mathcenter}{\@fleqnfalse}
\makeatother
\newlength{\bibitemsep}
\newlength{\bibparskip}
\let\oldthebibliography\thebibliography
\renewcommand\thebibliography[1]{%
  \oldthebibliography{#1}%
  \setlength{\parskip}{\bibitemsep}%
  \setlength{\itemsep}{\bibparskip}%
}

\usetikzlibrary{positioning, fit, calc}   
\tikzset{block/.style={draw, thick, text width=2cm ,minimum height=1.3cm, align=center},   
line/.style={-latex}     
}

\begin{document}

\makeatletter
\def\ps@pprintTitle{%
  \let\@oddhead\@empty
  \let\@evenhead\@empty
  \def\@oddfoot{\reset@font\hfil\thepage\hfil}
  \let\@evenfoot\@oddfoot
}
\makeatother

\begin{frontmatter}

\title{Nonlocal cooperative behaviour, psychological effects, and collective decision-making: an exemplification with predator-prey models
}


\author[inst1]{Sangeeta Saha}
\author[inst1]{Swadesh Pal}
\author[inst1,inst2]{Roderick Melnik}

\affiliation[inst1]{organization={MS2 Discovery Interdisciplinary Research Institute, Wilfrid Laurier University, Waterloo, Canada},
            }
\affiliation[inst2]{organization={BCAM - Basque Center for Applied Mathematics, E-48009, Bilbao, Spain}}

\date{}

\begin{abstract}
In bio-social models, cooperative behaviour has evolved as an adaptive strategy, playing multi-functional roles. One of such roles in populations is to increase the success of survival and reproduction of individuals and their families or social groups. Moreover, collective decision-making in cooperative behaviour is an aspect that is used to study the dynamic behaviour of individuals within a social group. In this paper, we have focused on population dynamics by considering a predator-prey model as our main exemplification, where the generalist predator has adopted a cooperative hunting strategy while consuming their prey. In particular, we have analyzed the dynamic nature of the system when a nonlocal term is introduced in the cooperation. First, the Turing instability condition has been studied for the local model around the coexisting steady-state, followed by the Turing and non-Turing patterns in the presence of the nonlocal interaction term. This work is also concerned with the existence of travelling wave solutions for predator-prey interaction with the nonlocal cooperative hunting strategy. Such solutions are reported for local as well as for nonlocal models. We have characterized the invading speed of the predator with the help of the minimal wave speed of travelling wave solutions connecting the predator-free state to the co-existence state. The travelling waves are found to be non-monotonic in this system. The formation of wave trains has been demonstrated for an extended range of nonlocal interactions. Finally, the importance of psychological effects in shaping the dynamics of nonlocal collective behaviour is demonstrated with several representative examples. 
\end{abstract}


\end{frontmatter}


\section{Introduction} \label{sec:1}

Interaction between two or more people is a necessary component of cooperative activities. Hence, an analysis of social interaction features at the dyadic, group, and population levels is necessary to comprehend the mechanisms supporting cooperation. Numerous academic fields, including psychology, mathematics, and the social sciences, have conducted studies on cooperative behaviour. One area of mathematics where this type of process has been explored extensively is cooperative game theory \cite{wang2023competitive}. When an accumulation of organisms cooperates, they act or work together for their mutual or shared advantage. For example, evolutionary biology commonly uses the ideas of kin selection and reciprocal altruism to explain cooperation. In the meantime, group dynamics and interactions determine a group's level of cooperation in psychology. In bio-social dynamics, collaboration is a diverse phenomenon that is impacted by a variety of social and biological elements. For years, researchers have endeavoured to comprehend the mechanisms and motivations behind cooperation among individuals and their consequences for public policy, economics, social psychology, and other fields \cite{gokcekus2021exploring}.

From the latest epidemic to the financial sector, many different kinds of psychological repercussions may be seen in human society. Referring to a person's emotional, psychological, and social well-being, mental health is a crucial component of overall health \cite{wood2023beyond, houlihan2023emotion, lv2021panic, ambrosio2021beyond, carrero2021mathematical, kaiser2022scientific}. Individuals who endure prolonged periods of sadness, anxiety, stress, or even PTSD may see physiological changes in them, such as elevated cardiac reactivity or decreased heart blood flow. In the context of population dynamics also, which includes interactions between prey and predator, several psychological consequences are now becoming visible. Evaluating the psychological consequences that predation has on the prey population is a crucial component of these models. These consequences may include anxiety, tension, and behavioural adjustments that alter the predator-prey dynamic. A prey species' primary goal is to find refuge and avoid harm from predators, which causes them to handle their feelings of fear, panic, PTSD, and even help each other. On the other hand, the predator species aims to pursue success, growth, collaboration, and defence. Zanette et al. showed that PTSD-like behavioural and brain abnormalities may happen to wild animals, suggesting that PTSD is not abnormal and that long-lasting consequences of fear of predators, which may have an impact on fertility and survival, are common in the natural world \cite{zanette2019predator}. Social influence has a role in the process of collective decision-making, influencing the outcome as a whole. Collaborative decision-making has a crucial role in several domains, including healthcare, social challenges, infrastructure development, environmental protection, and ecological systems. For instance, an interactive discussion is required to make choices on environmental protection, such as preserving endangered species and cutting carbon emissions etc. It also applies to the brain network model in the neuroscientific field \cite{thieu2023social}. Nevertheless, in biological systems, organisms use a variety of tactics to either attack their prey or elude predators; cooperative hunting, herd/schooling behaviour, group defence, and other tactics rely on collective consensus rather than individual judgment. \par

Bio-social dynamics is the study of how social and biological factors interact to shape social behaviour and personality traits \cite{tadic2020modeling, tadic2017mechanisms, tadic2021microscopic, tadic2021self, tadic2023evolving, dankulov2015dynamics, mitrovic2022analysis}. The approaches to studying these interactions in various application contexts range from time-series analysis and self-organized critical dynamics tools to multiscale models. Psychological factors play a major role in this dynamic because they affect how a population perceives, interprets, and responds to its social and biological surroundings. Gokcekus et al. describe how various routes to cooperation can be tested within the social network framework \cite{gokcekus2021exploring}. For instance, cooperative behaviours can be categorized as altruistic or mutualistic depending on whether the behaviour incurs a net cost to the direct fitness of the actor or not. Cooperation can be both local and nonlocal. Local cooperation, for instance, takes place amongst members of a certain group or within a defined geographic region, whereas nonlocal cooperation happens over greater distances or between people who are not geographically adjacent. However, cooperation occurs when people work together to achieve a similar objective in both situations, but the extent and kind of cooperation may differ. So, local or nonlocal cooperation depends on the scope, scale, and range of interaction involved. Social dynamics and ecological considerations are involved in adapting predator-prey models to bio-social systems. These models offer a more thorough comprehension of the intricate interactions seen in ecosystems. predator-prey models can contain social interactions, such as communication, competition, or cooperation, to better understand how these interactions affect the behaviour of both species. However, in order to overcome adverse circumstances, a species may display altruistic behaviours within its community. Furthermore, game-theoretic notions may be added to the predator-prey model to enable the study of evolutionary dynamics and tactics. Additionally, a predator-prey model that incorporates spatial aspects shows how the movements of both predators and prey affect the system. So, as an exemplification of a bio-social system, we can focus on the population dynamics here by considering a predator-prey interaction. 

Individuals of social animals exhibit widespread cooperative activity in various biological systems. When two prerequisites are met, cooperation is permissible. First, there are certain outcomes from a common activity that all participants must partake in, and second, the benefits to the participants must outweigh the costs of the common action. In nature, a variety of cooperative behaviours may be seen, including the creation of groups to decide on actions, cooperative breeding, cooperative hunting, predator inspection, defence against predator attacks, attending to the needs of injured group members, etc. In population dynamics and conservation biology, the study of mutually beneficial interactions (such as cooperation) among social animals is an area of study that is expanding quickly. Cooperative hunting is widespread among vertebrates (e.g., lions \cite{stander1992cooperative}, African wild dogs \cite{creel1995communal}, chimpanzees \cite{boesch1994cooperative, boesch2002cooperative}, wolves \cite{schmidt1997wolf}, ants \cite{dejean2010arboreal}, crocodiles \cite{dinets2015apparent}, primates \cite{kappeler2006cooperation}, etc.) Different animals adopt different strategies for making a group. For instance, group hunting in lionesses generally involves a formation whereby some lionesses circle prey while others wait for prey to move towards them \cite{stander1992cooperative}. A computational model of lion-hyena interaction developed in the work of Rajagopalan et al. to understand the evolution of mobbing behaviours of spotted hyenas while attacking a lion \cite{rajagopalan2019factors}. Jarvey et al., in their experiment, found that group hunting of spotted hyenas increased as tolerance increased and as the relative payoff from cooperative hunting increased, and higher-ranking species prefer to group hunt than lower-ranking species under despotic sharing conditions \cite{jarvey2022effects}. Roy et al., in their work, described that Galapagos sea lions cooperate with each other while going for their prey, Amberstripe scad, which is a schooling, fast-swimming semipelagic fish \cite{de2021cooperation}. Striped marlin, a marine group-hunting predator from the billfish family, uses different sequences while hunting their prey \cite{hansen2023mechanisms, hansen2022mechanisms}. The hunting strategies of 61 species of mammals, birds, vertebrates, and invertebrates were listed by Packer and Ruttan \cite{packer1988evolution}. They included information on which prey or kinds of prey to predate, whether to prey on one or several prey at once and the proportion of their hunting success that they obtained from various sources.
In biological systems, cooperation plays a vital role and is essential to animal social life \cite{jarvey2022effects, umrao2023bifurcation}. Some predators in ecosystems use cooperative hunting techniques to improve their chances of success and terrorize their victims. The dynamics of cooperative hunting in a McCann and Yodzis food network model with three distinct species were assessed by Duarte et al. \cite{duarte2009chaos}. Berec \cite{berec2010impacts} investigated predator-prey foraging facilitation in a modified Rosenzweig-MacArthur model and found that it could disturb predator-prey stability. Using an extended Lotka-Volterra model that considers predator-hunting cooperation, Alves and Hilker \cite{alves2017hunting} investigated how cooperation affects predator-prey dynamics and have concluded that the persistence of the predator population can be improved by cooperative hunting. Pal et al. have worked on a predator-prey model with the incorporation of predator's cooperative hunting where the growth of prey species is affected due to fear of predation \cite{pal2019fear}. Saha and Samanta formulated a predator-prey model in the presence of predator's cooperative hunting and prey's defence mechanism \cite{saha2020prey}, where it is shown that a higher cooperative hunting rate decreases the stable region of population coexistence. There are already some published articles where the significance of cooperative hunting has been explored in predator-prey interactions \cite{umrao2023bifurcation, yousef2022hunting}, whereas there are some other articles dealing with predator-prey models with nonlocal interactions \cite{fetecau2013nonlocal, massaccesi2017nonlocal, amorim2020predator}.

In the context of population biology, travelling wave solutions describe an evolving zone of transition of a particular species from lower density to higher density or vice versa \cite{dunbar1983travelling}. This kind of solution has been seen to exist in a number of ecological groups in the actual world. For example, by selecting parameter values that suit the field data, Dubois \cite{dubois1975model} examined the significance of travelling wave solutions for the oceanic plankton population and showed that waves propagate across their whole habitat. Wyatt described a similar kind of wave-like migration of plankton patches in the Southern Bight \cite{wyatt1973biology}. By reviewing the data available for the years 1972–1983 in 19 forestry districts of Finland and for the years 1900–1935 in 53 areas of France, Ranta and Kaitala \cite{ranta1997travelling} demonstrated the occurrence of travelling waves in the vole population. Lambin et al. detected the existence of travelling waves in a cyclic field vole population in northern Britain \cite{lambin1998spatial}. According to their estimate, the wave travelled at an average speed of 19 km/year from north to east in the direction of $78^{\circ}$. On the other hand, MacKinnon et al. \cite{mackinnon2001scale} examined field data on the vole population collected over 2.5 years from 147 sites in northern England and estimated a travelling wave of vole density that travelled through these sites at an average speed of 14 km/year heading in the direction of $66^{\circ}$ from the north. Similar to this, \cite{bierman2006changes} found a connection between travelling waves and the delayed density dependency structure for the cyclic field vole population. 

Furthermore, travelling waves were shown to exist for larch budmoth outbreaks in central Europe by Bj{\o}rnstad et al. \cite{bjornstad2002waves}. The travelling waves' average speed in the European Alps, measured from the west, is estimated to be 219.8 km/year, travelling from $67^{\circ}$ to $80^{\circ}$. It was covered in \cite{johnson2004landscape} how landscape geometry shapes travelling waves in outbreaks of larch budmoth. According to the authors, inhabitants in low-connected environments have a delayed shift from peaks to troughs. Connectivity-rich landscapes see sharp population drops, and as a result, they become focal points for travelling waves. Additionally, Johnson et al. investigated how landscape mosaic affected recurrent travelling waves in outbreaks of larch budmoth \cite{johnson2006landscape}. They discovered that in highly connected landscapes, over-compensatory density dependency serves as an underlying mechanism to create recurrent travelling waves. In the work of Ranta et al., the variance in population density for hare and lynx populations was explained in terms of travelling waves \cite{ranta1997dynamics}. Furthermore, evidence for a demographic travelling wave in the Red Grouse population in northeast Scotland was shown by Moss et al. \cite{moss2000spatial}.
Numerous researchers have worked for many years to study theoretically whether travelling waves and wave-train solutions exist in interacting population systems. The presence of travelling wave solutions for spatially extended Lotka–Volterra type models of predator-prey interactions was demonstrated by Dunbar \cite{dunbar1983travelling}. Later, he demonstrated the existence of such solutions for the Lotka–Volterra type predator-prey model with a logistic growth of prey by demonstrating a heteroclinic orbit linking two stable states in $\mathbb{R}^{4}$ \cite{dunbar1984traveling}. Moreover, he proved that a diffusive Rosenzweig–MacArthur model has also travelling wave and wave-train solutions \cite{dunbar1986traveling}, and to make the analysis simpler, he assumed the prey species immobile by choosing the diffusion coefficient associated with them as zero. Nevertheless, Huang et al. later examined the identical model in which both species have been considered mobile \cite{huang2003existence}. Gardner established travelling wave solutions of a diffusive predator-prey system using the connection index, a version of the Conley index \cite{gardner1984existence, gardner1991stability}. One may find further research on the presence of travelling wave and wave-train solutions for various diffusive predator-prey system types in \cite{ai2017traveling, bennett2017periodic, safuan2016travelling, tian2018traveling, wu2018traveling, zhao2018note, ducrot2021large}, as well as in other application areas \cite{alfaro2019traveling, wen2017bifurcations}. \par

These days, exploring nonlocal interaction has become an emerging pivotal topic, and studying the role of cooperative hunting has also grown an interest among researchers. While some of the research works have been mentioned above, observing the nonlocal interaction through cooperative hunting has not been analyzed yet. Though we are concerned with exploring the involvement of nonlocal interactions in bio-social dynamics, we have restricted our domain in this work to population dynamics by formulating a predator-prey interaction, and the main intention is to explore the impact of psychological effects affecting the system in terms of cooperative hunting through a local as well as a nonlocal approach. Here, a predator-prey relationship is presented in which the cooperative hunting strategy adopted by the predators affects the prey's growth. Along with that, the carrying capacity of the predator species is chosen to be a prey-dependent function, but the predator is considered to be a specialist one as it is not dependent only on the targeted prey. There is an indirect psychology that works among predators to save themselves from being extinct in this case. The main focus here is to analyze how this cooperation among predators affects the overall dynamics of the system. However, we have not restricted ourselves to the analysis of the local model, but we have also explored how the incorporation of the nonlocal cooperation term affects the overall dynamic behaviour. Not only that, but a detailed analysis of the existence of travelling wave and wave-train solutions for the spatio-temporal and nonlocal models is also investigated. The rest of the paper is organized as follows: Firstly, we have proposed a local predator-prey system with cooperative hunting and analyzed its dynamic nature in Section \ref{sec:2}. In Section \ref{sec:3}, we have incorporated one-dimensional diffusion terms in the system under periodic boundary conditions and analyzed the local stability of the spatio-temporal model. The model is reintroduced with nonlocal interaction in Section \ref{sec:4}, and stability analysis for the corresponding system has been performed. The travelling wave solutions connecting the equilibrium points have also been studied in this section. All the analytical findings have been validated through numerical simulations in Section \ref{sec:5}, and conclusions are given in Section \ref{sec:6}.

\section{Dynamics of bio-social interactions with a simplified model and its analysis} \label{sec:2}

Bio-social models mainly focus on the biological and social dimensions to understand different phenomena of human and non-human beings. On the other hand, understanding cooperative behaviour in bio-social models involves considering a combination of psychological, sociological, and biological factors. In this work, we have shrunken the domain a bit by focusing on the population dynamics and choosing a particular psychological effect from a lot to observe its endowment in the dynamics. In order to adapt predator-prey models to bio-social systems, we have to consider the social dynamics, interaction, cooperation, and even human involvement in the interactions between predators and prey, along with ecological considerations. Through these models, we get a comprehensive idea of intricate interactions in the ecosystem. 

Several researchers extensively study the dynamical complexities of the interacting predator-prey model to understand the long-term behaviour of the species. Most predator-prey models are based upon the classical Lotka-Volterra model, where predators' prey consumption rate is the predator's growth rate with a conversion factor. However, in this paper, we have considered the predator species to be generalist by discarding the situation of the predator's preference towards the targeted prey only. However, it is also taken into consideration that the carrying capacity of the predator species in the environment will be dependent on the mentioned prey, and they grow according to the modified Leslie-Gower model \cite{leslie1960properties}. As the predator becomes cooperative while hunting the prey population, it significantly impacts their growth, which is why we have considered the prey's involvement in the carrying capacity even for the generalist predator. Leslie introduced that the environmental carrying capacity $(K_{V})$ of the predators is proportional to the prey biomass $(U)$ i.e. $K_{V}=p_{1}U$, for some positive constant $p_{1}$, which is known as the conversion factor of prey into predators \cite{leslie1948some}. Also, a positive constant term is added with $K_{V}$ to avoid a mathematical singularity when the prey population becomes zero. From a biological point of view, in the case of severe scarcity of prey, predators can switch over to other populations (alternative food), but their growth will be restricted due to the non-availability of their favourite food. Also, adding this type of positive constant introduces a maximum decrease rate, which stands for environmental protection. Thus, we consider a modified Leslie-Gower model in a homogeneous environment where predator consumes their prey with Holling type-II functional response as follows:

\begin{equation}\label{eq:eq1}\tag{A}
\begin{aligned}
\frac{dU}{dT}&=r_{1}U\left(1-\frac{U}{K_{1}}\right)-\frac{eV}{1+h_{1}eU},\ \ U(0)>0 \\
\frac{dV}{dT}&=V\left(c_{1}-\frac{\gamma_{1}V}{m_{1}+p_{1}U}\right), \ \ v(0)>0
\end{aligned}  
\end{equation} 

In this model, the prey and predator biomass is denoted by $U$ and $V$, respectively, where the prey species grows with an intrinsic growth rate $r_{1}$, and carrying capacity $K_{1}$ in the environment. The parameters $e$ and $h_{1}$ signify the encounter rate of predators with their prey and the predator's handling time of a prey individual. The growth rate of predator population is noted by $c_{1}$, whereas $p_{1}$ is the conversion factor of prey into predator biomass. As the predator looks for a secondary food source in scarcity of their targeted prey, so $m_{1}$ represents a positive constant related to the alternative/additional food. 

In this work, we have mainly provided an exemplification of cooperative behaviour based on predator-prey dynamics, where it is taken into consideration that the predator species cooperate among themselves while attacking and consuming their targeted prey species. This example is instructive since cooperation in a species is a collective decision-making characteristic of a species where mutual understanding plays a key role. It involves psychological effects that lead the whole group to work collectively instead of making solitary decisions. As a result, the developed collective strategy of cooperation ultimately helps the whole cluster in many different functions, e.g., to gain more food. 
Some species show this type of behaviour while going for their prey; for example, wolves and lionesses show cooperative behaviours during hunting \cite{ripple2000historic, stander1992cooperative}. By attacking in this way, the predators could have a significant impact on their prey. The above observations motivate us to study the dynamics of the predator-prey system, including cooperative hunting by predators. There is literature where the predator-prey models with cooperative hunting of predators have been studied, but we have extended our model with nonlocal interaction, which is studied in the later part. The overall aim of the present study is to explore the following issues:
\begin{itemize}
    \item[(a)] The impact of cooperative hunting on stabilizing or destabilizing the dynamics of the system when it is introduced through a local as well as a nonlocal term?
    \item[(b)] How does the travelling wave train connect the coexisting equilibrium with the boundary points in the presence of the predator's cooperation?
\end{itemize}

Before proposing the model with the implementation of predator's cooperative hunting, let us briefly outline Berec’s approach to cooperative behaviours among predators \cite{berec2010impacts}. Berec, in his work, discussed cooperative phenomena through a Holling type-II functional response using two approaches. Firstly, when the encounter rate of predators $(e)$ increases with increasing predator density and handling time $(h)$ is chosen as constant, i.e., 
\begin{align*}
e(V)=\frac{e_{0}}{(a+V)^{w}},\ \ h(V)=h,\ a\geq 0\ \mbox{and}\ w<0.
\end{align*}
Then, the type-II functional response with this $e$ and $h$ is called the encounter-driven functional response. Secondly, if the encounter rate of predators $e$ is constant and the handling time $h$ decreases with the increase of predator biomass, i.e.,
\begin{align*}
e(V)=e, \ \ h(v)=h_{0}(a+V)^{w},\ a\geq 0\ \mbox{and}\ w<0.
\end{align*}
In this case, the type-II functional response with this $e$ and $h$ is called handling-driven functional response. In this work, we have confined ourselves to Berec’s encounter-driven functional response. In particular, we have chosen $w=-1$, and so $e(V)=e_{0}(a+V)=\alpha_{1}+\beta_{1}V$, where $\alpha_{1}$ corresponds to $e_{0}a$ and $\beta_{1}$ corresponds to $e_{0}$. Therefore, model (\ref{eq:eq1}) with the Berec’s encounter-driven functional response boils down as follows:
\begin{equation}\label{eq:eq2}\tag{B}
\begin{aligned}
\frac{dU}{dT}&=r_{1}U\left(1-\frac{U}{K_{1}}\right)-\frac{(\alpha_{1}+\beta_{1}V)UV}{1+h_{1}(\alpha_{1}+\beta_{1}V)U},\ \ U(0)>0 \\
\frac{dV}{dT}&=V\left(c_{1}-\frac{\gamma_{1}V}{m_{1}+p_{1}U}\right), \ \ v(0)>0
\end{aligned}  
\end{equation} 
It should be noted that in the absence of hunting cooperation $(\beta_{1}=0)$, the above model (\ref{eq:eq2}) will be as same as the modified Leslie-Gower model (\ref{eq:eq1}). Now, in order to reduce the number of parameters, we re-scale the variables, and the system (\ref{eq:eq2}) turns into as follows:
\begin{equation}\label{eq:det1}
\begin{aligned}
\frac{du}{dt}&=u(1-u)-\frac{c(1+\alpha v)uv}{m+(1+\alpha v)u},\ \ u(0)>0 \\
\frac{dv}{dt}&=sv\left(1-\frac{\gamma v}{\beta +u}\right), \ \ v(0)>0
\end{aligned}  
\end{equation} 
The new variables and parameters are given as follows: $U=K_{1}u,\ V=K_{1}v/\alpha_{1},\ T=t/r_{1},\ c=1/(r_{1}\alpha_{1}),\ \alpha=(\beta_{1}K_{1})/\alpha_{1}^{2},\ m=1/(h_{1}K_{1}\alpha_{1}),\ s=c_{1}/r_{1},\ \gamma=\gamma_{1}/(p_{1}c_{1}\alpha_{1})$ and $\beta=m_{1}/(p_{1}K_{1})$. The system parameters are chosen to be positive.

\subsection{Positivity and boundedness}

\noindent The positivity and boundedness of the solutions are demonstrated in this part to confirm the well-behavedness of the temporal system (\ref{eq:det1}). The positivity of a system indicates that both species survive as the state variables $u(t)$ and $v(t)$ are the biomass of prey and predator population, respectively. The boundedness of every dynamical system can be seen as the inherent restriction on unbridled growth because of the finite resources available. Before proving the theorems, it should be noted that $\mathbb{R}_{+}^{2}=\{(u,v): u>0,\ v>0\}$.

\begin{theorem} \label{Theorem-2.1}
Solutions of system (\ref{eq:det1}), starting in $\mathbb{R}_{+}^{2}$, are positive and bounded with time.
\end{theorem}
\begin{proof}
Functions on the right-hand side of the system (\ref{eq:det1}) are continuous and locally Lipschitzian (as they are polynomials and rationals in $(u,v)$), so there exists a unique solution $(u(t),v(t))$ of the system with positive initial conditions $(u(0), v(0))>0$ on $[0,\tau],$ where $0<\tau<+\infty$. From the first, and second equation of (\ref{eq:det1}) we have
\begin{align*}
\frac{du}{dt}&=u\left[1-u-\frac{c(1+\alpha v)v}{m+(1+\alpha v)u}\right]=u\psi_{1}(u,v) \\
\Rightarrow u(t)&=u(0)\exp\left[\int^t_0 \psi_{1}(u(z),v(z))\,dz\right]> 0, \ \textrm{for} \ u(0)> 0.
\end{align*}
Similarly,
\begin{align*}
\frac{dv}{dt}&=v\left[s\left(1-\frac{\gamma v}{\beta+u}\right)\right]=v\psi_{2}(u,v)\\
\Rightarrow v(t)&=v(0)\exp\left[\int^t_0 \psi_{2}(u(z),v(z))\,dz\right]> 0, \ \textrm{for} \ v(0)> 0.
\end{align*}
So, the solutions of the system (\ref{eq:det1}) are feasible with time. Now, the first equation of (\ref{eq:det1}) gives:

\begin{equation*}
\begin{aligned}
\frac{du}{dt}&=u(1-u)-\frac{c(1+\alpha v)uv}{m+(1+\alpha v)u}\leq u(1-u) \\
\displaystyle \Rightarrow & \limsup_{t\rightarrow \infty}u(t)\leq 1.
\end{aligned}
\end{equation*}

\noindent Then, from the second equation, we have
\begin{equation*}
\begin{aligned}
\frac{dv}{dt}&= sv\left(1-\frac{\gamma v}{\beta+u}\right) \leq sv\left(1-\frac{\gamma v}{\beta+1}\right) = sv\left(1-\frac{v}{\left(\frac{\beta+1}{\gamma}\right)}\right)\\
\displaystyle \Rightarrow & \limsup_{t\rightarrow \infty}v(t)\leq \left(\frac{\beta+1}{\gamma}\right).
\end{aligned}
 \end{equation*}
So, the solutions of system (\ref{eq:det1}) enter into the region: 
$\displaystyle \mathbb{T} = \left\{(u,v)\in \mathbb{R}^{2}_{+}: 0< u \leq 1; 0< v \leq (\beta+1)/\gamma\right\}$.
 \end{proof}

\subsection{Equilibrium points and local stability analysis of system \eqref{eq:det1}}
\noindent In this section, we are going to determine different kinds of existing equilibrium points of the system (\ref{eq:det1}). In addition, some stability-related parametric requirements will be figured out for these points as well. First, we will analyse the dynamic nature of the boundary equilibrium points, and then we will discuss the interior equilibrium point(s).

\subsubsection{Analysis of the boundary equilibrium points}
From the nullclines, it is obtained that the system has one trivial equilibrium point $E_{0}=(0,0)$, along with one predator-free equilibrium point $E_{1}=(1,0)$ and one prey-free equilibrium point $E_{2}=(0,\widetilde{v}) = (0, \beta/\gamma)$. To derive the local stability conditions for an equilibrium point, either boundary or interior, we need to look for the eigenvalues of the corresponding Jacobian matrix. For the temporal system (\ref{eq:det1}), the Jacobian matrix is as follows:
\mathcenter
\begin{equation}\ \label{eq:det2}
\textbf{J}=  \left(
\begin{array}{cc}
a_{11} & a_{12}  \\
a_{21} & a_{22}
\end{array}
\right),
\end{equation}
where $\displaystyle a_{11}=1-2u-\frac{cm(1+\alpha v)v}{\{m+(1+\alpha v)u\}^{2}},\ a_{12}=-\frac{cmu\alpha v}{\{m+(1+\alpha v)u\}^{2}}-\frac{cu(1+\alpha v)}{\{m+(1+\alpha v)u\}}$, \\
$\displaystyle a_{21}=\frac{s\gamma v^{2}}{(\beta+u)^{2}}$ and $\displaystyle a_{22}=s\left(1-\frac{2\gamma v}{\beta+u}\right)$.

\begin{theorem} \label{Theorem-2.2}
In system (\ref{eq:det1}), $E_{0}$ is an unstable equilibrium point and $E_{1}$ is a saddle point.
\end{theorem}
\begin{proof}
The Jacobian matrices corresponding $E_{0}$ and $E_{1}$ are given as: 
\mathleft
\begin{equation*}
\textbf{J}|_{E_{0}}=\left(
\begin{array}{cc}
1 & 0 \\
0 & s
\end{array}
\right)\ \mbox{and}\ \textbf{J}|_{E_{1}}=\begin{pmatrix}
b_{11} & b_{12} \\
0 & b_{22}
\end{pmatrix}=\begin{pmatrix}
-1 & -\frac{c}{m+1} \\
0 & s
\end{pmatrix}.
\end{equation*}
So, the eigenvalues corresponding to $\textbf{J}|_{E_{0}}$ and $\textbf{J}|_{E_{1}}$ are obtained as $(1, s)$ and $(-1, s)$, respectively. Both the eigenvalues for $\textbf{J}|_{E_{0}}$ are positive which gives $E_{0}$ an unstable equilibrium point. On the other hand, the eigenvalues at $\textbf{J}|_{E_{1}}$ are of opposite signs, concluding that $E_{1}$ is a saddle point.
\end{proof}

\begin{note}
The states where either both prey and predators are absent or only prey persist in the system are not inherently stable. Constant growth of prey and predators, or other influences, could disrupt the equilibrium easily.
\end{note}

\begin{theorem} \label{Theorem-2.3}
$E_{2}$ is locally asymptotically stable (LAS) when $m\gamma^{2}<c\beta(\gamma+\alpha\beta)$ holds.
\end{theorem}
\begin{proof}
The Jacobian matrix of the system (\ref{eq:det1}) at $E_{2}$ is given by
\begin{equation*}
\textbf{J}|_{E_{2}}=\begin{pmatrix}
 b_{11} & 0 \\
b_{21} & b_{22}   
\end{pmatrix}=\begin{pmatrix}
1-\frac{c\beta(\gamma+\alpha\beta)}{m\gamma^{2}} & 0 \\
\frac{s}{\gamma} & -s   
\end{pmatrix},
\end{equation*}
and its eigenvalues are the roots of the equation: 
\mathcenter
\begin{equation*}
\lambda^{2}+C_{1}\lambda+C_{2}=0,
\end{equation*}
where $C_{1}=-(b_{11}+b_{22})$ and $C_{2}=b_{11}b_{22}$. The equation has roots with negative real parts if $C_{1},\ C_{2}>0$, and this occurs when $b_{11}<0$, i.e., $m\gamma^{2}<c\beta(\gamma+\alpha\beta)$.
\end{proof}

\begin{note}
The above condition indicates that the predator's cooperation is beneficial for getting more food, but a higher level of cooperation may lead to prey extinction.
\end{note}

\subsubsection{Analysis of an interior equilibrium point}

The proposed system has an interior equilibrium point $E^{*}=(u^{*},v^{*})$, which exists when both of the following non-trivial prey and predator nullclines hold in the interior of the first quadrant: 
\begin{align*}
 \overline{F}_{1}(u,v)=u(1-u)-\frac{c(1+\alpha v)uv}{m+(1+\alpha v)u}=0,\ \
 \overline{F}_{2}(u,v)=sv\left(1-\frac{\gamma v}{\beta+u}\right)=0.
\end{align*}

From $\overline{F}_{2}(u,v)=0$ we get $v^{*}=(\beta+u^{*})/\gamma$ and using this we obtain $u^{*}$ as the root of the following equation 
\begin{align}\label{eq:stabdet1}
   R(u)\equiv A_{1}u^{3}+A_{2}u^{2}+A_{3}u+A_{4}=0,
\end{align}
where 
$A_{1} =\gamma\alpha, A_{2}=\alpha(c-\gamma)+\gamma(\gamma+\alpha\beta), A_{3} = (c-\gamma)(\gamma+\alpha\beta)+c\alpha\beta+m\gamma^{2}$, and $A_{4}=c\beta(\gamma+\alpha\beta)-m\gamma^{2}$.

\begin{figure}[!htb]
    \centering
    \includegraphics[width=10cm,height=7cm]{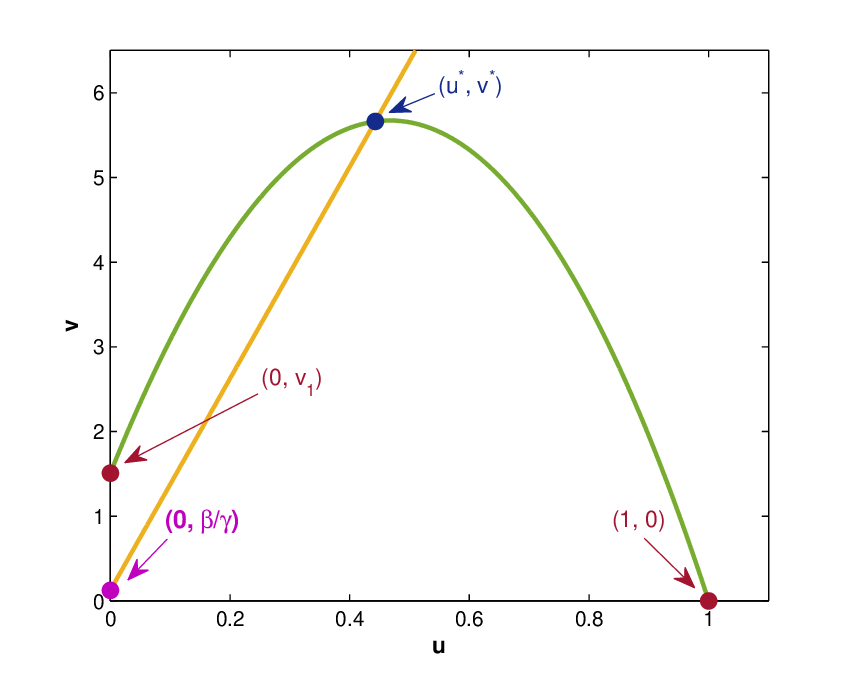}
    \caption{The nontrivial prey are predator nullcline is represented by ({\color{OliveGreen}\solidrule}) and ({\color{BurntOrange}\solidrule}) colored curves when $m\gamma^{2}>c\beta(\gamma+\alpha\beta)$. The positive interior equilibrium point is denoted by $E^{*}=(u^{*}, v^{*})$. The nontrivial prey and predator nullclines intersect the vertical axis at $(0, v_{1})$ and $(0, \widetilde{v})$, respectively with $v_{1}=(-c+\sqrt{c^{2}+4cm\alpha})/2c\alpha$ and $\widetilde{v}=\beta/\gamma$. Parameter values are mentioned in Section \ref{sec:5}. 
} \label{fig:null1}
\end{figure}

 \begin{table}
 \centering
\begin{tabular}{|c||c|c|c|c||c|}
\hline
  & $A_{1}$ & $A_{2}$ & $A_{3}$ & $A_{4}$ & No. of roots \\
 \hline
 \hline
& + & + & + & + & No roots \\
\multirow{-2}{*}{$c>\gamma$} & + & + & + & - & 1 root \\
\hline
& + & + & - & + & at least 0 or at most 2 roots \\
   & + & - & + & + & at least 0 or at most 2 roots \\
   & + & - & - & + & at least 0 or at most 2 roots \\
   & + & + & - & - & 1 root \\
   & + & - & + & - & at least 1 or at most 3 roots \\
 \multirow{-6}{*}{$c<\gamma$} & + & - & - & - & 1 root \\  
\hline
\end{tabular}
\caption{Possibilities of existence of positive interior equilibrium points}
 \end{table}

Here, $v^{*}>0$ only when $u^{*}>0$. Now, Fig. \ref{fig:null1} shows that the nontrivial prey nullcline intersects the vertical axis at $(0,v_{1})$ where $v_{1}=(-c+\sqrt{c^{2}+4cm\alpha})/2c\alpha$, and the nontrivial predator nullcline intersects the vertical axis at $(0,\widetilde{v})$ where $\widetilde{v}=\beta/\gamma$. There are two possibilities of occurrence such as $(0,\widetilde{v})$ lies below or above $(0,v_{1})$. Now, $\widetilde{v}<v_{1}$ implies $m\gamma^{2}>c\beta(\gamma+\alpha\beta)$, i.e., $A_{4}<0$ whereas $\widetilde{v}>v_{1}$ indicates $A_{4}>0$. Theorems \ref{Theorem-2.2} and \ref{Theorem-2.3} have already stated that the fulfilment of the condition $m\gamma^{2}<c\beta(\gamma+\alpha\beta)$ is enough for the instability of all the boundary equilibria, and this happens when $(0, v_{1})$ lies above $(0, \widetilde{v})$ [see Fig. \ref{fig:null1}]. However, the equation (\ref{eq:stabdet1}) has a positive root in $(0, v_{1})$ indicating the existence of a unique positive equilibrium point, and hence, we have made a remark below to state the condition under which the proposed system (\ref{eq:det1}) has a unique coexistence equilibrium point. 

\begin{remark}
If $m\gamma^{2}>c\beta(\gamma+\alpha\beta)$, the system (\ref{eq:det1}) has a unique interior equilibrium point except for the cases $A_{2}<0$ and $A_{3}>0$ along with $R(u_{1})>0$ and $R(u_{2})<0$ where $u_{1}$ and $u_{2}$ are the positive roots of the equation $R'(u)=0$ satisfying $0<u_{1}<u_{2}<1$.
\end{remark}

The above analysis implies that both populations survive in the system when the overall consumption of predators with each other's cooperation fails to suppress the maximum prey biomass. Furthermore, by examining the non-trivial prey and predator nullclines, it is discovered from a geometric perspective that the potential number of coexisting equilibrium points may be one or more. But, for the present data set [see Section \ref{sec:5}], we have obtained only one distinct coexistence steady state in this system.

\begin{theorem} \label{th2.4}
$E^{*}$ is locally asymptotically stable (LAS) when $D_{1},\ D_{2}>0$ hold (mentioned in the proof).
\end{theorem}
\begin{proof}
\mathleft
\begin{equation*}
\textrm{For}\ E^{*}=(u^{*}, v^{*}):\ \
\textbf{J}(E^*)=\textbf{J}|_{E^{*}}=\left(
\begin{array}{cc}
a_{11} & a_{12} \\
a_{21} & a_{22}
\end{array}
\right),
\end{equation*}
where $\displaystyle a_{11}=-u^{*}+\frac{c(1+\alpha v^{*})^{2}u^{*}v^{*}}{\{m+(1+\alpha v^{*})u^{*}\}^{2}},\ a_{12}=-\frac{cmu^{*}\alpha v^{*}}{\{m+(1+\alpha v^{*})u^{*}\}^{2}}-\frac{cu^{*}(1+\alpha v^{*})}{\{m+(1+\alpha v^{*})u^{*}\}}$, \\
$\displaystyle a_{21}=\frac{s\gamma v^{*2}}{(\beta+u^{*})^{2}}=\frac{s}{\gamma}$ and $\displaystyle a_{22}=\frac{-s\gamma v^{*}}{\beta+u^{*}}=-s$.
The characteristic equation corresponding to $\textbf{J}|_{E^{*}}$ is given as follows: 
\mathcenter
\begin{equation}\label{eq:2.4}
\lambda^{2}+D_{1}\lambda+D_{2}=0,
\end{equation}
where $D_{1}=-\mbox{tr}(\textbf{J}(E^{*}))=-(a_{11}+a_{22})$ and $D_{2}=\det(\textbf{J}(E^{*}))=a_{11}a_{22}-a_{12}a_{21}$. So, by Routh- Hurwitz criteria, the equilibrium point will be locally asymptotically stable if the equation has roots with negative real parts, i.e., $D_{1}>0$ and $D_{2}>0$, i.e., 
\begin{equation} \label{eq:2.4_1}
    a_{11}+a_{22}<0\ \textrm{and}\ a_{11}a_{22}>a_{12}a_{21}.
\end{equation}
\end{proof}

\begin{theorem} \label{theorem-2.7}
If $m\gamma^{2}>c\beta(\gamma+\alpha\beta)$ holds, then the interior equilibrium point $E^{*}$ is globally asymptotically stable in the positive quadrant $\mathbb{R}^{2}_{+}$
\end{theorem}

\begin{proof}
For $E^{*}=(u^{*}, v^{*})$: let $H(u,v)=\frac{1}{uv},\ \overline{F}_{1}(u,v)=u(1-u)-\frac{c(1+\alpha v)uv}{m+(1+\alpha v)u}$ \\
and $\overline{F}_{2}(u,v)=sv\left(1-\frac{\gamma v}{\beta+u}\right)$. \\
So, $H(u, v)>0$ in the interior of the positive quadrant of $u$-$v$ plane. Hence we get
\mathleft
\begin{equation*}
\begin{aligned}
\Delta(u, v)=\frac{\partial}{\partial u}(H\overline{F}_{1})+\frac{\partial}{\partial v}(H\overline{F}_{2}) &= -\frac{1}{v}-\frac{c(1+\alpha v)^{2}}{\{m+(1+\alpha v)u\}^{2}}-\frac{s\gamma}{u(\beta+u)}
&<0,\ \ \ \forall (u,v)\in \mathbb{T}
\end{aligned}
\end{equation*}
\mathcenter
So, by the Bendixson-Dulac criterion, no limit cycle exists in the positive quadrant of the $u$-$v$ plane. Also, all the boundary equilibrium points are unstable, which makes $E^{*}$ the only stable attractor of the system. It indicates that the interior of the first quadrant of the $u$-$v$ plane becomes the basin of attraction of $E^{*}$. Henceforth, $E^{*}$ becomes globally asymptotically stable. It completes the proof. 
\end{proof}

\begin{remark}
Let the unique coexisting state of the system exist. The region $\mathbb{T}$ is positively invariant by the semiflow generated by the system (\ref{eq:det1}) and contains all the non-negative equilibrium points of the system. This positive semiflow in $\mathbb{R}^{2}_{+}$ admits a global attractor, denoted by $\mathcal{A}_{\mathbb{R}^{2}_{+}}$ which lies in $\mathbb{T}$.
\end{remark}

Now, let us analyse the description of $\mathcal{A}_{\mathbb{R}^{2}_{+}}$. For this, we first discuss the existence of an interior attractor by considering the regions 
$$\mathbb{B}_{u}^{+}=\{(u,v)\in \mathbb{R}^{2}_{+}: v=0\}\ \ \textrm{and}\ \ \mathbb{B}_{v}^{+}=\{(u,v)\in \mathbb{R}^{2}_{+}: u=0\}$$
and the state-space (disjoint) decomposition
$\mathbb{R}^{2}_{+}=\mbox{int}(\mathbb{R}^{2}_{+})\cup (\mathbb{B}_{u}^{+} \cup \mathbb{B}_{v}^{+})$

In the following lemma, we have used this notion of global attractor that was first used in the literature written by Hale \cite{hale1988asymptotic, hale2000dissipation}. Some results on global attractors can be obtained in the work of Magal and Zhao \cite{magal2005global} and Magal \cite{magal2009perturbation} also.

\begin{lemma} \label{Lem-2.2}
Let the interior point of the system (\ref{eq:det1}) exist. Then the semiflow generated by the system lies in $\mathbb{R}^{2}_{+}$, has a global attractor $\mathcal{A}_{\mathbb{R}^{2}_{+}}$, which is a connected, compact subset and it attracts all compact subsets of $\mathbb{R}^{2}_{+}$.
\end{lemma}

\begin{note}
It is checked that the global attractors in $\mathbb{B}_{u}^{+}$ and $\mathbb{B}_{v}^{+}$ are, respectively
$$\mathcal{A}_{\mathbb{B}_{u}^{+}}=\{(u,v)\in \mathbb{R}^{2}_{+}: u \in [0, 1],\ v=0\}\ \ \textrm{and}\ \ \mathcal{A}_{\mathbb{B}_{v}^{+}}=\left\{(u,v)\in \mathbb{R}^{2}_{+}: u=0,\ v \in [0, \beta/\gamma\right]\}.$$

So, both $\mathcal{A}_{\mathbb{B}_{u}^{+}}$ and $\mathcal{A}_{\mathbb{B}_{v}^{+}}$ contain two equilibria in $\mathbb{B}_{u}^{+}$ and $\mathbb{B}_{v}^{+}$, respectively along with the heteroclinic orbits joining those equilibrium points. 
\end{note}

\noindent\textbf{Proof of Lemma \ref{Lem-2.2}:}
The existence of boundary attractors follows from the invariance of $\mathcal{A}_{\mathbb{B}_{u}^{+}}$ and $\mathcal{A}_{\mathbb{B}_{v}^{+}}$ together with dissipativeness of variable stated in Theorem \ref{Theorem-2.1}. Here we want to show the existence of interior attractor $\mathcal{A}_{\mbox{int}(\mathbb{R}^{2}_{+})}$. For this, it is sufficient to show the uniform persistence of state-space decomposition $\left(\mathbb{B}_{u}^{+} \cup \mathbb{B}_{v}^{+}; \mbox{int}(\mathbb{R}^{2}_{+})\right)$ \cite{hale1989persistence}, i.e., there exists a constant $\Theta>0$ such that for each $\displaystyle (u_{0}, v_{0})\in [0, \infty)^{2}$, 
$$\liminf_{t\rightarrow\infty}\min(u(t), v(t))\geq \Theta$$ 
Now, similar to the work of Hale and Waltman \cite{hale1989persistence}, in this model there are three equilibria $M_{1}=(0,0),\ M_{2}=(1,0)$ and $M_{3}=(0, \beta/\gamma)$ on the boundary $\left(\mathbb{B}_{u}^{+} \cup \mathbb{B}_{v}^{+}\right)$.
Hence, it is sufficient to show the local repulsive dynamics of these boundary equilibrium points with respect to the interior region $\mbox{int}(\mathbb{R}^{2}_{+})$.

\noindent Suppose it does not hold. Then there exists some sufficiently small $0<\epsilon<\epsilon_{1}$ such that 
$$u_{0}>0,\ v_{0}>0,\ u(t)+v(t)<\epsilon,\ \forall t>0,$$
where $\epsilon_{1}=\{-(c+m)+\sqrt{(c+m)^{2}+4c\alpha}\}/2c\alpha$.
The first equation of system (\ref{eq:det1}) gives
\begin{align*}
\frac{du}{dt} &=u(1-u)-\frac{c(1+\alpha v)uv}{m+(1+\alpha v)u} \geq u\left[1-\epsilon-\frac{c\epsilon}{m}(1+\alpha \epsilon)\right]\geq 0
\end{align*}
This gives $u(t)\rightarrow \infty$ as $t\rightarrow \infty$, which is a contradiction because of dissipativeness. \\
Next, let us assume that there exists some small $0<\epsilon<(1+\beta)/(1+\gamma)$ such that 
$$|u(t)-1|+v(t)<\epsilon,\ \forall t\geq 0$$
From the second equation of system (\ref{eq:det1}) we have
\begin{align*}
    \frac{dv}{dt}&= sv\left(1-\frac{\gamma v}{\beta+u}\right) \geq sv\left[1-\frac{\gamma \epsilon}{(\beta+1)-\epsilon}\right]\geq 0.
\end{align*}
So, we can deduce that $v(t)\rightarrow \infty$ as $t\rightarrow \infty$, which is a contradiction. \\
Lastly, let us assume that there exists some small $0<\epsilon<\epsilon_{2}$ such that 
$$u(t)+\Big|v(t)-\frac{\beta}{\gamma}\Big|<\epsilon,\ \forall t\geq 0,$$
where $\epsilon_{2}=[-\{2c\alpha\beta+\gamma(c+m)\}+\sqrt{-\{2c\alpha\beta+\gamma(c+m)\}^{2}+4c\alpha\{m\gamma^{2}-c\beta(\gamma+\alpha\beta)\}}]/2c\alpha\gamma>0$.
From the second equation of system (\ref{eq:det1}) we have
\begin{align*}
   \frac{du}{dt} &=u(1-u)-\frac{c(1+\alpha v)uv}{m+(1+\alpha v)u} \geq u\left[1-\epsilon-\frac{c}{m}\left(\frac{\beta}{\gamma}+\epsilon\right)\left\{1+\alpha\left(\frac{\beta}{\gamma}+\epsilon\right)\right\}\right]\geq 0
\end{align*}
which is a contradiction as $u(t)\rightarrow \infty$ as $t\rightarrow \infty$ in this case. This completes the proof.

\begin{corollary}
The result holds for all those trajectories that start from either $\mathbb{B}_{u}^{+}$, or $\mathbb{B}_{v}^{+}$, or $int(\mathbb{R}^{2}_{+})$. It means the semiflow lies in $\mathbb{B}_{u}^{+}$, and $\mathbb{B}_{v}^{+}$, has respective compact and connected global attractor $\mathcal{A}_{\mathbb{B}_{u}^{+}}$ and $\mathcal{A}_{\mathbb{B}_{v}^{+}}$ which attracts the compact subsets of $\mathbb{B}_{u}^{+}$ and $\mathbb{B}_{v}^{+}$, respectively. On the other hand, if a system-generated semiflow lies in  $int(\mathbb{R}^{2}_{+})$, it has a global attractor $\mathcal{A}_{int(\mathbb{R}^{2}_{+})}$ which attracts the compact subsets of $int(\mathbb{R}^{2}_{+})$.
\end{corollary}



\subsection{Local bifurcations around the equilibrium points of system (\ref{eq:det1})}
The local bifurcations around the equilibrium points are analyzed mainly with the help of Sotomayor's theorem and Hopf's bifurcation theorem. In the system, if the stability condition of any of the equilibrium points violates in such a way that the corresponding determinant becomes $0$, giving a simple zero eigenvalue, then there will occur transcritical bifurcation, and we can observe the exchange of stability in that bifurcation threshold. The following theorem will state the condition where such bifurcation can be observed in $E_{1}$.

Let $\textbf{V}= (v_{1},v_{2})^{T}$ and $\textbf{W}= (w_{1},w_{2})^{T}$, respectively be the eigenvectors of $\textbf{J}|_{\left(eq.\ point\right)}$ and $\textbf{J}|^{T}_{\left(eq.\ point\right)}$ for a zero eigenvalue at the equilibrium point.\\
Let $\overline{\textbf{F}}=(\overline{F}_{1},\overline{F}_{2})^{T},$ where 
\mathleft
\begin{align*}
 \overline{F}_{1}=u(1-u)-\frac{c(1+\alpha v)uv}{m+(1+\alpha v)u},\ \
 \overline{F}_{2}=sv\left(1-\frac{\gamma v}{\beta+u}\right)
\end{align*}

\begin{theorem}
System (\ref{eq:det1}) undergoes a transcritical bifurcation around $(0,\widetilde{v})$ at $\displaystyle c_{[TC]}=m\gamma^{2}/\beta(\gamma+\alpha\beta)$, choosing $c$ as the bifurcating parameter.
\end{theorem}
\begin{proof}
  \begin{equation*}
\textbf{J}|_{E_{2}}=\begin{pmatrix}
 b_{11} & 0 \\
b_{21} & b_{22}  
\end{pmatrix}=\begin{pmatrix}
1-\frac{c\beta(\gamma+\alpha\beta)}{m\gamma^{2}} & 0 \\
\frac{s}{\gamma} & -s 
\end{pmatrix}.
\end{equation*}
The eigenvalues are the roots of the equation $\lambda^{2}-(b_{11}+b_{22})\lambda+b_{11}b_{22}=0$, which gives the roots with negative real parts when $b_{11}<0$. Let $c_{[TC]}$ be the value of $c$ such that $m\gamma^{2}=c\beta(\gamma+\alpha\beta)$ so that $\textbf{J}|_{E_{2}}$ has a simple zero eigenvalue at $c_{[TC]}$. So, at $c=c_{[TC]}:$
\begin{equation*}
\textbf{J}|_{E_{2}}=\begin{pmatrix}
0 & 0 \\
b_{21} & b_{22}  
\end{pmatrix}.
\end{equation*}
The calculations give $\textbf{V}=(v_{1},v_{2})^{T}=(s,s/\gamma)^{T}$ and $\textbf{W}=(1,0)^{T}.$ Therefore,
\mathleft
\begin{equation*}
\begin{aligned}
\Omega_{1}&= \textbf{W}^{T}.\overline{\textbf{F}}_{c}(E_{2}, c_{[TC]})=\frac{-(1+\alpha v)uv}{[m+(1+\alpha v)u]}\bigg|_{E_{2}}=0, \\
\Omega_{2}&= \textbf{W}^{T}\left[D\overline{\textbf{F}}_{c}(E_{2}, c_{[TC]})\textbf{V}\right]=\frac{-s\beta(\gamma+\alpha\beta)}{m\gamma^{2}}\neq0 \\
\textrm{and}\ \Omega_{3}&= \textbf{W}^{T}\left[D^{2}\overline{\textbf{F}}(E_{2},c_{[TC]})(\textbf{V},\textbf{V})\right] =\frac{2}{m\gamma}(1-m\gamma)-\frac{2\gamma s^{2}}{\beta(\gamma+\alpha\beta)}(\gamma+2\alpha\beta) \neq0. \\
\end{aligned}
\end{equation*}
\mathcenter
Then, by Sotomayor's Theorem, the system undergoes a transcritical bifurcation around $E_{2}$ at $c=c_{[TC]}.$
\end{proof}

\begin{theorem} \label{Theorem-8.2}
System (\ref{eq:det1}) undergoes a saddle-node bifurcation for the bifurcation parameter $\alpha$ around $(u^{*},v^{*})$ when $R(u)=0$ (mentioned in equation (\ref{eq:stabdet1})) has a double root.
\end{theorem}
\begin{proof}
Let $u_{[sn]}$ be a double root of $R(u)=0$ such that $u_{[sn]}<1$. So, $R(u_{[sn]})=0=R^{'}(u_{[sn]})$ and $R^{''}(u_{[sn]})\neq 0$. Let $\alpha_{[sn]}$ be the threshold value of $\alpha$ such that for $\alpha=\alpha_{[sn]},\ u_{[sn]}$ is a double root of $R(u)=0$. The non-trivial nullclines touch each other at $(u_{[sn]}, \alpha_{[sn]})\equiv E^{*}_{[sn]}.$ Let, $\overline{F}_{1}(u,v)=ug_{1}(u,v)$ and $\overline{F}_{2}(u,v)=vg_{2}(u,v)$. Here $\frac{dv^{(g_{1})}}{du}$ denotes the slope of $g_{1}(u,v)=0$ and $\frac{dv^{(g_{2})}}{du}$ denotes the slope of $g_{2}(u,v)=0$. \\
Also, $\frac{dv^{(g_{2})}}{du}\bigg|_{E^{*}_{[sn]}}=\frac{dv^{(g_{1})}}{du}\bigg|_{E^{*}_{[sn]}}$ where, $\frac{dv^{(g_{2})}}{du}\bigg|_{E^{*}_{[sn]}}=-\frac{\partial g_{2}/\partial u}{\partial g_{2}/\partial v}\bigg|_{E^{*}_{[sn]}}$ and $\frac{dv^{(g_{1})}}{du}\bigg|_{E^{*}_{[sn]}}=-\frac{\partial g_{1}/\partial u}{\partial g_{1}/\partial v}\bigg|_{E^{*}_{[sn]}}$.
\begin{equation*}
\begin{aligned}
\textbf{J}|_{E^{*}}=\begin{pmatrix}
\frac{\partial\overline{F}_{1}}{\partial u} & \frac{\partial\overline{F}_{1}}{\partial v} \\
\frac{\partial\overline{F}_{2}}{\partial u} & \frac{\partial\overline{F}_{2}}{\partial v}
\end{pmatrix}=
\begin{pmatrix}
u\frac{dg_{1}}{du} & u\frac{dg_{1}}{dv} \\
v\frac{dg_{2}}{du} & v\frac{dg_{2}}{dv} \\
\end{pmatrix}&=
\begin{pmatrix}
a_{11} & a_{12} \\
a_{21} & a_{22}
\end{pmatrix}\\
&=
\begin{pmatrix}
-u^{*}+\frac{c(1+\alpha v^{*})^{2}u^{*}v^{*}}{\{m+(1+\alpha v^{*})u^{*}\}^{2}} & -\frac{cmu^{*}\alpha v^{*}}{\{m+(1+\alpha v^{*})u^{*}\}^{2}}-\frac{cu^{*}(1+\alpha v^{*})}{\{m+(1+\alpha v^{*})u^{*}\}} \\
\frac{s\gamma v^{*2}}{(\beta+u^{*})^{2}} & \frac{-s\gamma v^{*}}{\beta+u^{*}},
\end{pmatrix},
\end{aligned}
\end{equation*}
So, $Det\left(\textbf{J}|_{E^{*}_{[sn]}}\right)=uv\left(\frac{dg_{1}}{du}\frac{dg_{2}}{dv}-\frac{dg_{2}}{du}\frac{dg_{1}}{dv}\right)\bigg|_{E^{*}_{[sn]}}=0$ and the Jacobian matrix $\textbf{J}|_{E^{*}_{[sn]}}$ has a simple zero eigenvalue. The calculations give the right and left eigenvectors at $E^{*}_{[sn]}$ as $\textbf{V}=(-a_{12},a_{11})^{T}\ \textrm{and}\ \textbf{W}=(-a_{21},a_{11})^{T}$. So,
\mathleft
\begin{equation*}
\begin{aligned}
\Omega_{1}&= \textbf{W}^{T}.\overline{\textbf{F}}_{\alpha}(E^{*}, \alpha_{[sn]}) =\frac{cmuv^{2}a_{21}}{[m+(1+\alpha v)u]^{2}}\bigg|_{E^{*}_{[sn]}}\neq 0,
\end{aligned}
\end{equation*}
\begin{equation*}
\begin{aligned}
\textrm{and}\ \ \Omega_{2}&= \textbf{W}^{T}\left[D^{2}\overline{\textbf{F}}(E^{*}, \alpha_{[sn]})(\textbf{V},\textbf{V})\right]\\
&=a_{21}\left[2a_{12}^{2}-\frac{2cm[v^{*}(1+\alpha v^{*})^{2}a_{12}^{2}+\{(m+u^{*})(1+\alpha v^{*})+m\alpha v^{*}\}a_{11}a_{22}-\alpha u^{*}(m+u^{*})a_{11}^{2}]}{[m+(1+\alpha v^{*})u^{*}]^{3}}\right]\\
&-\frac{2s\gamma a_{11}}{(\beta+u^{*})^{3}}\left[a_{12}^{2}v^{*2}+(\beta+u^{*})v^{*}a_{11}a_{12}+(\beta+u^{*})^{2}a_{11}^{2}\right] \neq0.
\end{aligned}
\end{equation*}
\mathcenter
Thus, from Sotomayor's Theorem, system (\ref{eq:det1}) undergoes a non-degenerate saddle-node bifurcation around $E^{*}$ at $\alpha=\alpha_{[sn]}.$
\end{proof}

If any of the mentioned inequalities in (\ref{eq:2.4_1}) is violated, then the equilibrium point becomes unstable, and the system performs oscillatory or non-oscillatory behaviour. In fact, the system starts to oscillate around $(u^{*}, v^{*})$ if $D_{1}>0$ along with $D_{1}^{2}-4D_{2}<0$ as the eigenvalues will be in the form of the complex conjugate in this case. So, we get the following theorem:

\begin{theorem}
If $E^{*}$ exists with the feasibility conditions, then a simple Hopf bifurcation occurs at unique $\alpha=\alpha_{[H]}$, where $\alpha_{[H]}$ is the positive root of $D_{1}(\alpha)=0$ providing $D_{2}(\alpha_{[H]})>0$ (stated in equation (\ref{eq:2.4})).
\end{theorem}
\begin{proof}
At $\alpha=\alpha_{[H]}$, the characteristic equation of system (\ref{eq:det1}) at $E^{*}$ is $(\lambda^{2}+D_{2})=0$ and so, the equation has a pair of purely imaginary roots $\lambda_{1}=i\sqrt{D_{2}}$ and $\lambda_{2}=-i\sqrt{D_{2}}$ where $D_{2}(\alpha)$ is a continuous function of $\alpha$. \\
In the small neighbourhood of $\alpha_{[H]},$ the roots are $\lambda_{1}=p_{1}(\alpha)+ip_{2}(\alpha)$ and $\lambda_{2}=p_{1}(\alpha)-ip_{2}(\alpha)\ (p_{1},\ p_{2}\in \mathbb{R}$). \\
To show the transversality condition, we check 
$$\displaystyle \left(\frac{d}{d\alpha}[Re(\lambda_{i}(\alpha))]\right)\bigg|_{\alpha=\alpha_{[H]}}\neq 0,\ \mbox{for}\ i=1,2.$$ \\
Put $\lambda(\alpha)=p_{1}(\alpha)+ip_{2}(\alpha)$ in (\ref{eq:2.4}), we get
\mathcenter
\begin{equation} \label{eq:2.6}
(p_{1}+ip_{2})^{2}+D_{1}(p_{1}+ip_{2})+D_{2}=0.
\end{equation}
Differentiating with respect to $\alpha$, we get
\begin{equation*}
2(p_{1}+ip_{2})(\dot{p_{1}}+i\dot{{p_{2}}})+D_{1}(\dot{p_{1}}+i\dot{p_{2}})+\dot{D_{1}}(p_{1}+ip_{2})+\dot{D_{2}}=0.
\end{equation*}
Comparing the real and imaginary parts from both sides, we have \\
\begin{equation} \label{eq:2.7}
(2p_{1}+D_{1})\dot{p_{1}}-(2p_{2})\dot{p_{2}}+(\dot{D_{1}}p_{1}+\dot{D_{2}}) = 0, \\
\end{equation}
\begin{equation} \label{eq:2.8}
(2p_{2})\dot{p_{1}}+(2p_{1}+D_{1})\dot{p_{2}}+\dot{D_{1}}p_{2}=0.
\end{equation}
Solving we get, $\displaystyle \dot{p_{1}}=\frac{-2p_{2}^{2}\dot{D_{1}}-(2p_{1}+D_{1})(\dot{D_{1}}p_{1}+\dot{D_{2}})}{(2p_{1}+D_{1})^{2}+4p_{2}^{2}}$. \\
At, $p_{1}=0,\ p_{2}=\pm \sqrt{D_{2}}:\ \displaystyle \dot{p_{1}}=\frac{-2\dot{D_{1}}D_{2}-D_{1}\dot{D_{2}}}{D_{1}^{2}+4D_{2}}\neq 0$. Hence, this completes the proof.
\end{proof}

\begin{note}
System (\ref{eq:det1}) admits a unique stable periodic orbit surrounding the interior equilibrium point $E^{*}$ and the system has no other periodic orbit if $c(1+\alpha v^{*})^{2}u^{*}v^{*}>(s+u^{*})\{m+(1+\alpha v^{*})u^{*}\}^{2}$ holds.
\end{note}

Now, in the following theorem, we have proved that there are heteroclinic orbits that join the boundary equilibrium points $E_{1}$ and $E_{2}$ with the stable coexistence state $E^{*}$. The result is numerically validated also in the later part, and it helps us to show that the travelling waves are connecting $E^{*}$ with the boundary points instead of trivial equilibrium $E_{0}$.

\begin{theorem}\label{theorem-2.9}
Let $m\gamma^{2}>c\beta(\gamma+\alpha\beta)$ hold such that a unique interior equilibrium point of the system (\ref{eq:det1}) exists. Then the system has a unique heteroclinic orbit $(u, v)$ joining $(1,0)$ to the boundary of the interior attractor $\mathcal{A}_{\mbox{int}(\mathbb{R}^{2}_{+})}$. Also, there will be another unique heteroclinic orbit joining $(0, \beta/\gamma)$ to the boundary of the interior attractor $\mathcal{A}_{\mbox{int}(\mathbb{R}^{2}_{+})}$. The global attractor $\mathcal{A}_{\mathbb{R}^{2}_{+}}$ is composed of four disjoint parts 
 $$\mathcal{A}_{\mathbb{R}^{2}_{+}}=[0,1]\times\{0\} \bigcup \{0\}\times [0, \beta/\gamma] \bigcup \{(u(t),v(t)), t\in \mathbb{R}\} \bigcup \mathcal{A}_{\mbox{int}(\mathbb{R}^{2}_{+})}$$
\end{theorem}

\begin{proof}
The whole proof of the theorem contains four parts. First, we need to derive the existence of heteroclinic orbits with the help of the connectedness argument for the global attractor. Next, we prove the heteroclinic orbits start from the stationary points $E_{1}$ and $E_{2}$. Lastly, we conclude the uniqueness of each of the heteroclinic orbits by center unstable manifold argument \cite{ducrot2021large, ducrot2013multiple}.  
\begin{itemize}
    \item[(a)] Connectedness: The largest global attractor $\mathcal{A}_{\mathbb{R}^{2}_{+}}$ is
connected since it attracts the convex subset $\mathbb{T}$. It follows that the
projection of $\mathcal{A}_{\mathbb{R}^{2}_{+}}$ on the horizontal and vertical axis is a compact interval. The global attractor $\mathcal{A}_{\mathbb{R}^{2}_{+}}$ contains the interior global attractor $\mathcal{A}_{\mbox{int}(\mathbb{R}^{2}_{+})}$ which is compact, connected and locally stable. The global attractor $\mathcal{A}_{\mathbb{R}^{2}_{+}}$ also contains the boundary attractors $\mathcal{A}_{\mathbb{B}_{u}^{+}}$ and $\mathcal{A}_{\mathbb{B}_{v}^{+}}$.
The connectedness of $\mathcal{A}_{\mathbb{R}^{2}_{+}}$ and compactness of $\mathcal{A}_{\mbox{int}(\mathbb{R}^{2}_{+})},\ \mathcal{A}_{\mathbb{B}_{u}^{+}}$ and $\mathcal{A}_{\mathbb{B}_{v}^{+}}$ imply 
$$\mathcal{A}_{\mathbb{R}^{2}_{+}}-\left(\mathcal{A}_{\mbox{int}(\mathbb{R}^{2}_{+})} \bigcup \mathcal{A}_{\mathbb{B}_{u}^{+}} \bigcup \mathcal{A}_{\mathbb{B}_{v}^{+}}\right) \neq \phi$$

Moreover, by using Theorem 3.2 from the work conducted by Hale and Waltman \cite{hale1989persistence}, we deduce that for each point $(u,v) \in \mathcal{A}_{\mathbb{R}^{2}_{+}}-\left(\mathcal{A}_{\mbox{int}(\mathbb{R}^{2}_{+})} \bigcup \mathcal{A}_{\mathbb{B}_{u}^{+}} \bigcup \mathcal{A}_{\mathbb{B}_{v}^{+}}\right)$ the alpha and limit sets satisfy the following 
$$\alpha(u,v) \in \mathcal{A}_{\mathbb{B}_{u}^{+}} \bigcup \mathcal{A}_{\mathbb{B}_{v}^{+}} \ \textrm{and} \ \omega(u,v) \in \mathcal{A}_{\mbox{int}(\mathbb{R}^{2}_{+})}$$

Finally, since the boundary attractor has a Morse decomposition
$M_{1}=\{(0, 0)\},\ M_{2}=\{(1, 0)\}$ and $M_{3}=\{(0, \beta/\gamma)\}$, we have either $$\alpha(u,v)=M_{1},\ \textrm{or}\ \alpha(u,v)=M_{2},\ \textrm{or}\ \alpha(u,v)=M_{3},\ \forall (u,v)\in \mathcal{A}_{\mbox{int}(\mathbb{R}^{2}_{+})}-\left(\mathcal{A}_{\mathbb{R}^{2}_{+}} \bigcup \mathcal{A}_{\mathbb{B}_{u}^{+}}\bigcup \mathcal{A}_{\mathbb{B}_{v}^{+}}\right)$$. 

\item[(b)] Non-existence of heteroclinic orbit starting from $E_{0}$: 
Let us assume a heteroclinic orbit exists that starts from $(0,0)$. The second equation of system (\ref{eq:det1}) gives 
\begin{align*}
    \frac{dv}{dt}&=sv\left(1-\frac{\gamma v}{\beta+u}\right), \ \forall t,t_{0}\in \mathbb{R} \\
 \textrm{which gives}\ \  v(t)&=sv(t_{0})\exp\left\{\int_{t_{0}}^{t}\left(1-\frac{\gamma v(s)}{\beta+u(s)}\right)ds\right\}
\end{align*}
Since, $u(t_{0})>0,\ v(t_{0})>0$ and there exists $T<0$ such that $u(t)+v(t)<\epsilon,\ \forall t<T<0$. Then we have 
$$v(t)>sv(t_{0})\exp\left\{\int_{t_{0}}^{t}\left(1-\frac{\gamma\epsilon}{\beta+\epsilon}\right)ds\right\}$$
For $0<\epsilon<\beta/(\gamma-1),\ v(t)\rightarrow\infty$ as $t\rightarrow\infty$, which is a contradiction of boundedness of the global attractor. 

\item[(c)] Existence and uniqueness of heteroclinic orbit starting from $E_{1}$: 
We need to prove the uniqueness only. From $\displaystyle \textbf{J}|_{E_{1}}$ we have the eigenvalues are: $\lambda_{1}=-1<0$ and $\lambda_{2}=s>0$, and the corresponding eigenspaces are 
\begin{align*}
E_{\lambda_{1}}&=\{(u,v)\in \mathbb{R}^{2}: v=0\}\ \textrm{and}\ E_{\lambda_{2}}=\left\{(u,v)\in \mathbb{R}^{2}: u-1=-\frac{c}{(1+s)(1+m)}v\right\}.
\end{align*}
Here, $\mathbb{R}^{2}=E_{\lambda_{1}}\oplus E_{\lambda_{2}}$. The center-unstable manifold at $E_{1}$ is one dimensional. Let $\psi_{cu}:E_{\lambda_{2}}\rightarrow E_{\lambda_{1}}$ be a $C^{1}$ center-unstable manifold and consider the one dimensional manifold defined by 
$$M_{cu}:=\{\chi_{cu}+\psi_{cu}(\chi_{cu}): \chi\in E_{\lambda_{2}}\}$$
The manifold is invariant under the semiflow generated by the system at $(1,0)$. Also, $D_{\chi_{cu}}\psi_{cu}(0)=0$. It means the manifold $M_{cu}$ is tangent to $E_{\lambda_{2}}$ at $(1,0)$. \\
Moreover, we know that there exists $\epsilon>0$ such that $M_{cu}$ contains all negative orbits of the semiflow generated by the system lying in $\displaystyle B_{\mathbb{R}^{2}}((1,0),\epsilon)$ for all $t<0$.

Suppose, there exists two heteroclinic orbits $H_{1}$ and $H_{2}$ starting from $E_{1}$ to $\mathcal{A}_{\mbox{int}(\mathbb{R}^{2}_{+})}$ such that 
$$H_{1}=(u_{1}(t), v_{1}(t))_{t\in \mathbb{R}}\subset \mbox{int}(\mathbb{R}^{2}_{+}),\  \textrm{and}\ H_{2}=(u_{2}(t), v_{2}(t))_{t\in \mathbb{R}}\subset \mbox{int}(\mathbb{R}^{2}_{+})$$
Now, $\displaystyle \lim_{t\rightarrow -\infty}(u_{1}(t), v_{1}(t))=(1,0)$ and $\displaystyle \lim_{t\rightarrow -\infty}(u_{2}(t), v_{2}(t))=(1,0)$. Without loss of generality, one may assume that $(u_{1}(t), v_{1}(t))_{t\leq 0}\subset B_{\mathbb{R}^{2}}((1,0),\epsilon)$ and $(u_{2}(t), v_{2}(t))_{t\leq 0}\subset B_{\mathbb{R}^{2}}((1,0),\epsilon)$. It means that $(u_{1}(t), v_{1}(t))_{t\leq 0}\subset M_{cu}$ and $(u_{2}(t), v_{2}(t))_{t\leq 0}\subset M_{cu}$. 

Let, $\Pi_{\lambda_{1}}$ and $\Pi_{\lambda_{2}}$ be the linear projectors from $\mathbb{R}^{2}$ to $E_{\lambda_{1}}$ and $E_{\lambda_{2}}$, respectively. Then there will be $t_{1}<0$ and $t_{2}<0$ such that 
\begin{align*}
\Pi_{\lambda_{2}}(u_{1}(t), v_{1}(t))= \Pi_{\lambda_{2}}(u_{2}(t), v_{2}(t))
    &\Rightarrow \psi_{cu}(\Pi_{\lambda_{2}}(u_{1}(t), v_{1}(t)))= \psi_{cu}(\Pi_{\lambda_{2}}(u_{2}(t), v_{2}(t))) \\
    &\Rightarrow (u_{1}(t), v_{1}(t)) = (u_{2}(t), v_{2}(t))
\end{align*}
Also, from the uniqueness of solution of system (\ref{eq:det1}), we get
$(u_{1}(t_{1}+.), v_{1}(t_{1}+.))=(u_{2}(t_{2}+.), v_{2}(t_{2}+.))$ which implies $H_{1}=H_{2}$. Hence, the heteroclinic orbit starting from $E_{1}$ is unique.

\item[(d)] Existence and uniqueness of heteroclinic orbit starting from $E_{2}$: 
From $\displaystyle \textbf{J}|_{E_{2}}$ we obtain the eigenvalues as $\lambda_{1}=-s<0$ and $\lambda_{2}=1-c\beta(\gamma+\alpha\beta)/m\gamma^{2}>0$, and the corresponding eigenspaces are 
\begin{align*}
E_{\lambda_{1}}&=\{(u,v)\in \mathbb{R}^{2}: u=0\}\ \textrm{and}\ E_{\lambda_{2}}=\left\{(u,v)\in \mathbb{R}^{2}: v-\frac{\beta}{\gamma}=\frac{s}{\gamma\left(s+1-\frac{c\beta(\gamma+\alpha\beta)}{m\gamma^{2}}\right)}u\right\}.
\end{align*}
It should be noted that, $\mathbb{R}^{2}=E_{\lambda_{1}}\oplus E_{\lambda_{2}}$. The center-unstable manifold at $E_{2}$ is one dimensional. Let $\omega_{cu}:E_{\lambda_{2}}\rightarrow E_{\lambda_{1}}$ be a $C^{1}$ center-unstable manifold and consider the one dimensional manifold defined by 
$$N_{cu}:=\{\kappa_{cu}+\omega_{cu}(\kappa_{cu}): \kappa\in E_{\lambda_{2}}\}$$
The manifold is invariant under the semiflow generated by the system at $E_{2}$. Also, $D_{\kappa_{cu}}\omega_{cu}(0)=0$, which implies the manifold $N_{cu}$ is tangent to $E_{\lambda_{2}}$ at $(0,\beta/\gamma)$. \\
Moreover, there exists $\epsilon>0$ such that $N_{cu}$ contains all negative orbits of the semiflow generated by the system lying in $\displaystyle B_{\mathbb{R}^{2}}((0,\beta/\gamma),\epsilon)$ for all $t<0$.

Suppose, there exists two heteroclinic orbits $J_{1}$ and $J_{2}$ starting from $E_{2}$ to $\mathcal{A}_{\mbox{int}(\mathbb{R}^{2}_{+})}$ such that 
$$J_{1}=(u_{1}(t), v_{1}(t))_{t\in \mathbb{R}}\subset \mbox{int}(\mathbb{R}^{2}_{+}),\  \textrm{and}\ J_{2}=(u_{2}(t), v_{2}(t))_{t\in \mathbb{R}}\subset \mbox{int}(\mathbb{R}^{2}_{+})$$
Now, $\displaystyle \lim_{t\rightarrow -\infty}(u_{1}(t), v_{1}(t))=(0,\beta/\gamma)$ and $\displaystyle \lim_{t\rightarrow -\infty}(u_{2}(t), v_{2}(t))=(0,\beta/\gamma)$. Without loss of generality, one may assume that $(u_{1}(t), v_{1}(t))_{t\leq 0}\subset B_{\mathbb{R}^{2}}((0,\beta/\gamma),\epsilon)$ and $(u_{2}(t), v_{2}(t))_{t\leq 0}\subset B_{\mathbb{R}^{2}}((0,\beta/\gamma),\epsilon)$. It means that $(u_{1}(t), v_{1}(t))_{t\leq 0}\subset N_{cu}$ and $(u_{2}(t), v_{2}(t))_{t\leq 0}\subset N_{cu}$. 

Let, $\Lambda_{\lambda_{1}}$ and $\Lambda_{\lambda_{2}}$ be the linear projectors from $\mathbb{R}^{2}$ to $E_{\lambda_{1}}$ and $E_{\lambda_{2}}$, respectively. Then there will be $t_{1}<0$ and $t_{2}<0$ such that 
\begin{align*}
\Lambda_{\lambda_{2}}(u_{1}(t), v_{1}(t))= \Lambda_{\lambda_{2}}(u_{2}(t), v_{2}(t))
    &\Rightarrow \omega_{cu}(\Lambda_{\lambda_{2}}(u_{1}(t), v_{1}(t)))= \omega_{cu}(\Lambda_{\lambda_{2}}(u_{2}(t), v_{2}(t))) \\
    &\Rightarrow (u_{1}(t), v_{1}(t)) = (u_{2}(t), v_{2}(t))
\end{align*}
Using the uniqueness property of solution of system (\ref{eq:det1}), we get
$(u_{1}(t_{1}+.), v_{1}(t_{1}+.))=(u_{2}(t_{2}+.), v_{2}(t_{2}+.))$ which implies $J_{1}=J_{2}$. Henceforth, the heteroclinic orbit, starting from $E_{2}$, is unique. 
\end{itemize}
\end{proof}

Till now, we have analyzed the temporal model (\ref{eq:det1}), the stability of its equilibrium points, and the existence of heteroclinic orbits joining the equilibria. The results show that the predator-free state will always be a saddle point, whereas the stability of the prey-free state will be dependent on parametric restrictions. The main focus here is to show the contribution of $\alpha$ to control the system dynamics. Though the deterministic model has shown some interesting results on how the cooperative behaviour of predators has a significant impact, now we will extend the work by considering the spatio-temporal model in the following section.  \par

Spatio-temporal reaction-diffusion models in bio-social dynamics are mainly used to get an insight into ecological patterns and interactions, population dynamics, infectious illness transmission, etc. Temporal models focus on the spatial variations of the species only; while on the other hand, spatio-temporal systems combine the spatial as well as the spatio-temporal structures, which help to look at the space-time variation by finding movement patterns that continue over time as well as spatial units. On the other hand, the nonlocal models allow long-range interactions among a species in the system. In this work, we are concerned with the population dynamics by looking into predator-prey interaction. Here, we have emphasized the nonlocal interaction of the predator species due to hunting cooperation, which is a psychological phenomenon affecting the predator's nature in order to achieve higher success at the time of hunting their prey. In Section \ref{sec:3} and \ref{sec:4}, the main aim is to explore how a system will change its dynamics when the species start to move in a direction, and the predators cooperate with each other situated at neighbouring positions.

\section{Inclusion of population movement through spatio-temporal model} \label{sec:3}
The distributions of populations, being heterogeneous, depend not only on time but also on the spatial positions in the habitat. So, it is natural and more precise to study the corresponding PDE problem.  In this work, we have included the spatial aspects into the temporal model (\ref{eq:det1}) and extended it to a coupled reaction-diffusion equations
over a bounded one-dimensional spatial domain $\Omega=[-L, L]\subset \mathbb{R}$ with closed boundary $\partial\Omega$ and $\overline{\Omega}=\Omega\cup \partial\Omega$.
\mathcenter
\begin{align} \label{eq:diff1}
\frac{\partial u}{\partial t}&=d_{1}\frac{\partial^{2}u}{\partial x^{2}}+u(1-u)-\frac{c(1+\alpha v)uv}{m+(1+\alpha v)u} \nonumber \\
\frac{\partial v}{\partial t}&=d_{2}\frac{\partial^{2}v}{\partial x^{2}}+sv\left(1-\frac{\gamma v}{\beta+u}\right)
\end{align}
subject to non-negative initial conditions, and we have chosen periodic boundary conditions. Here, the parameters $d_{1}$ and $d_{2}$ represent the diffusion coefficients of prey and predator species. For the periodic boundary condition, $u(-L)=u(L)$ and $v(-L)=v(L)$ hold at the boundaries.

\subsection{Linear stability analysis of model (\ref{eq:diff1}) around $(u^{*},v^{*})$}

We intend to find the condition for Turing instability. If the homogeneous steady state of the temporal model is locally stable to infinitesimal perturbation but becomes unstable in the presence of diffusion, a scenario of Turing instability occurs. For the temporal model, $(u^{*},v^{*})$ is locally asymptotically stable when $D_{1}=-(a_{11}+a_{22})>0$ and $D_{2}=a_{11}a_{22}-a_{12}a_{21}$ hold. Here, we apply heterogeneous perturbation around $E^{*}$ to obtain the criterion for instability of the spatio-temporal model. Let us perturb the homogeneous steady state of the local system (\ref{eq:diff1}) around $(u^{*},v^{*})$ by 
\mathcenter
\begin{align*}
    \begin{pmatrix}
        u \\
        v
    \end{pmatrix}=\begin{pmatrix}
        u^{*} \\
        v^{*}
    \end{pmatrix}+\epsilon \begin{pmatrix}
        u_{1} \\
        v_{1}
    \end{pmatrix}e^{\lambda t+ikx}
\end{align*}
with $|\epsilon|<<1$ where $\lambda$ is the growth rate of perturbation and $k$ denotes the wave number. Substituting these values in system (\ref{eq:diff1}), the linearization takes the form:
\begin{equation}\label{eq:3.2}
\textbf{J}_{k}\begin{bmatrix}
u_{1}\\
v{1}
\end{bmatrix}\equiv \begin{bmatrix}
 a_{11}-d_{1}k^{2}-\lambda & a_{12} \\
 a_{21} & a_{22}-d_{2}k^{2}-\lambda 
\end{bmatrix}
\begin{bmatrix}
    u_{1} \\
    v_{1}
\end{bmatrix}
=\begin{bmatrix}
        0 \\
        0
    \end{bmatrix}.
\end{equation}
where $a_{11},\ a_{12},\ a_{21}$ and $a_{22}$ are mentioned in the proof of Theorem \ref{th2.4}. We are interested in finding the non-trivial solution of the system (\ref{eq:3.2}), so $\lambda$ must be a zero of $\det(\textbf{J}_{k})=0$, where $\textbf{J}_{k}$ is the coefficient matrix of the system (\ref{eq:3.2}). Now 
\begin{align*}
    \det(\textbf{J}_{k})=\lambda^{2}-B(k^{2})\lambda+C(k^{2})
\end{align*}
with $B(k^{2})=\mbox{tr}(\textbf{J}(E^*))-(d_{1}+d_{2})k^{2},\ C(k^{2})=\det(\textbf{J}(E^*))-(d_{2}a_{11}+d_{1}a_{22})k^{2}+d_{1}d_{2}k^{4}$. So, $\det(\textbf{J}_{k})=0$ gives 
\begin{align*}
    \lambda_{\pm}(k^{2})=\frac{B(k^{2})\pm\sqrt{(B(k^{2}))^{2}-4C(k^{2})}}{2}
\end{align*}
The local asymptotic stability conditions of $(u^{*}, v^{*})$ in the temporal model give $\mbox{tr}(\textbf{J}(E^*))>0$ and $\det(\textbf{J}(E^*))>0$. 
So, for $k=0$, we have $Re(\lambda_{\pm}(0))<0$ due to the stability condition of $(u^{*}, v^{*})$. If $Re(\lambda(k))>0$ for some $k\neq 0$, then the homogeneous steady state $(u^{*}, v^{*})$ becomes unstable for heterogeneous perturbation. In this case, we must have either $B(k^{2})>0$ or $C(k^{2})<0$. But, $B(k^{2})<0$ for all $k$ when the temporal model is locally asymptotically stable. So, the homogeneous solution will be stable under heterogeneous perturbation when $C(k^{2})>0$ for all $k$. If the inequality is violated for some $k\neq 0$, the system is unstable.

Here $\displaystyle k^{2}_{min}=(d_{2}a_{11}+d_{1}a_{22})/2d_{1}d_{2}$ is the minimum value of $k^{2}$ for which $\displaystyle C(k^{2})$ will attain its minimum value, say 
$$\displaystyle C(k^{2})_{min}=(a_{11}a_{22}-a_{12}a_{21})-\frac{(d_{2}a_{11}+d_{1}a_{22})^{2}}{4d_{1}d_{2}}$$. 
This $\sqrt{k^{2}_{min}}$ is known as the critical wave number for Turing instability. And, the critical diffusion coefficient (Turing bifurcation threshold) $d_{1c}$ such that $C(k^{2})_{min}=0$ is given as 
\begin{equation} \label{eq:3.3}
d_{1c}=\frac{d_{2}(a_{11}a_{22}-2a_{12}a_{21})-\sqrt{d_{2}^{2}(a_{11}a_{22}-2a_{12}a_{21})^{2}-d_{2}^{2}a_{11}^{2}a_{22}^{2}}}{a_{22}^{2}}.  
\end{equation}

Moreover, to ensure the positivity of $k^{2}=k^{2}_{\min}$ at the Turing bifurcation threshold, we need to have $d_{2}a_{11}+d_{1}a_{22}>0$, i.e., $d_{1}<d_{2}$ needs to be satisfied for the Turing instability conditions. Hence, the self-diffusion coefficient of the prey population is less than that of the predator population for the model (\ref{eq:diff1}). 

The wavelength at the Turing bifurcation threshold is $\lambda_{m}=2\pi/k_{m}$ where $k_{m}$ is the wavenumber corresponding to the maximum real part of the positive eigenvalue. If the above necessary condition holds and $min_{k^2}<0$ with certain $k^{2}$ in the interval of $(\zeta^{-}, \zeta^{+})$ where 
\begin{equation} \label{eq:3.4_N}
\begin{aligned}
\zeta^{+}(d_{1})&=\frac{(d_{2}a_{11}+d_{1}a_{22})+\sqrt{(d_{2}a_{11}+d_{1}a_{22})^{2}-4d_{1}d_{2}\det(\textbf{J}(E^*))}}{2d_{1}d_{2}} \\
\zeta^{-}(d_{1})& =\frac{(d_{2}a_{11}+d_{1}a_{22})-\sqrt{(d_{2}a_{11}+d_{1}a_{22})^{2}-4d_{1}d_{2}\det(\textbf{J}(E^*))}}{2d_{1}d_{2}}.
\end{aligned}
\end{equation}
Then $(u^{*},v^{*})$ is an unstable homogeneous steady-state of system (\ref{eq:diff1}). Summarizing the conditions, we get the following theorem:

\begin{theorem} \label{Theorem-3.1}
Considering that the interior equilibrium point $(u^{*},v^{*})$ is locally asymptotically stable, if the following conditions hold
\begin{equation} \label{eq:3.5_N}
(d_{2}a_{11}+d_{1}a_{22})^{2}>4d_{1}d_{2}\det(\textbf{J}(E^*)),\ \  d_{1}<d_{2}
\end{equation}
and there exists a wave-number $k^{2}\in(\zeta^{-},\zeta^{+})$, then the positive constant steady-state solution $(u^{*},v^{*})$ of model (\ref{eq:diff1}) is Turing unstable.
\end{theorem}

The formation of stationary and non-stationary spatial patterns where the high prey density coexists with high predator density is observed when $d_{1}<d_{1c}$, but the coexisting equilibrium $(u^{*},v^{*})$ of the local model (\ref{eq:diff1}) remains stable under random heterogeneous perturbation when $d_{1}>d_{1c}$ and the cooperation rate does not exceed the Hopf threshold..


\section{Incorporation of nonlocal interaction} \label{sec:4}

In bio-social dynamics, nonlocal modelling indicates the incorporation of spatial and temporal interactions occurring at different scales and distances in biological and social systems. The nonlocal modelling here is very much related when we try to explore the interactions between biological processes (such as disease transmission or physiological responses) and social dynamics (such as behavioural response of population or societal structures). So, in this work, dealing with a predator-prey interaction with the predator's cooperative hunting strategy uplifted us to explore the significance of nonlocal interaction in the model. \par

Now, the system (\ref{eq:diff1}) uses the assumption that the predator located at the space point $x$ consumes the prey at the same point. However, in reality, a predator located at space point $x$ can help and cooperate with those predators who are located in some areas around this spatial point, which can be obtained by convolution term
$$U(x,t)=(\phi_{\delta}*v)(x,t)=\int_{-\infty}^{\infty}\phi_{\delta}(x-y)v(y,t)dy.$$
Here, the kernel function $\phi_{\delta}(y)$ describes the cooperation of the predator at the space point $x$ with the predator located at the space point $y$. Hence, the kernel $\phi_{\delta}$ is a function dependent on the position $x$. The first subscript $\delta$ is the range of nonlocal interaction. We assume the kernel function $\phi_{\delta}$ to be non-negative with compact support. Also, to preserve the same homogeneous steady-state solutions for both the local and nonlocal models, we assume that $\int_{-\infty}^{\infty}\phi_{\delta}(y)dy=1$. Now, predators' cooperation is limited by their biomass. We can apply this limitation to each space point independently. This means that predators located at space point $y$ cooperate with the predators located at space point $x_{1}$ independently of their concentration at another point $x_{2}$. Under this assumption, we obtain the rate of impact of fear on prey at the space point $x$ as
$$M(x,t)=\frac{c\{1+\alpha(\phi_{\delta}*v)\}uv}{m+\{1+\alpha(\phi_{\delta}*v)\}u}=\frac{c\{1+\alpha\int_{-\infty}^{\infty}\phi_{\delta}(x-y)v(y,t)dy\}uv}{m+\{1+\alpha\int_{-\infty}^{\infty}\phi_{\delta}(x-y)v(y,t)dy\}u}.$$
Implementing the nonlocal cooperative hunting term of predator species as well as the random motion of the population, we get the integro-differential equation model as
\begin{align} \label{eq:loc1}
\frac{\partial u}{\partial t}&=d_{1}\frac{\partial^{2}u}{\partial x^{2}}+u(1-u)-\frac{c\{1+\alpha(\phi_{\delta}*v)\}uv}{m+\{1+\alpha(\phi_{\delta}*v)\}u} \nonumber \\
\frac{\partial v}{\partial t}&=d_{2}\frac{\partial^{2}v}{\partial x^{2}}+sv\left(1-\frac{\gamma v}{\beta+u}\right)
\end{align}
with non-negative initial conditions and periodic boundary conditions.

\subsection{Local stability analysis of nonlocal model (\ref{eq:loc1})}

Both the local model (\ref{eq:diff1}) and the nonlocal model (\ref{eq:loc1}) show the same dynamics for homogeneous steady states. Let us consider $E^{*}=(u^{*},v^{*})$ as the coexisting homogeneous steady-state. Now, perturbing the system around $(u^{*},v^{*})$ by 
\begin{align*}
    \begin{pmatrix}
        u \\
        v
    \end{pmatrix}=\begin{pmatrix}
        u^{*} \\
        v^{*}
    \end{pmatrix}+\epsilon \begin{pmatrix}
        u_{1} \\
        v_{1}
    \end{pmatrix}e^{\lambda t+ikx}
\end{align*}
with $|\epsilon|<<1$ and substituting these values in system (\ref{eq:loc1}) the linearization takes the form:
\begin{equation} \label{eq:4.2}
\overline{\textbf{J}}_{k}\begin{bmatrix}
u_{1}\\
v{1}
\end{bmatrix}\equiv \begin{bmatrix}
 a_{11}-d_{1}k^{2}-\lambda & A_{12} \\
 a_{21} & a_{22}-d_{2}k^{2}-\lambda 
\end{bmatrix}
\begin{bmatrix}
    u_{1} \\
    v_{1}
\end{bmatrix}
=\begin{bmatrix}
        0 \\
        0
    \end{bmatrix}.
\end{equation}
where $a_{11},\ a_{12},\ a_{21}$ and $a_{22}$ are mentioned in the proof of Theorem \ref{th2.4} and 
$$\displaystyle A_{12}=-\left\{\frac{c(1+\alpha v^{*})u^{*}}{m+(1+\alpha v^{*})u^{*}}-\frac{cm\alpha u^{*}v^{*}}{\{m+(1+\alpha v^{*})u^{*}\}^{2}}\frac{\sin k\delta}{k\delta}\right\}=a_{12}+\frac{cm\alpha u^{*}v^{*}}{\{m+(1+\alpha v^{*})u^{*}\}^{2}}\left\{1-\left(\frac{\sin k\delta}{k\delta}\right)\right\}.$$
Now, we are interested in finding the non-trivial solution of the system (\ref{eq:4.2}), so $\lambda$ must be a zero of $\det(\overline{\textbf{J}}_{k})=0$, where $\overline{\textbf{J}}_{k}$ is the coefficient matrix of the system (\ref{eq:4.2}). Now $\det(\overline{\textbf{J}}_{k})=0$ where 
\begin{align*}
\det(\overline{\textbf{J}}_{k})=\lambda^{2}-\Gamma(k,d_{1},\delta)\lambda+\Delta(k,d_{1},\delta),
\end{align*}
where $\Gamma(k,d_{1},\delta)=(a_{11}+a_{22})-k^{2}(d_{1}+d_{2})$ and $\Delta(k,d_{1},\delta)=d_{1}d_{2}k^{4}-(a_{11}d_{2}+a_{22}d_{1})k^{2}+(a_{11}a_{22}-a_{21}A_{12})$ with 
So, $\det(\overline{\textbf{J}}_{k})=0$ will give
\begin{align*}
\lambda(k^{2})=\frac{\Gamma(k,d_{1},\delta)\pm\sqrt{(\Gamma(k,d_{1},\delta))^{2}-4\Delta(k,d_{1},\delta)}}{2}.
\end{align*}
The homogeneous steady-state $(u^{*},v^{*})$ is stable if $\Gamma(k,d_{1},\delta)<0$ and $\Delta(k,d_{1},\delta)>0$ holds for all $k$ for some fixed $d_{1}$ and $\delta$. Now, if the local model is stable, then we already have $\Gamma(k,d_{1},\delta)<0$ here. So, the instability of the coexisting homogeneous steady-state depends on the sign of $\Delta(k,d_{1},\delta)$ only. Moreover, Turing instability occurs if $\Gamma(k,d_{1},\delta)<0$ holds for all $k$ and $\Delta(k,d_{1},\delta)=0$ holds for a unique $k$. 
$\Gamma(k,d_{1},\delta)$ and $\Delta(k,d_{1},\delta)$ depending on the parameter $\delta$, and hence it plays an important role for the above instabilities. Therefore, we find these instability conditions by fixing $\delta$. We assume that the equilibrium point $E^{*}$ is locally asymptotically stable for the temporal model.

To find the condition of Turing instability of the nonlocal model, we need to find a value of $d_{1}$ such that $\Delta(k,d_{1},\delta)=0$ for a fixed $k$ and $\Gamma(k,d_{1},\delta)<0$ for all k. Also, for all $d_{1}$, we have got $\Delta(k,d_{1},\delta)>0$ when $k$ is sufficiently small as well as a large quantity. So, $\Delta(k,d_{1},\delta)=0$ holds for a unique $k$ when 
$$\Delta(k,d_{1},\delta)=0\ \ \textrm{and}\ \ \frac{\partial \Delta(k,d_{1},\delta)}{\partial k}=0$$ hold, i.e., 
\begin{subequations}\label{eq:4.3}
\begin{align} 
d_{1}&=\frac{d_{2}a_{11}k^{2}-(a_{11}a_{22}-a_{12}a_{21})+\frac{a_{21}cm\alpha u^{*}v^{*}}{\{m+(1+\alpha v^{*})u^{*}\}^{2}}\left(1-\frac{\sin k\delta}{k\delta}\right)}{k^{2}(d_{2}k^{2}-a_{22})} \label{eq:4.3a} \\
\mbox{and}\ 4d_{1}d_{2}k^{3}&-2(d_{2}a_{11}+d_{1}a_{22})k+\frac{a_{21}cm\alpha u^{*}v^{*}}{k\delta\{m+(1+\alpha v^{*})u^{*}\}^{2}}\left(\delta \cos k\delta-\frac{\delta \sin k\delta}{k\delta}\right)=0. \label{eq:4.3b}
\end{align}
\end{subequations}
From (\ref{eq:4.3}), eliminating $d_{1}$ leads to the following transcendental equation
\mathleft
\begin{equation} \label{eq:4.4}
\begin{aligned}
2a_{11}d_{2}^{2}k^{4}-2(2d_{2}k^{2}-a_{22})\left[(a_{11}a_{22}-a_{12}a_{21})-\frac{a_{21}cm\alpha u^{*}v^{*}}{\{m+(1+\alpha v^{*})u^{*}\}^{2}}\left\{1-\left(\frac{\sin k\delta}{k\delta}\right)\right\}\right]+\\
\frac{(d_{2}k^{2}-a_{22})a_{21}cm\alpha u^{*}v^{*}}{\{m+(1+\alpha v^{*})u^{*}\}^{2}}\left(\cos k\delta-\frac{\sin k\delta}{k\delta}\right)=0,
\end{aligned}
\end{equation}
which needs to be solved numerically for a fixed value of $\delta$ to obtain the critical wave number $k_{\min}^{T}$. Here, we may get multiple solutions of (\ref{eq:4.4}) for a large value of $\delta$. Out of these multiple values of $k$, we choose $k_{\min}^{T}$ for which $\Delta(k,d_{1},\delta)=0$ holds for a unique $k$. Substitution of this value in (\ref{eq:4.3a}) will give the critical diffusion coefficient $d_{1c}^{T}$. Here $d_{1}<d_{1c}^{T}$ leads to $\Delta(k,d_{1},\delta)<0$, so, Turing instability occurs for $d_{1}<d_{1c}^{T}$.

\subsection{Travelling wave front connecting the coexisting equilibrium $(u^{*}, v^{*})$}
Here, we consider the homogeneous steady-state oscillatory in time. Depending on $\delta$, either the most unstable wave number $k_{c}=0$, or the most unstable wave number $k_{c}>0$. 

Let us consider the case $k_{c}=0$. We try to obtain the solution of a linear system corresponding to the nonlocal model around $(u^{*},v^{*})$ of the form: 
$$\overline{u}=e^{\widehat{\lambda}t}e^{-Kx}\widetilde{u}, \qquad \overline{v}=e^{\widehat{\lambda}t}e^{-Kx}\widetilde{v}$$
Substituting $\overline{u}$ and $\overline{v}$ in linear system of (\ref{eq:loc1}) we get 
\mathcenter
\begin{align*}
\widehat{\lambda}\widetilde{u} &= d_{1}K^{2}\widetilde{u}+a_{11}\widetilde{u}+\overline{a}_{12}\widetilde{v} \\
\widehat{\lambda}\widetilde{v} &= d_{2}K^{2}\widetilde{v}+a_{21}\widetilde{u}+a_{22}\widetilde{v}
\end{align*}
where $\displaystyle a_{11}=-u^{*}+\frac{c(1+\alpha v^{*})^{2}u^{*}v^{*}}{\{m+(1+\alpha v^{*})u^{*}\}^{2}},\ \overline{a}_{12}=-\frac{c(1+\alpha v^{*})u^{*}}{m+(1+\alpha v^{*})u^{*}}-\frac{cm\alpha u^{*}v^{*}}{\{m+(1+\alpha v^{*})u^{*}\}^{2}}\left(\frac{\sinh{k\delta}}{k\delta}\right),\ a_{21}=\frac{s}{\gamma}$ and $\displaystyle a_{22}=-s$.

The non-trivial solution of the above system will be obtained from the following equation:

\[\widehat{\lambda}^{2}-\widehat{\Gamma}\widehat{\lambda}+\widehat{\Delta}=0\]

where $\displaystyle\widehat{\Gamma}=(a_{11}+a_{22})+(d_{1}+d_{2})K^{2}$ and $\widehat{\Delta}=d_{1}d_{2}K^{4}+(d_{2}a_{11}+d_{1}a_{22})K^{2}+(a_{11}a_{22}-a_{12}a_{21})$. Here, $\widehat\lambda=\widehat{\lambda_{r}}+i\widehat{\lambda_{i}}$ is a complex number. So, let us write $e^{\widehat\lambda t}e^{-Kx}=e^{i\widehat{\lambda_{i}}t}e^{-K\left[x-\frac{\widehat{\lambda_{r}}}{K}t\right]}$. The wave speed, here, is given by 

\[\overline{c}(K)=\frac{\widehat{\lambda_{r}}}{K}=\frac{\widehat{\Gamma}}{2K}=\frac{1}{2K}\left[(d_{1}+d_{2})K^{2}-u^{*}+\frac{c(1+\alpha v^{*})^{2}u^{*}v^{*}}{\{m+(1+\alpha v^{*})u^{*}\}^{2}}-s\right]\]

For given $\delta$, we find the local minimum and local maximum of $\overline{c}(K)$ when it exists. 
\mathleft
\begin{align*}
\mbox{Now,}\ \frac{d\overline{c}(K)}{dK}=\left(\frac{d_{1}+d_{2}}{2}\right)-\frac{1}{2K^{2}}\mbox{tr}(\textbf{J}(E^{*}))\ \ \mbox{and}\  \frac{d^{2}\overline{c}(K)}{dK^{2}}=\frac{2 \mbox{tr}(\textbf{J}(E^{*}))}{3\,! K^{3}}
\end{align*}

At a local maximum or minimum of $\displaystyle\overline{c}(K):\ \frac{d\overline{c}(K)}{dK}=0$. Solving the equation, we get a critical value $\displaystyle K_{\overline{c}}$, where $\displaystyle K_{\overline{c}}^{2}=\mbox{tr}(\textbf{J}(E^{*}))/(d_{1}+d_{2})$. At $K=K_{\overline{c}}$, we check the sign of $\displaystyle\frac{d^{2}\overline{c}(K)}{dK^{2}}$. 

So, when the interior state is an unstable focus, we have $\mbox{tr}(\textbf{J}(E^{*}))>0$, which leads to $\displaystyle \frac{d^{2}\overline{c}(K)}{dK^{2}}>0$. In this case, we will obtain a $\overline{c}(K)|_{\min}$. On the other hand, if the interior equilibrium turns out to be a stale point when $k_{c}=0$, there will be no such minimum value of $\overline{c}(K)$. Also, $\overline{c}(K)$ does not have a local minimum for $k_{c}>0$. So, if we fail to obtain any minimum value of $\overline{c}(K)$ by the above process, we use the complex number technique mentioned below. 

We try to obtain the solution of the linear system of (\ref{eq:loc1}) of form: $$\overline{u}=e^{\widetilde{\lambda}t}e^{-Kx+iqx}\widetilde{u}, \qquad \overline{v}=e^{\widetilde{\lambda}t}e^{-Kx+iqx}\widetilde{v}, \ \ \ K>0,$$
where $q$ is the wave number of spatial distribution. Let us consider $Q=q+iK$. Then, substituting $\overline{u}$ and $\overline{v}$ into the linear system of (\ref{eq:loc1}) we obtain the characteristics equation as 
$\widetilde{\lambda}^{2}-\widetilde{\Gamma}\widetilde{\lambda}+\widetilde{\Delta}=0,$
where $\displaystyle\widetilde{\Gamma}=(a_{11}+a_{22})-(d_{1}+d_{2})Q^{2}$ and $\widetilde{\Delta}=d_{1}d_{2}Q^{4}+(d_{2}a_{11}+d_{1}a_{22})Q^{2}+(a_{11}a_{22}-\widetilde{a}_{12}a_{21})$. Here, $\displaystyle \widetilde{a}_{12}=-\frac{c(1+\alpha v^{*})u^{*}}{m+(1+\alpha v^{*})u^{*}}-\frac{cm\alpha u^{*}v^{*}}{\{m+(1+\alpha v^{*})u^{*}\}^{2}}\left(\frac{\sin{Q\delta}}{Q\delta}\right)=a_{12}+\frac{cm\alpha u^{*}v^{*}}{\{m+(1+\alpha v^{*})u^{*}\}^{2}}\left(1-\frac{\sin{Q\delta}}{Q\delta}\right)$ and rest of the coefficients are same as mentioned above. So, $\widetilde{\Gamma}$ and $\widetilde{\Delta}$ both are complex as $Q=(q+iK)$. Here we get the real and imaginary parts of $\widetilde{\Gamma}$ and $\widetilde{\Delta}$ are given by: 
\mathleft
\begin{equation*}
\begin{aligned}
\widetilde{\Gamma_{r}}&=(a_{11}+a_{22})-(d_{1}+d_{2})(q^{2}-K^{2}),\ \ \widetilde{\Gamma_{i}}=-(d_{1}+d_{2})2qK, \\
\widetilde{\Delta_{r}}&=d_{1}d_{2}(q^{4}-6q^{2}K^{2}+K^{4})-(d_{2}a_{11}+d_{1}a_{22})(q^{2}-K^{2})+(a_{11}a_{22}-a_{12}a_{21})- \\
&\frac{a_{21}cm\alpha u^{*}v^{*}}{\{m+(1+\alpha v^{*})u^{*}\}^{2}}\left[1-\frac{\delta(q\sin{q\delta}\cosh{K\delta}+K\cos{q\delta}\sinh{K\delta})}{(q\delta)^{2}+(K\delta)^{2}}\right], \\
\widetilde{\Delta_{i}}&=\frac{a_{21}cm\alpha u^{*}v^{*}}{\{m+(1+\alpha v^{*})u^{*}\}^{2}}\left(\frac{\delta(q\cos{q\delta}\sinh{
K\delta}-K\sin{q\delta}\cosh{K\delta})}{(q\delta)^{2}+(K\delta)^{2}}\right)+d_{1}d_{2}(4qK)(q^{2}-K^{2})\\
&-2qK(d_{2}a_{11}+d_{1}a_{22})
\end{aligned}
\end{equation*}

Now, solving the characteristic equation, we have  
\mathcenter
\begin{align*}
\widetilde{\lambda}=\frac{\widetilde{\Gamma}\pm \sqrt{\widetilde{\Gamma}^{2}-4\widetilde{\Delta}}}{2}=\frac{(\widetilde{\Gamma_{r}}+i\widetilde{\Gamma_{i}})\pm \sqrt{\{(\widetilde{\Gamma_{r}}^{2}-\widetilde{\Gamma_{i}}^{2})-4\widetilde{\Delta_{r}}\} +i(2\widetilde{\Gamma_{r}}\widetilde{\Gamma_{i}}-4\widetilde{\Delta_{i}})}}{2}
\end{align*}
\mathleft
It gives the real part of $\widetilde{\lambda}$ as $\widetilde{\lambda_{r}}=Re(\widetilde{\lambda})$, where
$$\widetilde{\lambda_{r}}=\frac{\widetilde{\Gamma_{r}}\pm \sqrt{\widetilde{G}}}{2},\ \ 
\widetilde{G}=\frac{(\widetilde{\Gamma_{r}}^{2}-\widetilde{\Gamma_{i}}^{2}-4\widetilde{\Delta_{r}})+\sqrt{(\widetilde{\Gamma_{r}}^{2}-\widetilde{\Gamma_{i}}^{2}-4\widetilde{\Delta_{r}})^{2}+4(\widetilde{\Gamma_{r}}\widetilde{\Gamma_{i}}-2\widetilde{\Delta_{i}})^{2}}}{2}.$$

\noindent Also, $e^{\widetilde{\lambda}t}e^{-Qx}=e^{i(\widetilde{\lambda_{i}}t+qx)}e^{-K\left(x-\frac{\widetilde{\lambda_{r}}}{K}t\right)}$, and so, the wave propagation speed is $\displaystyle\overline{c}(K,q)=\widetilde{\lambda_{r}}/K$.

The wave speed depends on $q$ and $K$. We find saddle points of $\overline{c}$, i.e., $\displaystyle \frac{\partial \overline{c}}{\partial K}=0$ and $\displaystyle \frac{\partial \overline{c}}{\partial q}=0$ subject to $\displaystyle \frac{\partial^{2}\overline{c}}{\partial K^{2}}>0$ and $\displaystyle \frac{\partial^{2}\overline{c}}{\partial q^{2}}<0$. As $\widetilde{\lambda_{r}}$ is in transcendental form, we need to verify these conditions numerically and also check the positivity of wave speed. \par

\begin{note}
The wave propagational speed, for the local model (\ref{eq:diff1}), can be obtained from $\overline{c}(q, K)$ when $\delta\rightarrow 0$. Let the wave speed be $c_{l}(q, K)$ in this case, which is obtained by applying the complex number technique (mentioned above). Then we have
$$\displaystyle c_{l}(K,q)=\frac{\widetilde{\Gamma_{r}}\pm \sqrt{\widetilde{G}_{l}}}{2K},$$

where, $\displaystyle \widetilde{G}_{l}=\frac{(\widetilde{\Gamma_{r}}^{2}-\widetilde{\Gamma_{i}}^{2}-4\widetilde{\Delta_{rl}})+\sqrt{(\widetilde{\Gamma_{r}}^{2}-\widetilde{\Gamma_{i}}^{2}-4\widetilde{\Delta_{rl}})^{2}+4(\widetilde{\Gamma_{r}}\widetilde{\Gamma_{i}}-2\widetilde{\Delta_{il}})^{2}}}{2}$ along with 
\mathleft
\begin{equation*}
\begin{aligned}
\widetilde{\Gamma_{r}}&=(a_{11}+a_{22})-(d_{1}+d_{2})(q^{2}-K^{2}),\ \ \widetilde{\Gamma_{i}}=-(d_{1}+d_{2})2qK, \\
\widetilde{\Delta_{rl}}&=d_{1}d_{2}(q^{4}-6q^{2}K^{2}+K^{4})-(d_{2}a_{11}+d_{1}a_{22})(q^{2}-K^{2})+(a_{11}a_{22}-a_{12}a_{21}), \\
\widetilde{\Delta_{il}}&=d_{1}d_{2}(4qK)(q^{2}-K^{2})-2qK(d_{2}a_{11}+d_{1}a_{22})
\end{aligned}
\end{equation*}

\end{note}

Zhang and Jin, in their work, have already shown the existence of predator-invasion travelling wave solutions and prey-spread travelling wave solutions in the upstream and downstream directions \cite{zhang2017traveling}. Using those results, we can state the following theorems:

\begin{theorem} \label{Theorem-4.1}
Let the unique interior point $E^{*}$ exist. Then, system (\ref{eq:diff1}) has a positive travelling wave solution $(u(x+c_{1}t),v(x+c_{1}t))$ satisfying 
\mathcenter
\begin{equation*}
(u(-\infty),v(-\infty))=E_{1}\ \ \mbox{and}\ \ (u(+\infty),v(+\infty))=E^{*}
\end{equation*}
for $c_{1}\geq \overline{c}$.
\end{theorem}

\begin{theorem} \label{Theorem-4.2}
Suppose, both the axial equilibrium $E_{2}$ and interior point $E^{*}$ exist such that $m\gamma^{2}>c\beta(\gamma+\alpha\beta)$ holds. Then, system (\ref{eq:diff1}) has a positive travelling wave solution $(u(x+c_{1}t),v(x+c_{1}t))$ satisfying 
\mathcenter
\begin{equation*}
(u(-\infty),v(-\infty))=E_{2}\ \ \mbox{and}\ \ (u(+\infty),v(+\infty))=E^{*}
\end{equation*}
for $c_{1}\geq \overline{c}$.
\end{theorem}

The waves in Theorem \ref{Theorem-4.1} are known as predator-invasion travelling waves as the state before invasions is $(1,0)$ where only prey exist, and the state after invasions is $E^{*}$ where predators persist (see Figure \ref{fig:7}). Similarly, the waves in Theorem \ref{Theorem-4.2} are called prey-spread travelling waves since the state before invasions is $(0, \beta/\gamma)$ where only predators exist and the state after invasions is $E^{*}$ where the preys persist.

\subsection{Travelling wavefront connecting the predator-free state $(1, 0)$}
The predator-free equilibrium point $(1, 0)$ is a saddle point in this system. We intend to find a travelling wave solution that connects $E_{1}$ with some other equilibrium point. Linearizing the nonlocal model around $(1, 0)$ we get
\mathcenter
\begin{equation} \label{eq:tw4.5}
\begin{aligned}
\frac{\partial\overline{u}}{\partial t}&= d_{1}\frac{\partial^{2}\overline{u}}{\partial x^{2}}-\overline{u}-\frac{c}{1+m}\overline{v}, \\
\frac{\partial\overline{v}}{\partial t}&= d_{2}\frac{\partial^{2}\overline{v}}{\partial x^{2}}+s\overline{v}    
\end{aligned}
\end{equation}
We try to obtain the solution of system (\ref{eq:tw4.5}) in form 
\begin{align*}
\begin{cases}
    \overline{u}&=e^{\lambda t}e^{-Kx}\widetilde{u}, \\
    \overline{v}&=e^{\lambda t}e^{-Kx}\widetilde{v},
\end{cases}
\end{align*}
where $K>0$. Now, the eigenvalues of the corresponding characteristic equation be $\lambda_{1}=-1+d_{1}K^{2}$ and $\lambda_{2}=s+d_{2}K^{2}$.
We focus on the larger eigenvalue $\lambda_{2}$, and the wave propagation speed $\overline{c}$ is given by 
$$\overline{c}(K)=\frac{\lambda_{2}}{K}=d_{2}K+s/K$$
and $\overline{c}(K)|_{\min}=2\sqrt{sd_{2}}$. So, it is observed that the minimum speed of propagation $\overline{c}(K)|_{\min}$ increases for $d_{2}>0$.

\subsection{Travelling wavefront connecting the prey-free state $(0, \beta/\gamma)$}
In the system, the prey free equilibrium point $E_{2}$ is a saddle point if $m\gamma^{2}>c\beta(\gamma+\alpha\beta)$. Now, we try to find a travelling wave solution that connects the prey-free state with other equilibrium points. Linearization of the nonlocal model around $E_{2}$ gives
\mathcenter
\begin{equation} \label{eq:tw4.6}
\begin{aligned}
\frac{\partial\overline{u}}{\partial t}&= d_{1}\frac{\partial^{2}\overline{u}}{\partial x^{2}}+b_{11}\overline{u}, \\
\frac{\partial\overline{v}}{\partial t}&= d_{2}\frac{\partial^{2}\overline{v}}{\partial x^{2}}+\frac{s}{\gamma}\overline{u}-s\overline{v}    
\end{aligned}
\end{equation}
with $b_{11}=1-c\beta(\gamma+\alpha\beta)/m\gamma^{2}$. We seek the solution of system (\ref{eq:tw4.6}) in the form 
\begin{align*}
\begin{cases}
    \overline{u}&=e^{\lambda t}e^{-Kx}\widetilde{u}, \\
    \overline{v}&=e^{\lambda t}e^{-Kx}\widetilde{v}, 
\end{cases}
\end{align*}
where $K>0$. The eigenvalues of the corresponding characteristic equation be $\lambda_{1}=b_{11}+d_{1}K^{2}$ and $\lambda_{2}=-s+d_{2}K^{2}$.
So, the larger eigenvalue is $\lambda_{1}$ along with the wave propagation speed $\overline{c}$: $$\overline{c}(K)=\frac{\lambda_{1}}{K}=d_{1}K+b_{11}/K$$
and $\overline{c}(K)|_{\min}=2\sqrt{b_{11}d_{1}}$. So, it is observed that the minimum speed of propagation $\overline{c}(K)|_{\min}$ increases for $d_{1}>0$. 

In a predator-prey system, the travelling wave speed denotes the velocity at which the spatial patterns of prey and predator species expand across the environment. Here, we have obtained the minimum wave speeds for the travelling waves connecting $E_{1}$ and $E_{2}$ and observed that the wave speeds are influenced here by some of the system parameters as well as the self-diffusion parameters. In the numerical part below, we validate the analytical results and show how the predators' cooperative hunting rate impacts system dynamics.

\section{Numerical Results}\label{sec:5}

This section contains the dynamical scenarios of the temporal as well as the spatio-temporal model and analyses how the cooperative hunting strategy adopted by a predator species affects the overall dynamics. We have mentioned earlier that the generalist predator is considered in the model to see the overall predator-prey dynamics in the presence of alternating food sources for the predator. As a result, the considered system may have different types of equilibrium points, including the trivial, prey-free, and predator-free states. It is also mentioned that the system may have three coexisting equilibrium points depending on the number of positive roots of the cubic polynomial equation (\ref{eq:stabdet1}). In addition, we have fixed the system parameter values as $c=0.05, m=0.08, s=0.05, \gamma=0.08$, and $\beta=0.01$ to perform the simulation.

\subsection{Temporal behaviour of the system}

\begin{figure}[!htb]
     \centering
     \begin{subfigure}{0.32\textwidth}
         \centering
         \includegraphics[height=5cm]{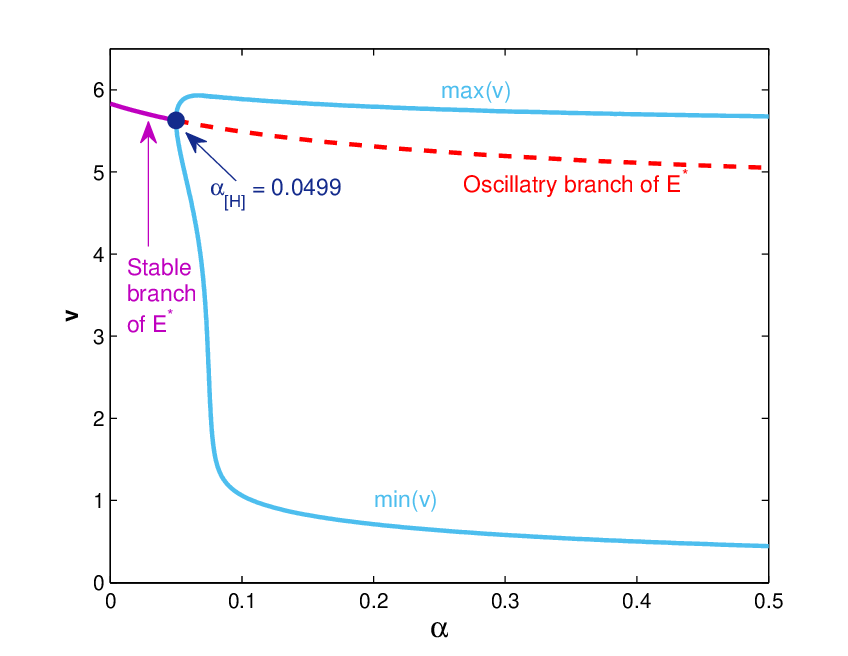}
         \caption{}\label{fig:1a}
     \end{subfigure}
     \begin{subfigure}{0.32\textwidth}
         \centering
         \includegraphics[height=5cm]{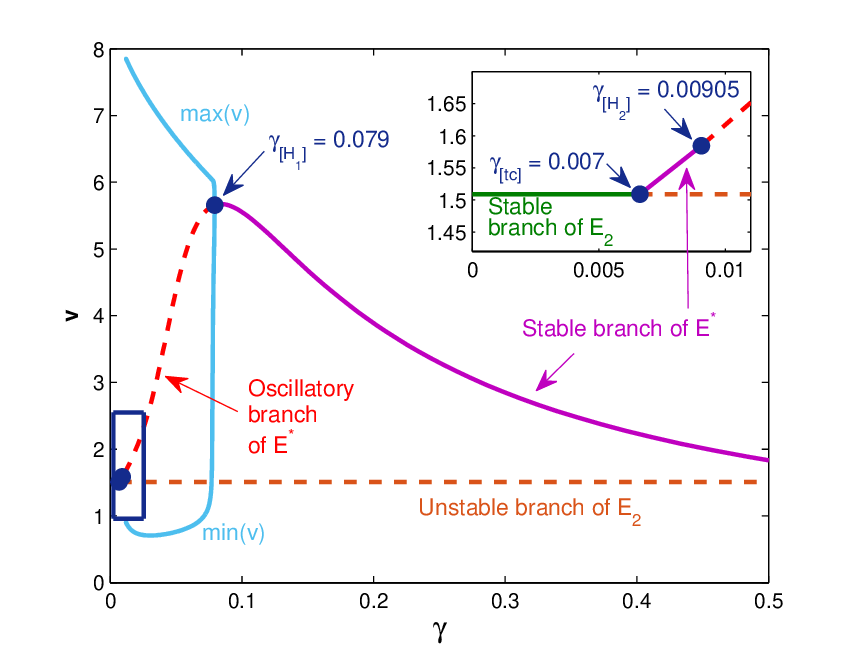}
         \caption{}\label{fig:1b}
     \end{subfigure}
    \begin{subfigure}{0.32\textwidth}
         \centering
         \includegraphics[height=5cm]{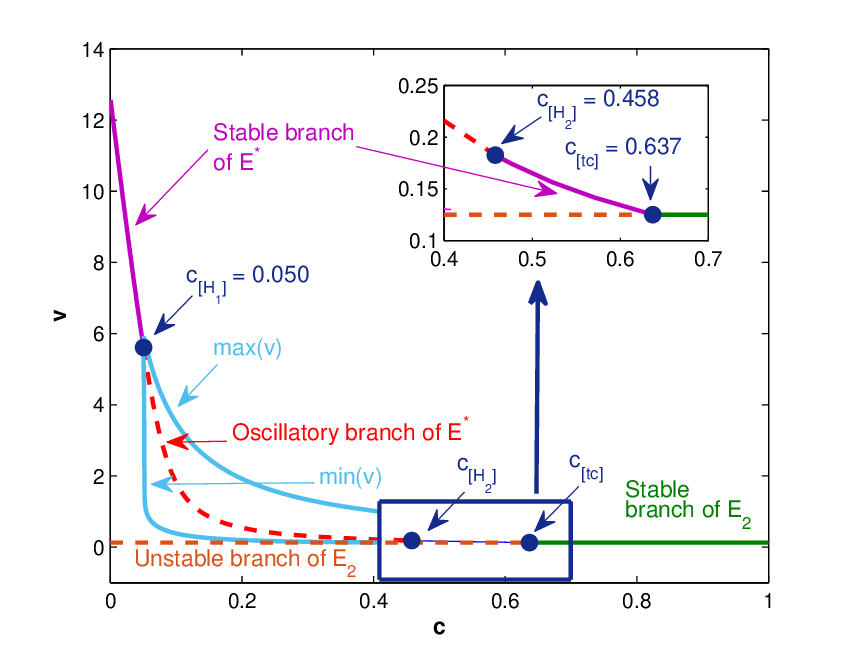}
         \caption{}\label{fig:1c}
     \end{subfigure}
\caption{Change of dynamical behaviour of temporal system (\ref{eq:det1}) with increasing (a) $\alpha$, (b) $\gamma$, and (c) $c$. } \label{fig:1}
\end{figure}

The system parameters used while formulating a model have their own importance in capturing the dynamic nature of the prey and predator populations. There are many psychological factors of both the prey and predators working in a system for coexistence. Predators' mutual cooperation while attacking prey is counted as one of such factors that help them to gain more success in hunting. Not only that, but their intra-specific competition because of lack of food or unfavourable environment also acts as another factor that significantly changes the overall dynamics of the interaction. However, in this work, our main intention is to explore the contribution of the predator's cooperative hunting strategy $(\alpha)$ on the system. We have demonstrated the temporal dynamics of the model by considering $\alpha$ as a bifurcation parameter in Fig. \ref{fig:1a}. In this case, both the population density changes with the increase of the cooperative hunting rate of predator $(\alpha)$. A stable coexistence state is observed up to a certain value of $\alpha$, and beyond that, the oscillatory coexistence appears through a supercritical Hopf bifurcation at $\alpha_{[H]}=0.0499$ with first Lyapunov coefficient $l_{1}=-0.049$. This scenario reveals that if the predators cooperate with each other maintaining a certain limit, then this cooperation benefits species cohabitation, but beyond that threshold value, it causes the system to oscillate, i.e., too much cooperation disturbs the population coexistence.

Though our primary goal is to investigate the model's dynamics by changing the cooperative hunting rate, we also show how the other factors influence the resultant dynamics. Here, the intra-specific competition of predator $(\gamma)$ plays a significant role in the system dynamics. For a fix $\alpha = 0.04$, the system does not have any coexisting equilibrium point for $\gamma = 0$; however, it emerges from the prey-free axial equilibrium $E_{2}$ through a transcritical bifurcation at $\gamma_{[tc]}=0.0066$ [see Fig. \ref{fig:1b}]. An increase in the value of $\gamma$ leads to a stability switch between $E_{2}$ and $E^{*}$ through a transcritical bifurcation at $\gamma_{[tc]}=0.0066$, and $E_{2}$ becomes unstable for $\gamma>\gamma_{[tc]}$. We get a stable interior equilibrium from this situation, but the stability does not continue for a larger range as the system shows oscillatory behaviour once $\gamma$ exceeds a certain threshold value. It is observed that the system exhibits a subcritical Hopf bifurcation around $E^{*}$ at $\gamma_{[H_{2}]}=0.00905$ with $l_{1}=0.155$. Although it becomes unstable for $\gamma>\gamma_{[H_{2}]}$, the further increment in $\gamma$ leads to a situation where the species coexist again as a steady state. The stability exchange occurs through a supercritical Hopf bifurcation at $\gamma_{[H_{1}]}=0.079$ with $l_{1}=-0.049$. This stable scenario continues for $\gamma>\gamma_{[H_{1}]}$. The figure implies that a certain amount of competition among predators is needed to make the system stable.


Along with $\alpha$ and $\gamma$, the truncation in prey population size due to predator's hunt $(c)$ has its importance, which is presented in Fig. \ref{fig:1c}. Both the prey and predator population live in the system in a stable state when $c$ is chosen to be sufficiently small, i.e., in a situation where the predator relies on the secondary food source more than the targeted prey, but this scenario changes with increasing $c$. The population starts to show oscillatory behaviour when $c$ crosses a threshold value through a supercritical Hopf bifurcation at $c_{[H_{1}]}=0.05$ with $l_{1}=-0.0503$, and $E^{*}$ becomes unstable for $c>c_{[H_{1}]}$. However, this unstable behaviour does not continue for a larger consumption rate as the coexisting state again becomes stable through another Hopf bifurcation. The figure depicts that the system exhibits a subcritical Hopf bifurcation at $c_{[H_{2}]}=0.458$ with $l_{1}=10.67$. The continuous increase in consumption lowers the prey count at a striking rate and eventually reaches a prey-free stable situation. We have observed an exchange in stability in prey-free state $(E_{2})$ and interior state $(E^{*})$ through a transcritical bifurcation at $c_{[tc]}=0.637$ so that $E_{2}$ remains stable when $c$ goes above $c_{[tc]}$. Here, Fig. \ref{fig:1} reveals the fact that the prey and predator populations coexist in a system when the predators compete among themselves to consume an adequate amount of their targeted prey. On the contrary, a small amount of cooperation while hunting also has a positive contribution to coexistence, but if the cooperation rate starts to increase at a higher level, it ultimately disturbs the stability of the whole system.

\begin{figure}[!htb]
     \centering
     \begin{subfigure}{0.4\textwidth}
         \centering
         \includegraphics[height=5.5cm]{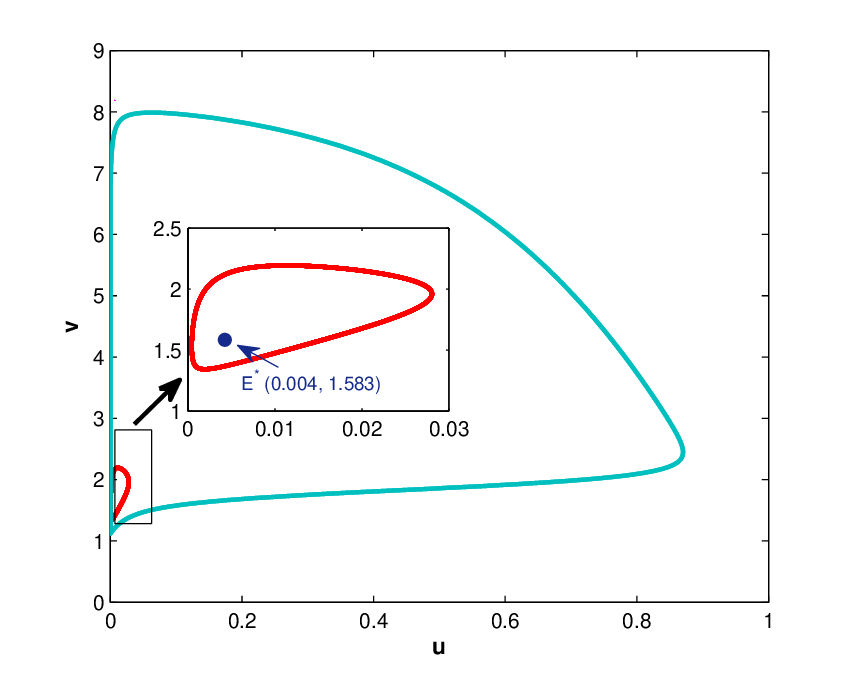}
         \caption{}\label{fig:2a}
     \end{subfigure}
    \begin{subfigure}{0.4\textwidth}
         \centering
         \includegraphics[height=5.5cm]{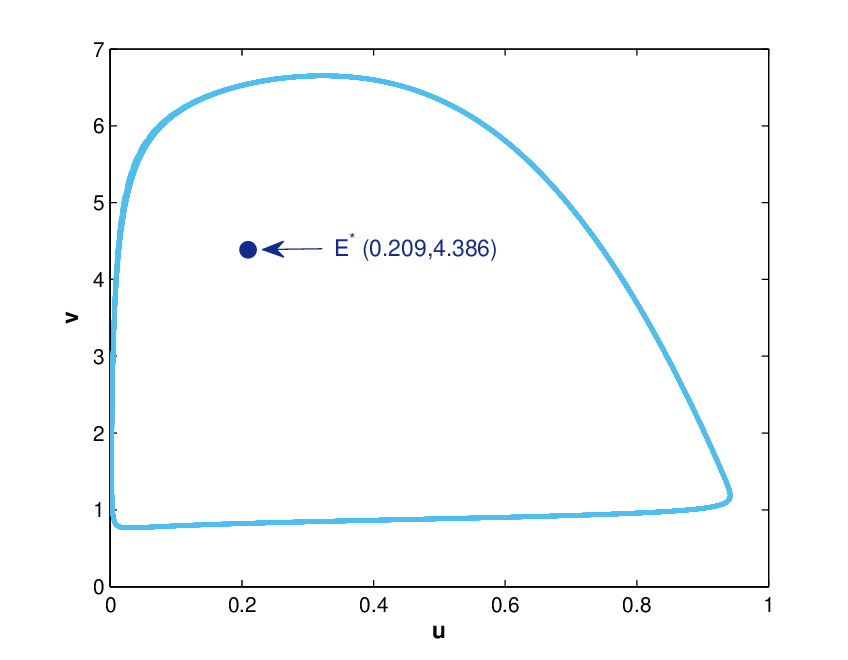}
         \caption{}\label{fig:2b}
     \end{subfigure}
\caption{Phase portrait of system (\ref{eq:det1}) for (a) $\gamma=0.009$ and (b) $\gamma=0.05$. The unstable and stable limit cycles are denoted by ({\color{red} \solidrule}) and ({\color{cyan} \solidrule}) colored curves, respectively. In Figure \ref{fig:2a}, multiple limit cycles are observed around stable $E^{*}$, and in Figure \ref{fig:2b}, a stable limit cycle occurs around unstable $E^{*}$.} \label{fig:2}
\end{figure}

As mentioned earlier, both subcritical and supercritical Hopf bifurcations have been observed while regulating $\gamma$, and so, we have presented the phase portraits of the proposed temporal system (\ref{eq:det1}) in Fig. \ref{fig:2} for different intra-specific competition rates of predators. From Fig. \ref{fig:1b}, it can be stated that a trajectory will converge to the stable interior equilibrium $E^{*}$ for $\gamma<\gamma_{[H_{2}]}$. But there will be an unstable limit cycle surrounding the stable interior point for an increased value of $\gamma$ (near the Hopf threshold $\gamma_{[H_{2}]}$). Fig. \ref{fig:2a} shows the existence of such an unstable limit cycle around $E^{*}$ for $\gamma=0.009$. In this case, the interior point will act as a stable point around which the unstable limit cycle occurs, but we have observed a stable limit cycle also surrounding the unstable cycle. This situation indicates the existence of multiple limit cycles in the system when $\gamma$ lies in $(\gamma_{[H_{2}]},\ \gamma_{[H_{1}]})$. Now, with the increase of $\gamma$, the unstable limit cycle starts to shrink, and ultimately it merges with the coexisting point, making it unstable. In this situation, $E^{*}$ loses its stability, and Fig. \ref{fig:2b} depicts that there is a stable limit cycle around the unstable $E^{*}$ when $\gamma$ is chosen as $0.05$. Again, for further increments, the limit cycle starts to shrink and ultimately disappears and settles down to a stable coexistence state for $\gamma>\gamma_{[H_{1}]}$. 

\begin{figure}[!htb]
     \centering
     \begin{subfigure}{0.3\textwidth}
         \centering
         \includegraphics[width=6.2cm,height=5.0cm]{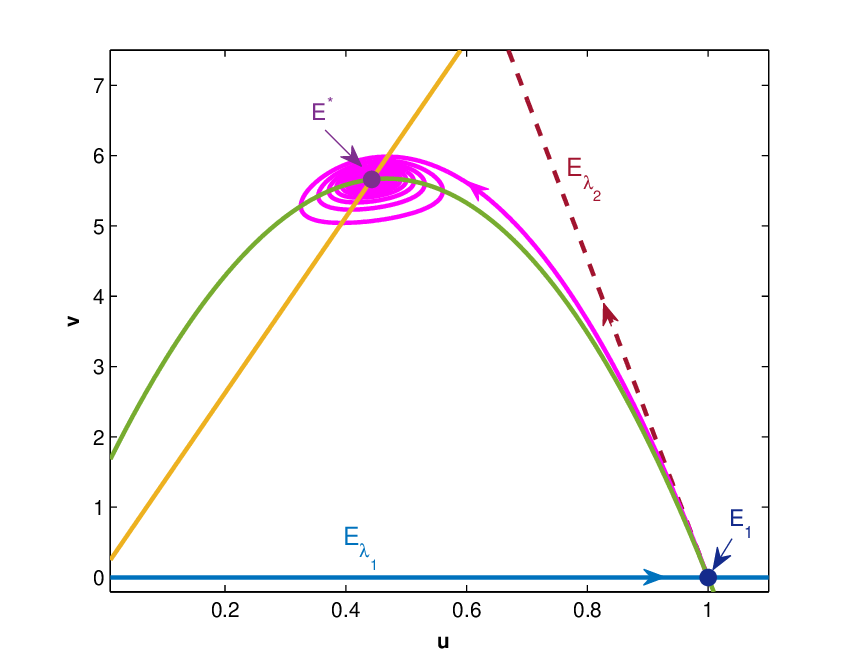}
         \caption{}\label{fig:3a}
     \end{subfigure}
     \begin{subfigure}{0.3\textwidth}
         \centering
         \includegraphics[width=6.2cm,height=5.0cm]{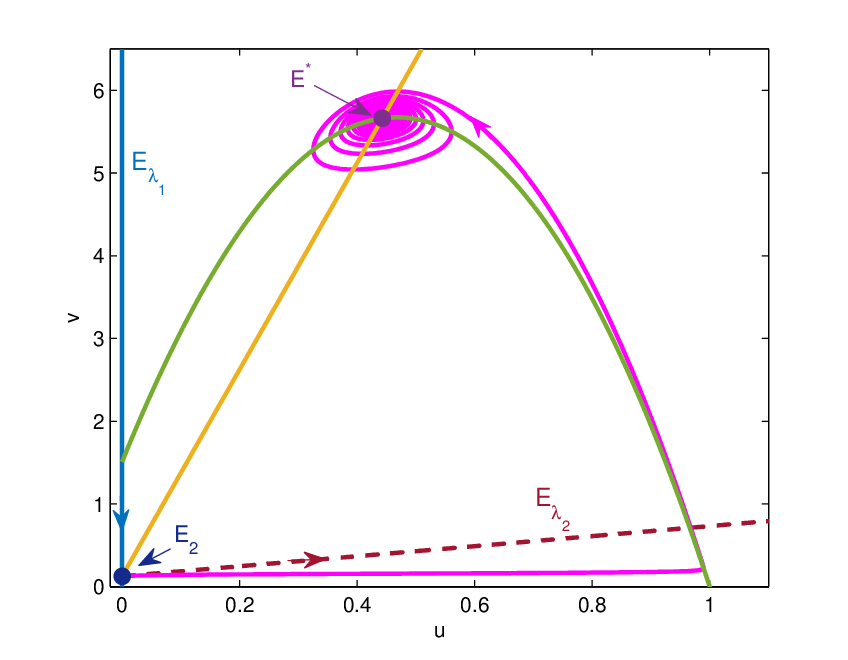}
         \caption{}\label{fig:3b}
     \end{subfigure}
     \begin{subfigure}{0.3\textwidth}
         \centering
         \includegraphics[width=6.2cm,height=5.0cm]{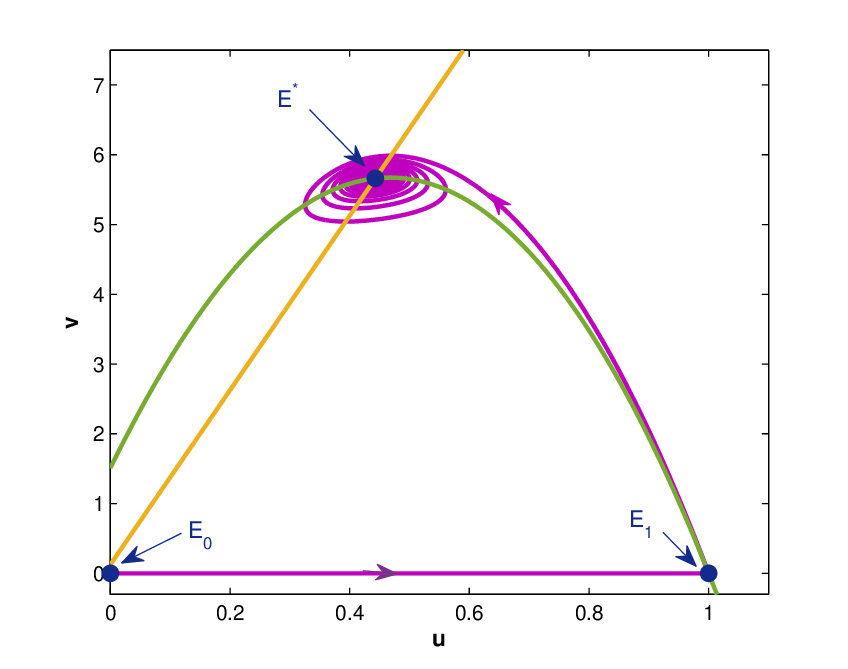}
         \caption{}\label{fig:3c}
     \end{subfigure}
\caption{Heteroclinic orbits joining the boundary equilibrium points $E_{1}$ and $E_{2}$ to stable interior equilibrium $E^{*}$ for $\alpha=0.04$. The non-trivial prey and predator nullclines are presented by ({\color{OliveGreen}\solidrule}) and ({\color{BurntOrange}\solidrule}) colored curves. The stable eigenspace in each figure is noted by ({\color{Aquamarine}\solidrule}) colored curve, while the unstable eigenspace is shown by ({\color{Maroon}\protect\dashedrule}) colored curve. } \label{fig:3}
\end{figure}

Though all of $\alpha, \gamma$, and $c$ have their influence in the system, we mainly intend to look into those dynamic changes that occur because of predator's cooperative hunting. The stability analysis already reveals that the predator-free equilibrium acts as a saddle point, whereas the stability of the prey-free state relies on a parametric restriction. In Theorem \ref{theorem-2.9}, we have already proven that in this system, there will be heteroclinic orbits joining $E_{1}$ and $E^{*}$ as well as $E_{2}$ and $E^{*}$. To support the statement, we have portrayed Fig. \ref{fig:3} where the existence of heteroclinic orbits joining the axial equilibrium points with $E^{*}$ is depicted. In Fig. \ref{fig:3a}, it shown that the trajectory, starting from $E_{1}$, ultimately converges to the homogeneous steady-state $E^{*}$ and the stable and unstable eigenspaces corresponding to $E_{1}$ are respectively represented as $E_{\lambda_{1}}=\{(u, v)\in\mathbb{R}^{2}: v=0\}$ and $E_{\lambda_{2}}=\{(u, v)\in\mathbb{R}^{2}: u+0.044v=1\}$. Moreover, for $\alpha=0.04$, the prey-free state $E_{2}$ acts as a saddle point, and Fig. \ref{fig:3b} shows that there is a heteroclinic orbit joining $E_{2}$ and $E^{*}$ in this situation where the stable and unstable eigenspaces are found as $E_{\lambda_{1}}=\{(u, v)\in\mathbb{R}^{2}: u=0\}$ and $E_{\lambda_{2}}=\{(u, v)\in\mathbb{R}^{2}: v-0.125=0.643u\}$. Startlingly, Theorem \ref{theorem-2.9} also states that there will be no direct connection from $E_{0}$ to $E^{*}$, but these will be connected through $E_{1}$ only, which is reflected in Fig. \ref{fig:3c}, where it is shown that if a trajectory, starts from $E_{0}$, first tends to $E_{1}$ and then converges to $E^{*}$. From a biological point of view, the continuous growth of prey species propels them to saturate at their maximum biomass, but because of the presence of the predator and their cooperation in the environment, the prey biomass decreases to make a steady coexistence. It is the saddle nature of $E_{1}$ that the predator-free state connects the population extinction and persistence in the environment. 

\begin{figure}[!htb]
     \centering
     \begin{subfigure}{0.4\textwidth}
         \centering
         \includegraphics[height=5.5cm]{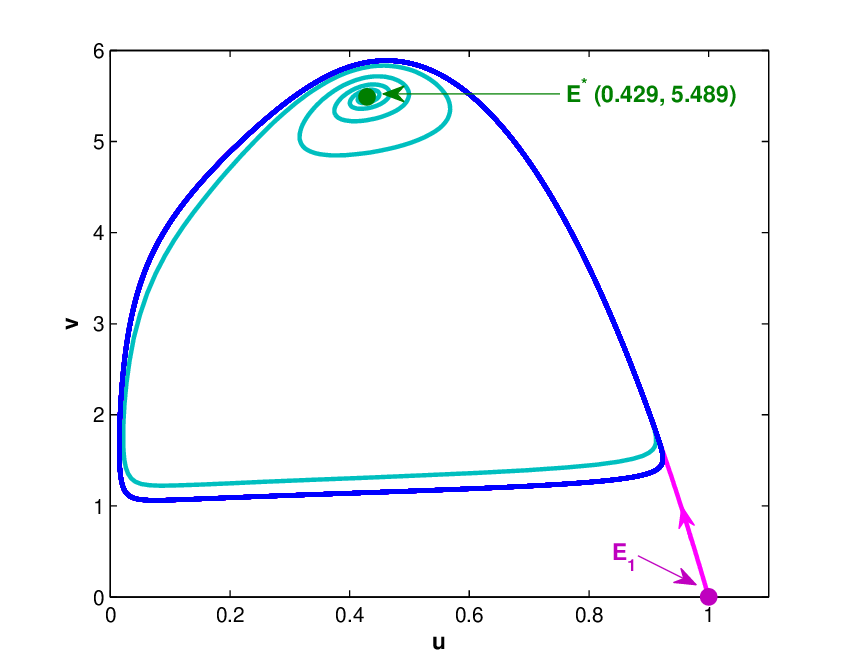}
         \caption{}\label{fig:4a}
     \end{subfigure}
     \begin{subfigure}{0.4\textwidth}
         \centering
         \includegraphics[height=5.5cm]{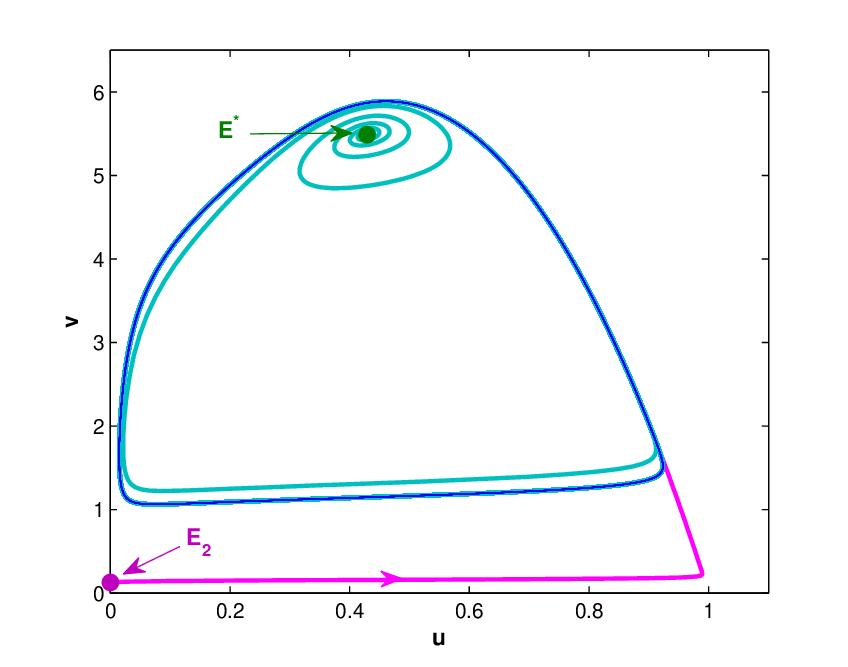}
         \caption{}\label{fig:4b}
     \end{subfigure}
\caption{The heteroclinic orbits joining the boundary equilibrium point $E_{1}$ and $E_{2}$ to a stable limit cycle. Here, $E^{*}$ noted by (\tikzcircle[OliveGreen, fill=OliveGreen]{3pt}) is an unstable point. } \label{fig:4}
\end{figure}

These heteroclinic orbits not only make connections between two equilibrium points but can also connect an equilibrium point with a limit cycle. We have already mentioned earlier that the population can show oscillatory behaviour once the cooperative hunting rate exceeds a threshold value through Hopf bifurcation [see Fig. \ref{fig:1}]. A stable coexisting equilibrium point occurs at $\alpha=0.04$, but for an increasing value of $\alpha$, it becomes unstable, and a stable limit cycle is generated through a Hopf bifurcation. In Fig. \ref{fig:4}, we have shown the existence of heteroclinic orbits for $\alpha=0.1$ joining the boundary equilibrium points $E_{1}$ [see Fig. \ref{fig:4a}] and $E_{2}$ [see Fig. \ref{fig:4b}] and the interior limit cycle which is a stable periodic orbit surrounding the unstable coexisting equilibrium $E^{*}$.

\subsection{Turing and non-Turing patterns}
We have outstretched the work by analyzing the pattern formation and travelling wave solutions of the corresponding spatio-temporal and nonlocal model and the significance of predators' cooperation in it. We mainly look into how this psychological factor impacts the spatial movement of the species. Before going to the travelling wave solution, let us first discuss the Turning patterns from scratch. Turing instability is one of the main factors that help to find non-homogeneous stationary patterns. To find such Turing instability conditions, the coexisting homogeneous steady-state has to be locally asymptotically stable. In the temporal model, there exists a Hopf bifurcation threshold $\alpha_{[H]}$ below which the coexisting equilibrium point is found to be stable, and the system shows periodic dynamics while exceeding this threshold [see Fig. \ref{fig:3a}].

\begin{figure}[!htb]
     \centering
     \begin{subfigure}{0.4\textwidth}
         \centering
         \includegraphics[width=7.5cm,height=5.5cm]{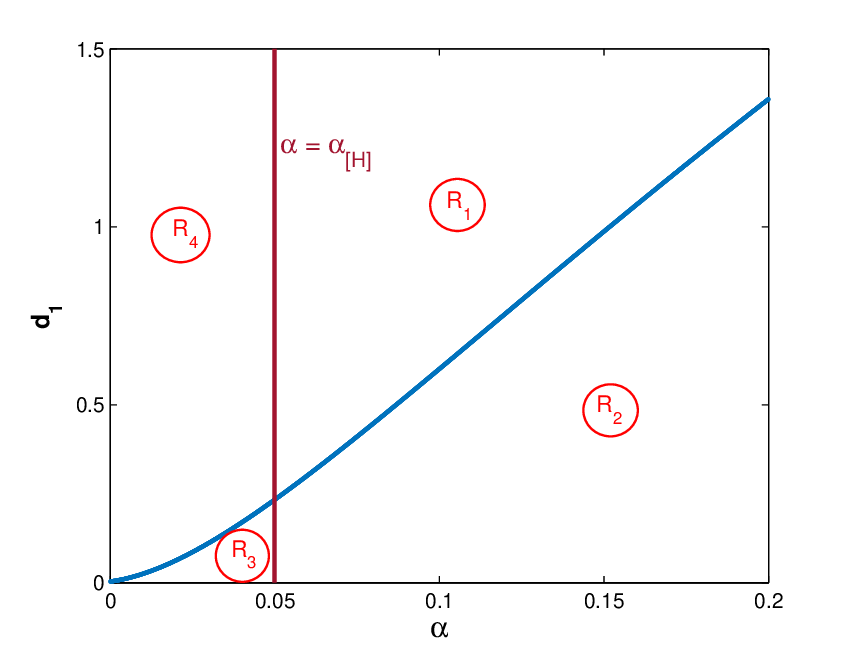}
         \caption{}\label{fig:5a}
     \end{subfigure}
     \begin{subfigure}{0.4\textwidth}
         \centering
         \includegraphics[width=7.7cm,height=5.5cm]{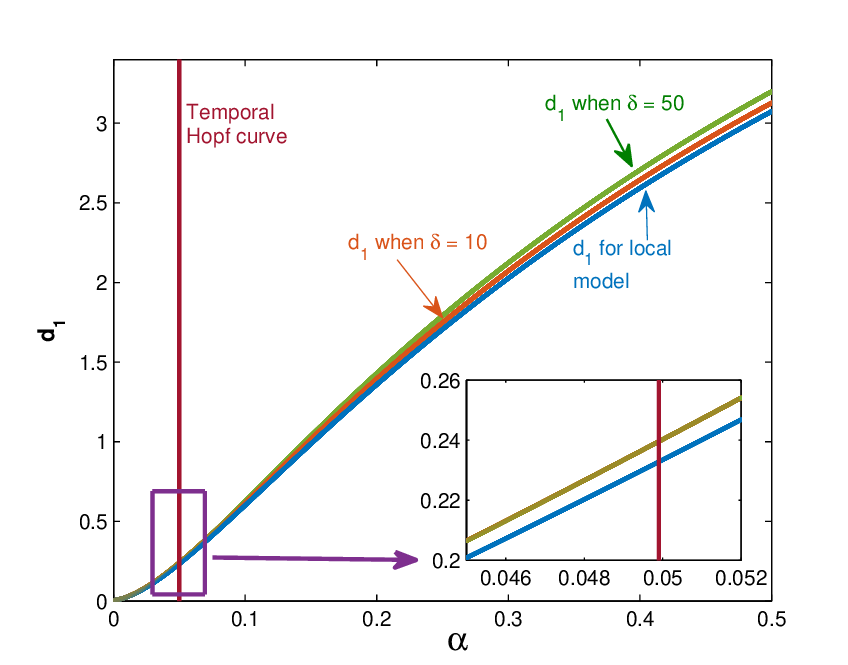}
         \caption{}\label{fig:5b}
     \end{subfigure}
\caption{(a) Temporal-Hopf and Turing bifurcation curve in the $\alpha$-$d_{1}$ plane for the local model. (b) Plot of Turing curves for local and nonlocal models. The ({\color[HTML]{A2142F} \rule[1mm]{8mm}{3pt}}) color curve denotes the temporal Hopf curve whereas the Turing curve for the local model is denoted by ({\color[HTML]{0072BD} \rule[1mm]{8mm}{3pt}}) color curve. The ({\color[HTML]{D95319} \rule[1mm]{8mm}{3pt}}) and ({\color[HTML]{77AC30} \rule[1mm]{8mm}{3pt}}) color curves are the Turing curves for the nonlocal model, i.e., for $\delta=10$ and $\delta=50$, respectively.} \label{fig:5}
\end{figure}

For a fixed $\alpha$, we obtain the Turing bifurcation threshold as $d_{1c}$, and we plot this set of Turing bifurcation thresholds for different values of $\alpha$ in the $\alpha$-$d_{1}$ plane [see Fig. \ref{fig:5a}]. In addition, we plot the temporal Hopf curve, which intersects the Turing curve and divides the whole region into four sub-regions. Now, the temporal model is stable for $\alpha<\alpha_{[H]}$ and oscillatory when $\alpha>\alpha_{[H]}$. The Turing domain $(R_{3})$ and homogeneous solution $(R_{4})$ lie on the left of the Hopf curve $\alpha=\alpha_{[H]}$, below and above the Turing curve, respectively. The bottom region $(R_{2})$ on the right of temporal Hopf curves is the Turing-Hopf domain, while the upper region $(R_{1})$ is the Hopf domain. So, for $\alpha>\alpha_{[H]}$, the spatio-temporal model can generate non-homogeneous stationary patterns and oscillatory solutions when $d_{1}<d_{1c}$ and $d_{1}>d_{1c}$ hold. But, for $\alpha<\alpha_{[H]}$, it can produce stationary Turing patterns or homogeneous solutions depending on $d_{1}<d_{1c}$ or $d_{1}>d_{1c}$. Here, we have chosen different values of $\alpha$ and $d_{1}$ from each of $R_{1}, R_{2}, R_{3}$ and $R_{4}$ domains to see the dynamical behaviour of the local spatio-temporal model. On the other hand, we have shown the influence of nonlocal interaction $(\delta)$ on the Turing threshold in Fig. \ref{fig:5b}. For the nonlocal model, the temporal Hopf bifurcation curve is independent of $\delta$. However, it is seen that the Turing bifurcation curve depends on $\delta$. We plot the Turing bifurcation curves for different values of the range of nonlocal interaction in this figure.  

We have observed an interesting scenario regarding the Turing curve for the nonlocal model. On the left of the temporal Hopf curve, the nonlocal Turing curve moves downwards for an increasing value of $\delta$, but it always remains above the Turing curve corresponding to the local model, i.e., $d_{1c}<d_{1c}^{T}|_{\delta=50}<d_{1c}^{T}|_{\delta=10}$ when $\alpha\leq \alpha_{[H]}$. However, when $\alpha$ exceeds the Hopf threshold, the Turning curve for the nonlocal model shifts upwards with the increase of $\delta$, i.e., $d_{1c}<d_{1c}^{T}|_{\delta=10}<d_{1c}^{T}|_{\delta=50}$ when $\alpha>\alpha_{[H]}$. For example, the Turing bifurcation threshold at $\alpha=0.04$ is $d_{1c}=0.1695$ for the local model, whereas for the nonlocal model, the threshold becomes $d_{1c}^{T}=0.1741$ when $\delta=10$ and $d_{1c}^{T}=0.17397$ when $\delta=50$, but at $\alpha=0.4$ the thresholds become $2.59$, $2.64$ and $2.70$, respectively. Overall, the introduction of nonlocal interaction through the cooperative hunting strategy expands the Turing domain a bit, which consequently increases the chance of the occurrence of stationary population patterns. 

\begin{figure}[!htb]
\centering
     \begin{subfigure}{0.4\textwidth}
         \centering
         \includegraphics[width=7cm]{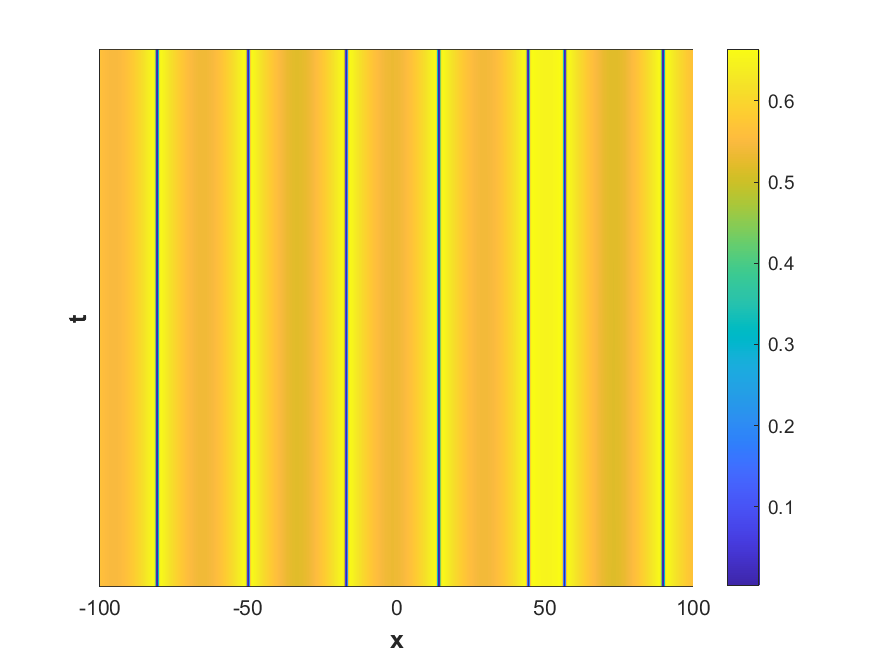}
         \caption{$\delta=0$}\label{fig:6a}
     \end{subfigure}
     \begin{subfigure}{0.4\textwidth}
         \centering
         \includegraphics[width=7cm]{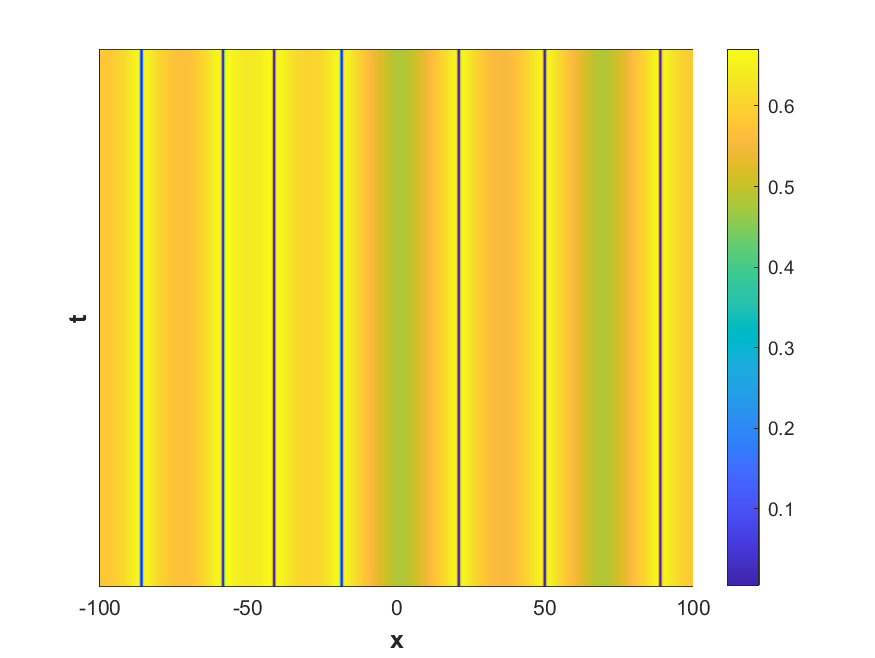}
         \caption{$\delta=50.0$}\label{fig:6b}
     \end{subfigure}
\centering
\caption{Turing patterns corresponds to the prey population $(u)$ of the (a) local and (b) nonlocal models for $(\alpha, d_{1})=(0.04, 0.01)$ with the mentioned initial condition around the homogeneous steady-state $E^{*}$.}\label{fig:6}
\end{figure}

\begin{figure} 
     \begin{subfigure}{0.3\textwidth}
         \centering
         \includegraphics[width=6.2cm,height=4.0cm]{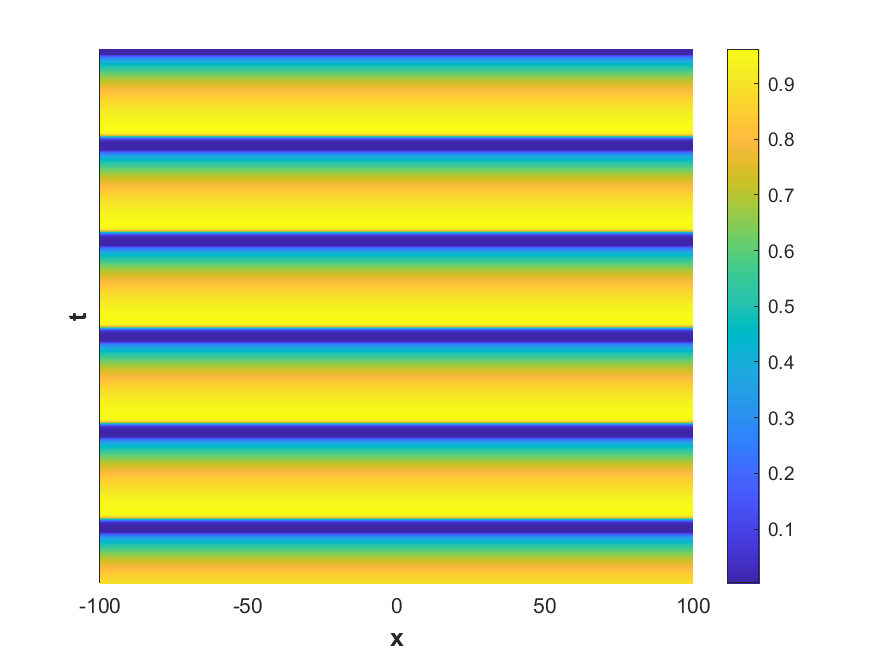}
         \caption{$(\delta, d_{1})=(0.0, 3.5)$}\label{fig:7a}
     \end{subfigure}
     \begin{subfigure}{0.3\textwidth}
         \centering
         \includegraphics[width=6.2cm,height=4.0cm]{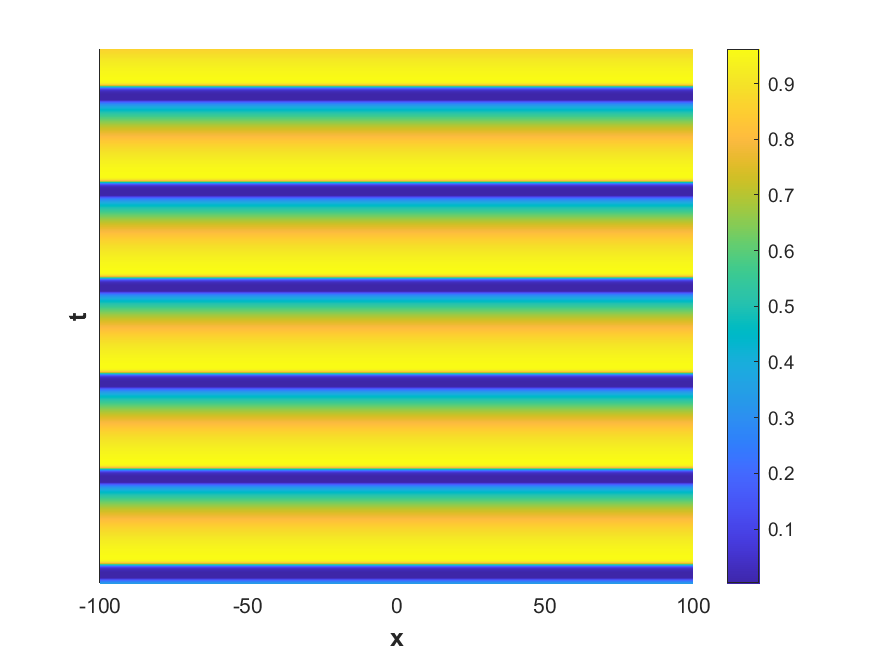}
         \caption{$(\delta, d_{1})=(0.0, 0.01)$}\label{fig:7b}
     \end{subfigure}
     \centering
     \begin{subfigure}{0.3\textwidth}
         \centering
         \includegraphics[width=6.2cm,height=4.0cm]{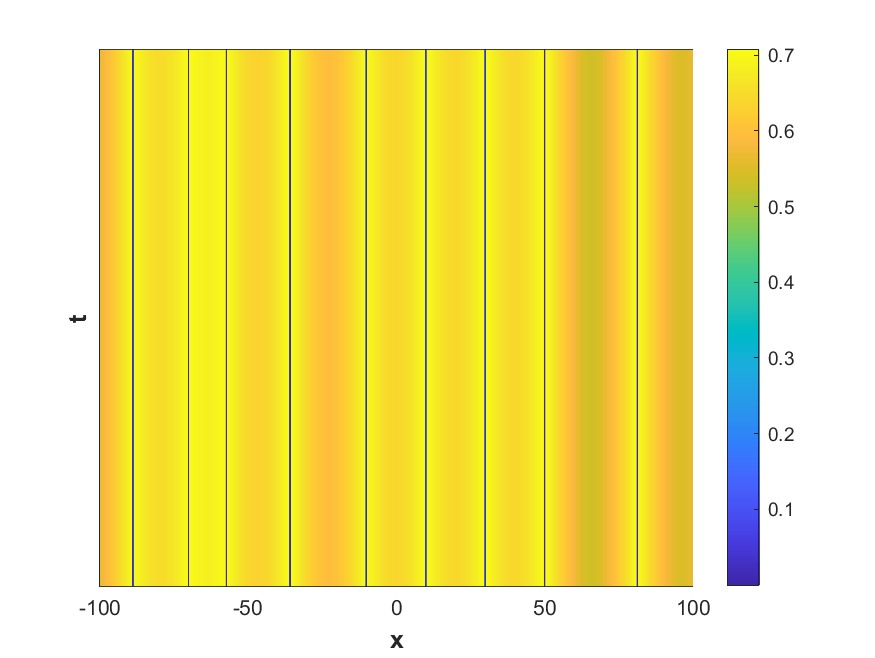}
         \caption{$(\delta, d_{1})=(0.0, 0.001)$}\label{fig:7c}
     \end{subfigure}
   \begin{subfigure}{0.3\textwidth}
         \centering
         \includegraphics[width=6.2cm,height=4.0cm]{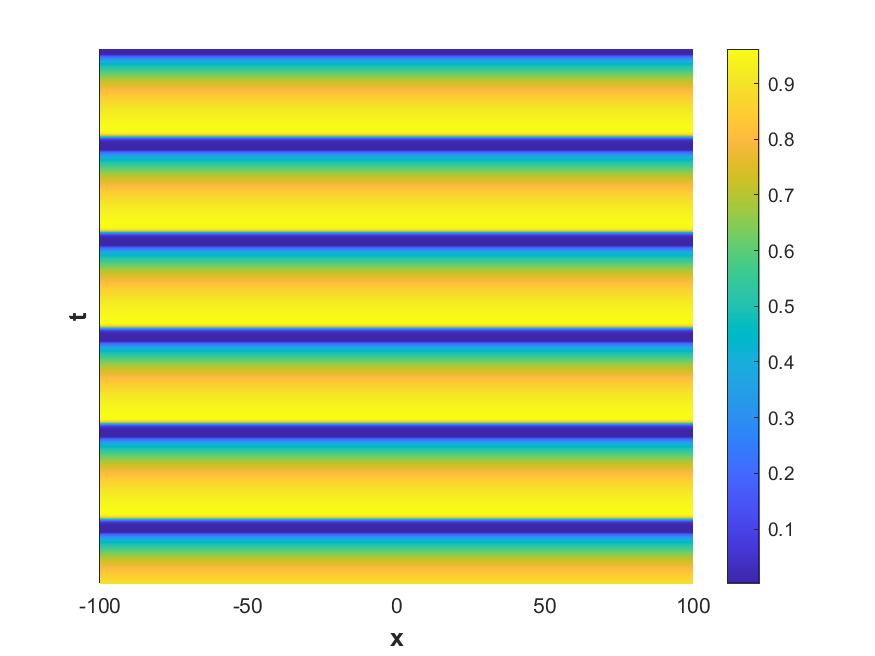}
         \caption{$(\delta, d_{1})=(50.0, 3.5)$}\label{fig:7d}
     \end{subfigure}
    \begin{subfigure}{0.3\textwidth}
         \centering
         \includegraphics[width=6.2cm,height=4.0cm]{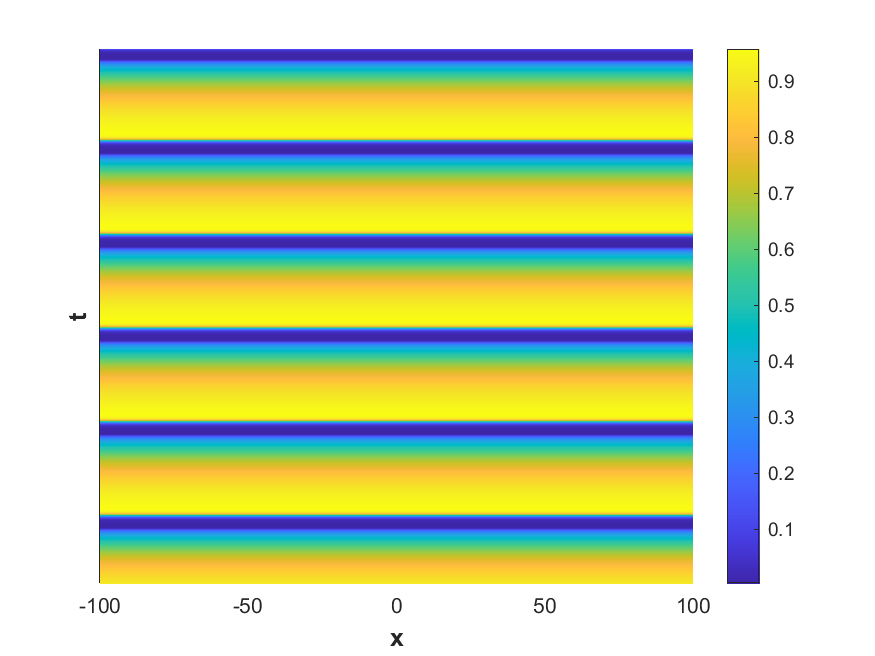}
         \caption{$(\delta, d_{1})=(50.0, 0.01)$}\label{fig:7e}
     \end{subfigure}
    \begin{subfigure}{0.3\textwidth}
         \centering
         \includegraphics[width=6.2cm,height=4.0cm]{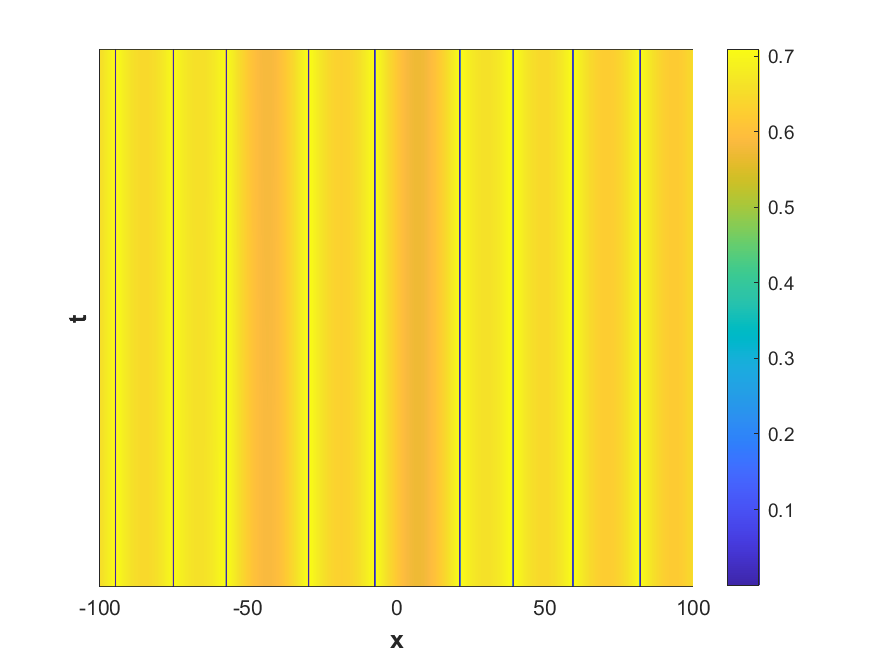}
         \caption{$(\delta, d_{1})=(50.0, 0.001)$}\label{fig:7f}
     \end{subfigure}
\caption{ Turing patterns corresponds to the prey population $(u)$ of the local and nonlocal models for $\alpha=0.4$ with the mentioned initial condition around the homogeneous steady-state $E^{*}$. The upper panel corresponds to the local model, whereas the lower row represents the nonlocal model. } \label{fig:7}
\end{figure}

To find the Turing and non-Turing patterns for the local and nonlocal models, we have chosen the spatial domain as $[-L, L]\equiv[-100,100]$, with periodic boundary conditions. In addition, a heterogeneous perturbation is given around the coexisting homogeneous steady-state as the initial conditions:
\begin{equation} \label{eq:5.1}
     u(x_{j},0)=u^{*}+\epsilon\xi_{j}, \ \ \ v(x_{j},0)=v^{*}+\epsilon\psi_{j},
\end{equation}
where $|\epsilon|=10^{-5}<<1$ and $\xi_{j}$ and $\psi_{j}$ are Gaussian white noise $\delta$-correlated in space. The dynamical behaviour of the proposed local and nonlocal models is explored in Figs. \ref{fig:6} and \ref{fig:7}. It should be noted that the nonlocal model (\ref{eq:loc1}) turns into a local model (\ref{eq:diff1}) if the range of nonlocal interaction $\delta$ tends to $0$. 

For $\alpha=0.04(<\alpha_{[H]})$ the temporal model (\ref{eq:diff1}) has the feasible interior equilibrium point $E^{*} = (0.443, 5.662)$, which is locally asymptotically stable. And, the Turing bifurcation threshold $d_{1c}$ is found to be $0.1695$ when $d_{2}=10$. Now, when $d_{1}>d_{1c}$ holds in the left of the temporal Hopf curve, the stationary homogeneous solution can be obtained in the stable domain $(R_{4})$. On the other hand, if $d_{1}$ is chosen from the Turing curve in this domain, then we will find homogeneous solutions even under the heterogeneous perturbations around the coexisting steady-state.

The Turing pattern for the local model (\ref{eq:diff1}) with prey diffusion coefficient $d_{1}=0.01$ is plotted in Fig. \ref{fig:7a}, and the corresponding situations for the nonlocal model is depicted in Figs. \ref{fig:7b}. It should be noted that the predator's diffusive coefficient $d_{2}$ is chosen as $10$ while performing the numerical simulation. In Figs. \ref{fig:7a}, \ref{fig:7d}, the oscillatory solutions are shown for local as well as nonlocal models when $\alpha>\alpha_{[H]}$ and $d_{1}$ is chosen from Hopf domain. Here, we will get the solutions that will be homogeneous in space but oscillatory with time. 

Previous research has demonstrated that when a specialist predator is accounted for in a predator-prey model, time-dependent spatial patterns emerge when the diffusion parameter is picked from a region slightly within the temporal Hopf domain \cite{petrovskii2003quantification}. In the Turing-Hopf domain, the Turing behaviour mainly dominates and creates stationary patterns. However, the Hopf behaviour also dominates and produces oscillatory patterns. These oscillatory solutions can be found in this domain in a small region near the Turing curve. Generally, non-homogeneous stationary patterns exist in most parts of the Turing-Hopf domain for the local model, which is depicted in Fig. \ref{fig:7c} for $(\alpha, d_{1}=0.4, 0.001)$. This scenario is observed not only for the local model but in the presence of nonlocal interaction also [see Fig. \ref{fig:7f}]. But, if we choose the value of $d_{1}$ from very near the Turning curve in the Turing-Hopf region, we will get an oscillatory solution in both local and nonlocal models [see Figs. \ref{fig:7b}, \ref{fig:7e}].

Now, in a predator-prey interaction, the predators are mainly focused on achieving their goals and improving their abilities to strike quickly at their opponents, and cooperating with each other while hunting is one of them. The prey, on the other hand, is worried about how to dodge the predator's attack and is constantly reacting instead of acting. So, it is their mindset that helps to form Turing patterns over space with time. For instance, Figs. \ref{fig:6} and \ref{fig:7} reveal that when the predator's cooperation rate is very low $(\alpha<\alpha_{[H]})$, the prey species form patches that are heterogeneous over space, but when the cooperation becomes higher, we get population patches that oscillate with time.

\subsection{Travelling wave solutions:}

Apart from observing the heterogeneous patterns, we intend to analyze the travelling wave solutions of the local as well as nonlocal models to observe the species invasion in the system and how the predator's cooperative hunting rate influences the species' incursion. This study will help us to understand how the species channel their psychology to act and react in a predator-prey interaction through travelling wave formation and connecting different equilibrium states.
In order to study the travelling wave solutions of the local and nonlocal models, we connect different combinations of the homogeneous solutions $(u_{1},v_{1})$ and $(u_{2},v_{2})$ by the followings:
\begin{equation} \label{eq:5.2}
    u(x,0)=\begin{cases}
        u_{1}~~ \text{if}~ |x|<L_{1}, \\
        u_{2}~~  \text{otherwise},
    \end{cases} \ \ v(x,0)=\begin{cases}
        v_{1}~~ \text{if}~ |x|< L_{1}, \\
        v_{2}~~  \text{otherwise},
    \end{cases}
\end{equation}
where $L_{1}$ is a positive constant and $|x|\leq L$ with $2L$ as the length of the spatial domain. In addition, the conditions (\ref{eq:5.2}) have been considered as the initial conditions for the numerical simulations of the travelling waves.

\begin{figure}[!htb]
     \centering
     \begin{subfigure}{0.4\textwidth}
         \centering
         \includegraphics[width=7cm,height=4.5cm]{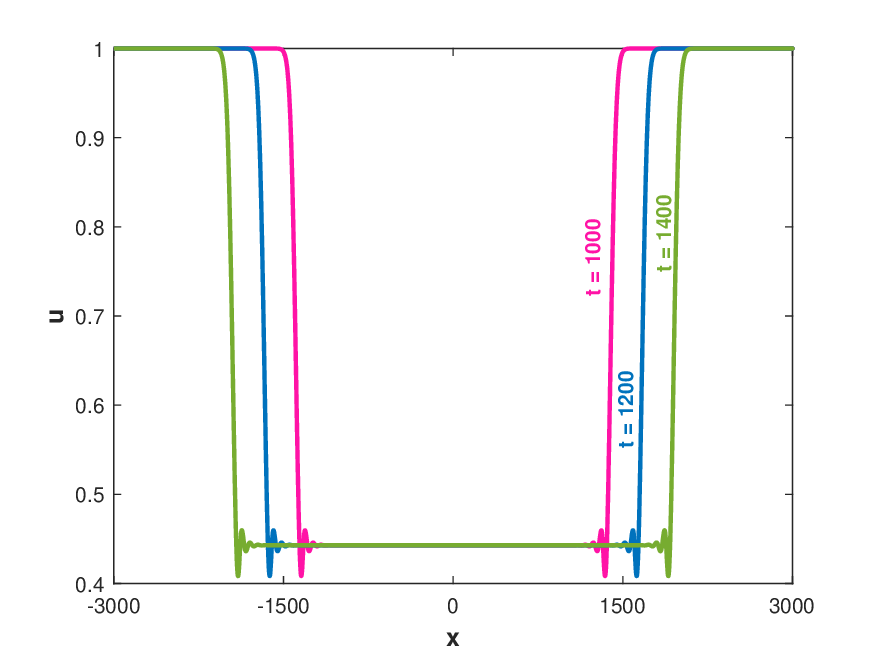}
         \caption{}\label{fig:8a}
     \end{subfigure}
     \begin{subfigure}{0.4\textwidth}
         \centering
         \includegraphics[width=7cm,height=4.5cm]{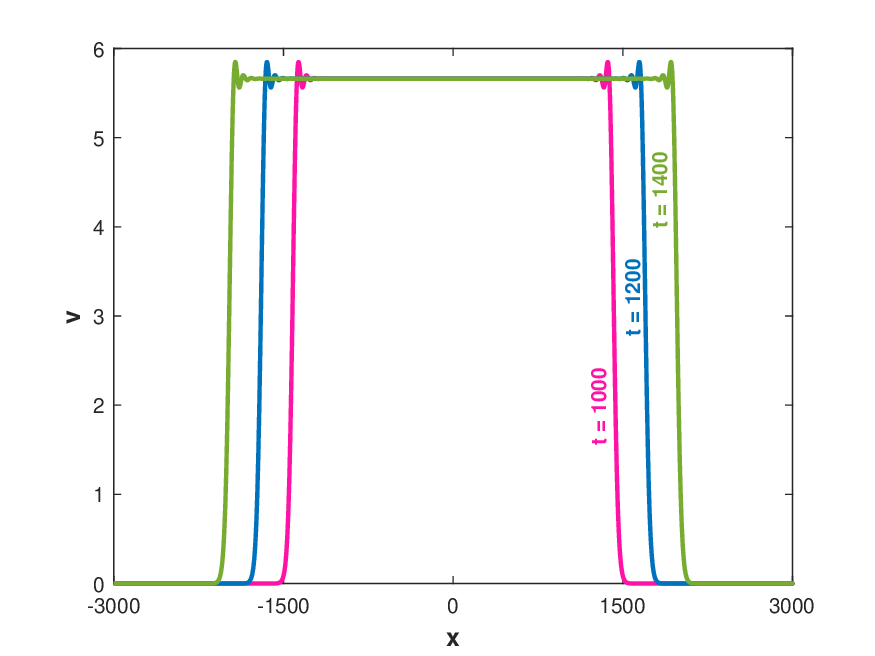}
         \caption{}\label{fig:8b}
     \end{subfigure}
\caption{Travelling wavefront corresponding to the (a) prey species $(u)$ and (b) predator species $(v)$ of the local model (\ref{eq:diff1}) connecting $E_{1}$ and stable equilibrium $E^{*}$ for $\alpha=0.04$ and $d_{1}=1$.} \label{fig:8}
\end{figure}

First, we find the travelling wave solution for the local model connecting the homogeneous solutions $E_{1}$ and $E^{*}$. Figure \ref{fig:8} depicts the travelling wave solution for the parameter values $\alpha=0.04$ and $d_{1}=1$ for three different time instances. In this case, we have observed that the travelling wave moves smoothly towards the boundaries without changing its shape. It is because the homogenous steady state $E^{*}$ is stable for the chosen parameter values. The symmetry of the travelling wave appears due to the symmetric initial conditions, and $L_{1}$ is chosen as 100 while performing the simulation. In addition, this travelling wave is found to be non-monotonic as the eigenvalues for the Jacobian matrix around the equilibrium point $E^{*}$ are complex conjugate with negative real parts. This result indicates that if prey and predator species are introduced in some areas of the habitat while the density of prey species is kept at its saturated level in the rest of the domain, then after a certain time, the predator species will be distributed over the domain and forces the prey species to reduce their density in all over the domain. Moreover, both the species survive for all the time.

\begin{figure}[!htb]
     \centering
     \begin{subfigure}{0.4\textwidth}
         \centering
         \includegraphics[width=7.5cm,height=5.5cm]{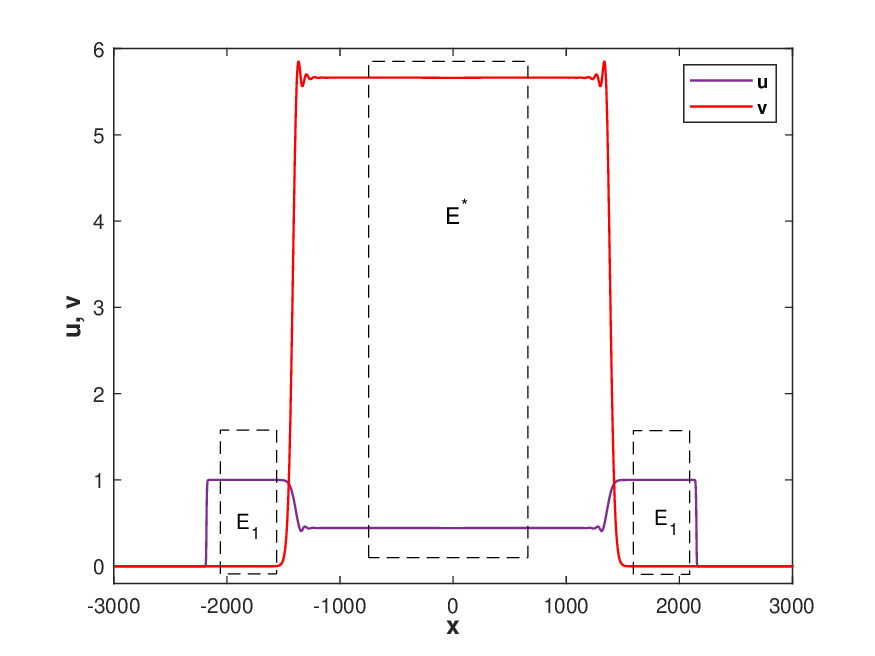}
         \caption{}\label{fig:9a}
     \end{subfigure}
     \begin{subfigure}{0.4\textwidth}
         \centering
         \includegraphics[width=7.7cm,height=5.5cm]{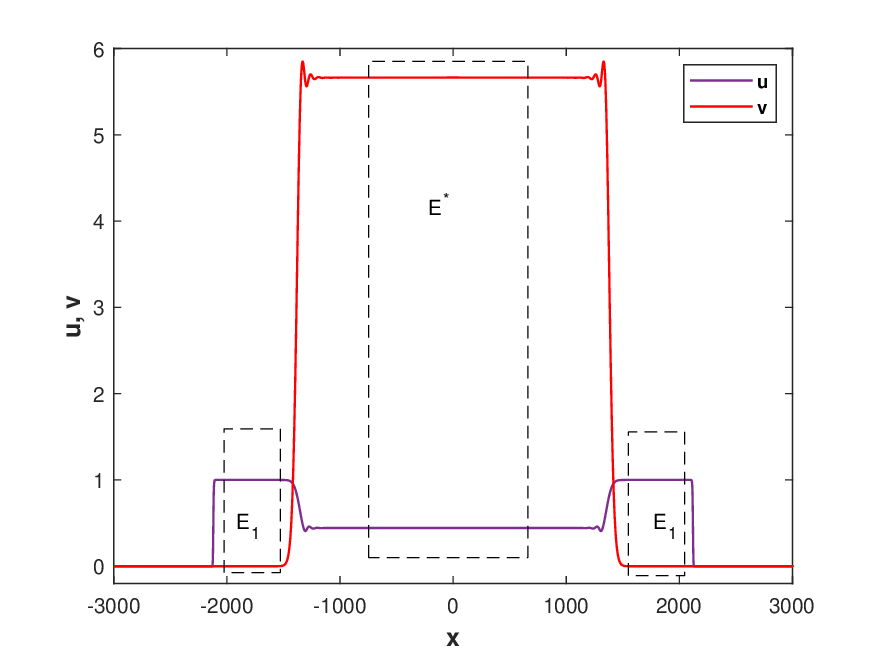}
         \caption{}\label{fig:9b}
     \end{subfigure}
\caption{The travelling of species in the (a) local model (\ref{eq:diff1}) and (b) nonlocal model (\ref{eq:loc1}).} \label{fig:9}
\end{figure}

Next, in Fig. \ref{fig:9}, it is observed that the travelling wave connecting $E_{0}$ and $E^{*}$ joins through $E_{1}$ only, i.e., we do not get any travelling wave that connects the states from population extinction and persistence. In fact, in Fig. \ref{fig:3}, it is already shown that if a trajectory starts at $E_{0}$, it will first tend to $E_{1}$ before converging to the coexistence state $E^{*}$. This situation arises not only for the local spatio-temporal model but also for the nonlocal model [see Fig. \ref{fig:9b}]. If the prey and predator species are introduced in some areas of a domain where the prey biomass is at its saturated (maximum) level, then the presence of predator will compel the prey to reduce their biomass with time. There is a chance of both species going extinct because of the sudden inadequacy of prey, or species coexistence. It is observed that the speeds at which the maximum prey biomass moves towards extinction and coexistence are different, and in fact, the speed is higher when it goes extinct. As a consequence, after some time, the prey biomass attains its maximum level at the boundary of the domain, but the predator starts to decline while going towards the boundary. This implies that predators follow normal invasions, but prey does not. The figure signifies that because of the predator's cooperative hunting and its generalist nature, the prey species, from maximum biomass, declines but coexists with the predator species instead of going to extinction.

\begin{figure}[!htb]
     \centering
     \begin{subfigure}{0.4\textwidth}
         \centering
         \includegraphics[width=7cm,height=4.5cm]{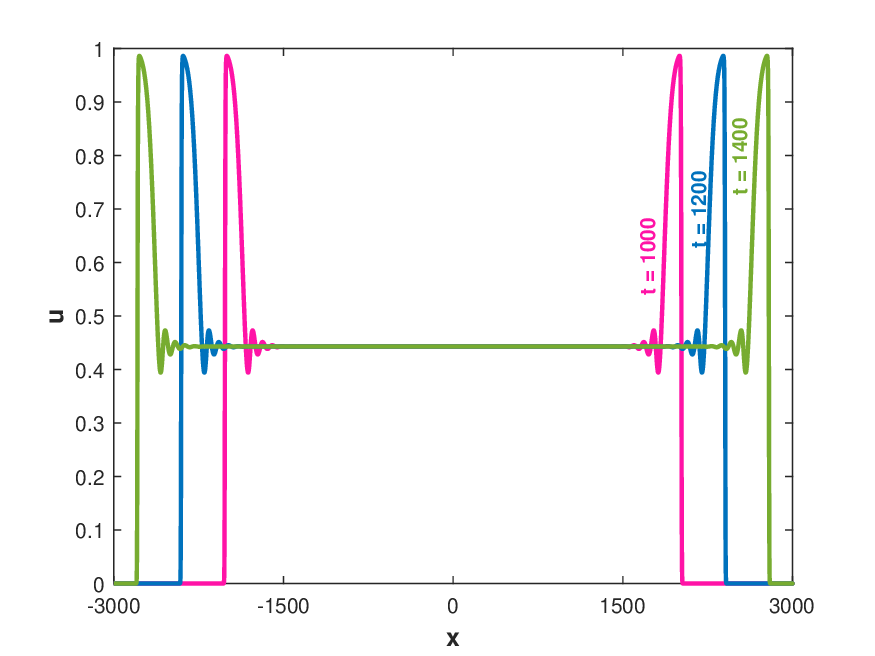}
         \caption{}\label{fig:10a}
     \end{subfigure}
     \begin{subfigure}{0.4\textwidth}
         \centering
         \includegraphics[width=7cm,height=4.5cm]{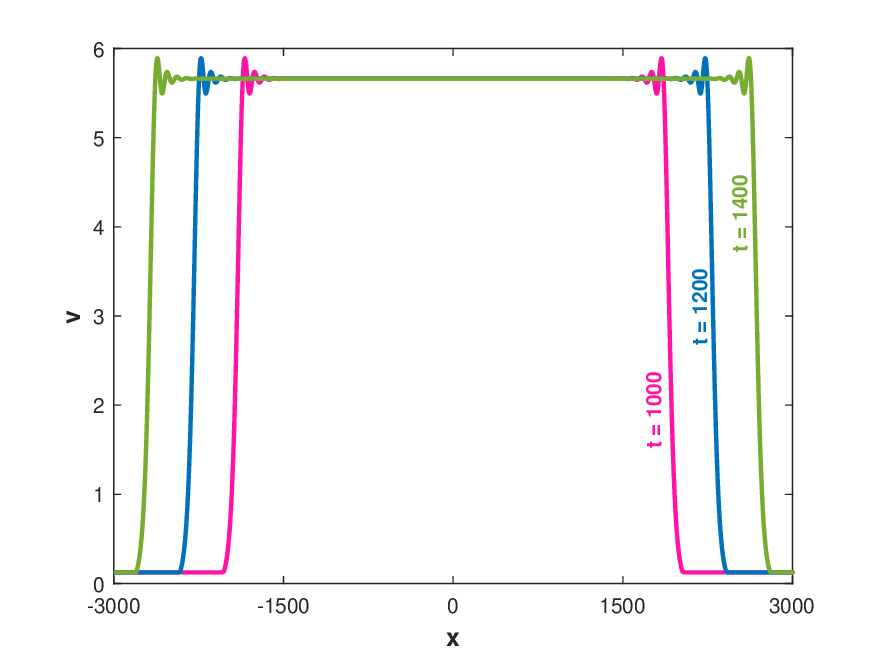}
         \caption{}\label{fig:10b}
     \end{subfigure}
\caption{Travelling wavefront corresponding to the (a) prey species $(u)$ and (b) predator species $(v)$ of the local model (\ref{eq:diff1}) connecting $E_{2}$ and stable equilibrium $E^{*}$ for $\alpha=0.04$ and $d_{1}=1$.} \label{fig:10}
\end{figure}

Now, we find the travelling wave solution connecting $E_{2}$ and $E^{*}$ for $\alpha=0.04$. In this case, $E_{2}$ is as a saddle point, and $E^{*}$ is locally asymptotically stable. Figure \ref{fig:10} depicts the travelling wave solutions for the local model when $d_{1}=1$. This scenario portrays that the travelling waves move smoothly toward the boundaries for the advancement of time, and local stability of $E^{*}$ is the main reason for the occurrence. In the figure, a transition zone from the prey-free state $E_{2}$ is observed to an increased level of prey density $u^{*}$ and a reduced level of predator density $v^{*}$. Moreover, this travelling wave is non-monotonic in this case too. From the ecological point of view, if prey and predator species are introduced in some areas of the habitat while the predator biomass is at a low level in the rest of the domain due to a shortage of targeted prey, then after a certain time, the predator species will be distributed over the domain. It is observed that the predator biomass increases significantly because of their cooperative hunting, but as they are provided with a secondary food source, the prey species also get a chance to grow up in the domain. Ultimately, the species coexist all the time. 

\begin{figure}[!htb]
     \centering
     \begin{subfigure}{0.4\textwidth}
         \centering
         \includegraphics[width=7cm,height=4.5cm]{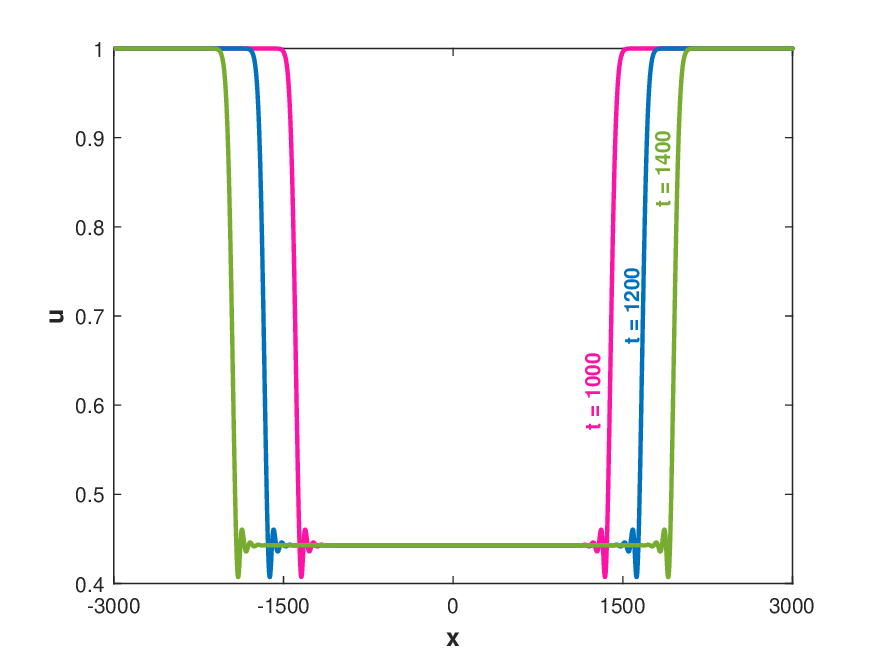}
         \caption{}\label{fig:11a}
     \end{subfigure}
     \begin{subfigure}{0.4\textwidth}
         \centering
         \includegraphics[width=7cm,height=4.5cm]{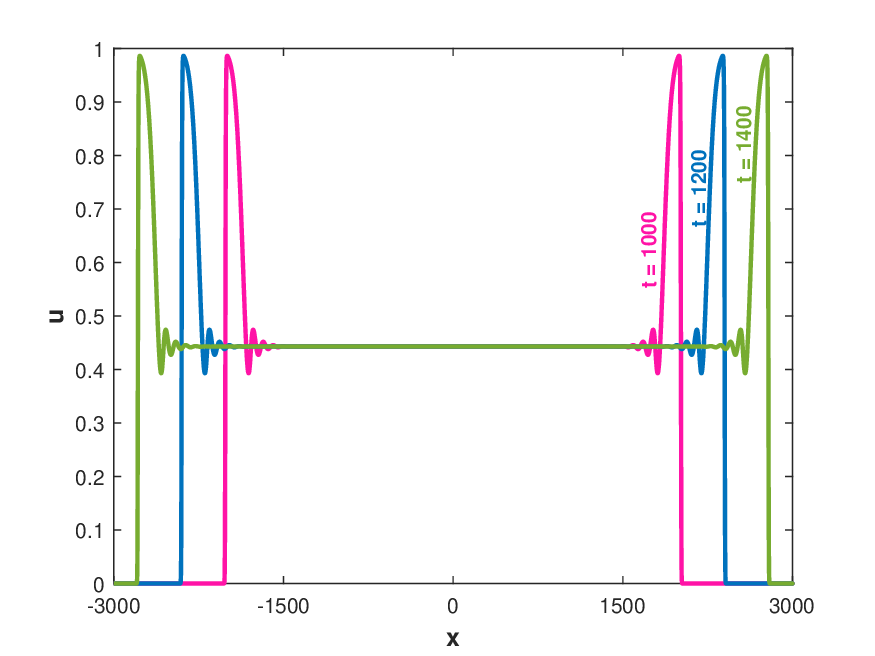}
         \caption{}\label{fig:11b}
     \end{subfigure}
\caption{Travelling wavefront corresponding to the prey species $(u)$ of the nonlocal model (\ref{eq:loc1}) connecting (a) $E_{1}$ and (b) $E_{2}$ with the stable equilibrium $E^{*}$ for $\alpha=0.04,\ \delta=50$ and $d_{1}=1$.} \label{fig:11}
\end{figure}

\begin{figure}[!htb]
     \centering
     \begin{subfigure}{0.4\textwidth}
         \centering
         \includegraphics[width=7cm,height=4.5cm]{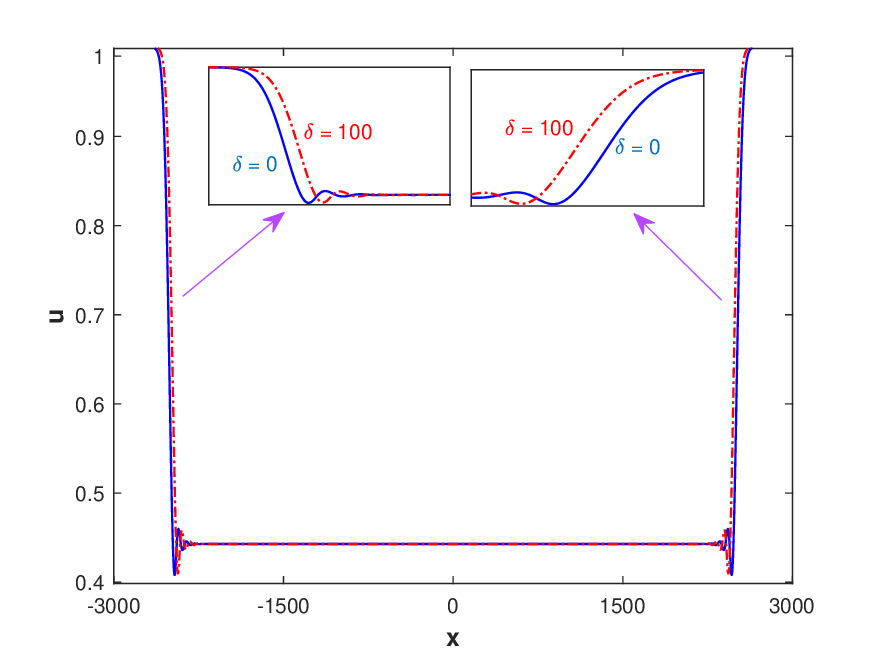}
         \caption{}\label{fig:12a}
     \end{subfigure}
     \begin{subfigure}{0.4\textwidth}
         \centering
         \includegraphics[width=7cm,height=4.5cm]{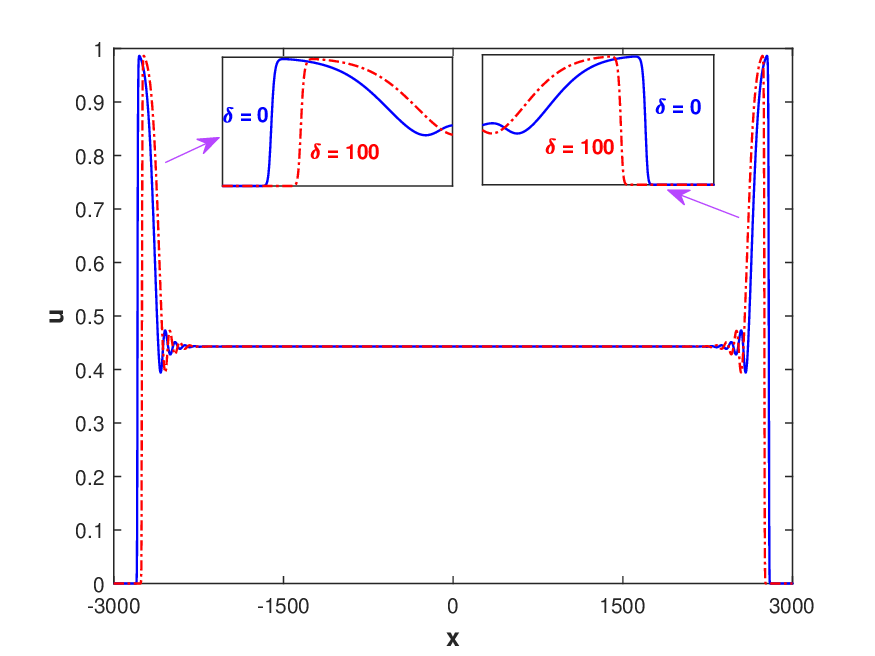}
         \caption{}\label{fig:12b}
     \end{subfigure}
\caption{Travelling wave solution corresponds to the prey species $(u)$ of the nonlocal model (\ref{eq:loc1}) for $(\alpha, d_{1})=(0.04, 1)$ connecting (a) $E_{1}=(1,0)$ and $E^{*}=(0.443,5.662)$; and (b) $E_{2}=(0,0.125)$ and $E^{*}=(0.443,5.662)$ for different values of nonlocal interactions.} \label{fig:12}
\end{figure}

Here, we see the effect of nonlocal interaction on the travelling wave solution. In Fig. \ref{fig:11}, we have plotted the evolution of the prey population $(u)$ for $\delta=50.0$ with the initial condition (\ref{eq:5.2}), and it is observed that the travelling wave solution exists for the nonlocal model connecting $E_{1}$ and $E^{*}$ as well as $E_{2}$ and $E^{*}$. It is depicted that the travelling wave moves smoothly towards the boundary when $d_{1}$ is chosen from the HSS domain. In Fig. \ref{fig:11a}, a transition zone is shown where the prey biomass is reduced from 1 to $u^{*}$, and the predator biomass is increased from 0 to $v^{*}$ with time. On the other hand, Fig. \ref{fig:11b} shows that when prey and predator species are introduced in some regions of the domain while a small amount of predator is there in the remaining region, then both the prey and predator increase from 0 to $u^{*}$ and $\beta/\gamma$ to $v^{*}$, respectively with time. The influence of nonlocal interaction is portrayed in Fig. \ref{fig:12} by comparing the travelling wave solution for the local and nonlocal models. Fig. \ref{fig:12a} shows the connection between $E_{1}$ and $E^{*}$ whereas Fig. \ref{fig:12b} depicts the connection between $E_{2}$ and $E^{*}$. The scenario reveals that the travelling waves take a longer time to distribute in the spatial domain for a higher nonlocal range of interaction. In addition, for $\alpha=0.04$, the estimated wave propagation speed of the wavefront connecting the predator-free state is $\overline{c}(K)|_{\min}=1.414$, while the wave propagation speed connecting the prey-free state is $\overline{c}(K)|_{\min}=1.92$. On the other hand, using numerical computation, we obtain the speed of the wavefront connecting the predator-free state is $1.403$ for the local model and $1.401$ when $\delta = 100$. In addition, the speed of the wavefront connecting the prey-free state is $1.9325$ for the local model and $1.9312$ when $\delta=100$. Hence, the computational and predicted wave propagation speeds are very close to each other.  

\begin{figure}[!htb]
     \centering
     \begin{subfigure}{0.4\textwidth}
         \centering
         \includegraphics[width=7cm,height=4.5cm]{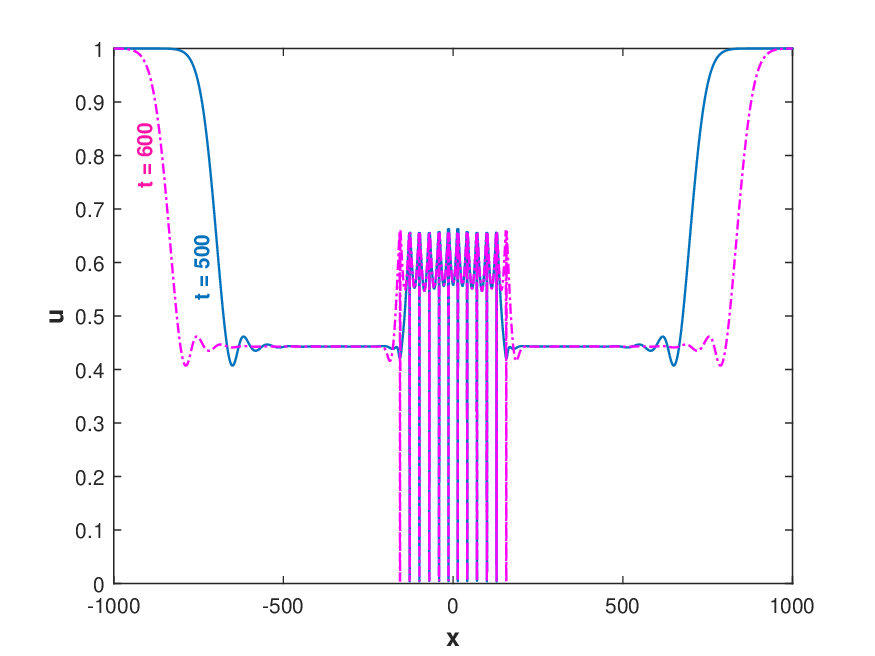}
         \caption{}\label{fig:13a}
     \end{subfigure}
     \begin{subfigure}{0.4\textwidth}
         \centering
         \includegraphics[width=7cm,height=4.5cm]{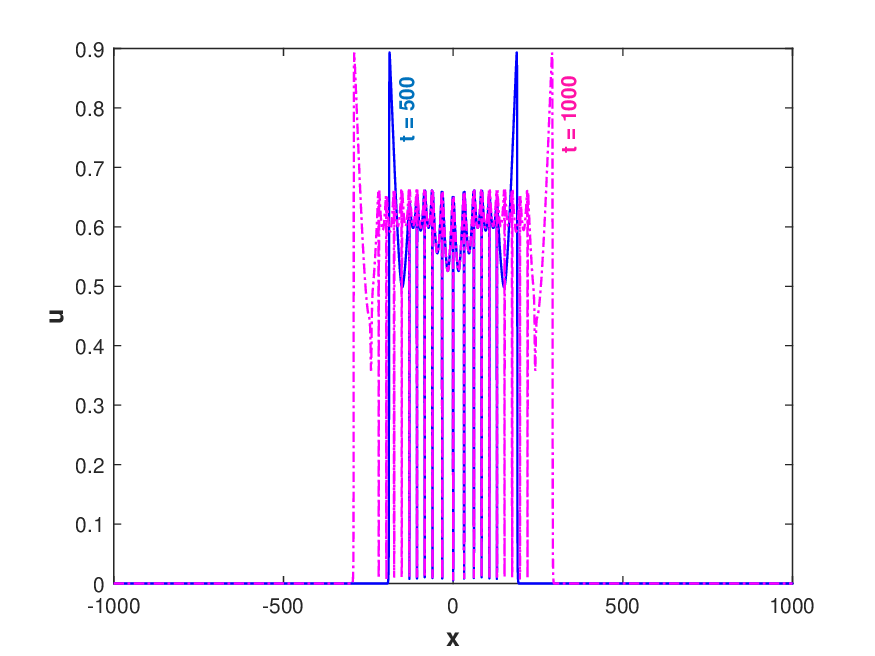}
         \caption{}\label{fig:13b}
     \end{subfigure}
\caption{Travelling wave solution corresponds to the prey species $(u)$ for the local model (\ref{eq:diff1}) for $\alpha=0.04$ and $d_{1}=0.01$ connecting (a) $(1,0)$ and homogeneous steady-state $(0.443,5.662)$; and (b) $(0,0.125)$ and $(0.443,5.662)$. } \label{fig:13}
\end{figure}

\begin{figure}[!htb]
     \centering
     \begin{subfigure}{0.4\textwidth}
         \centering
         \includegraphics[width=7cm,height=4.5cm]{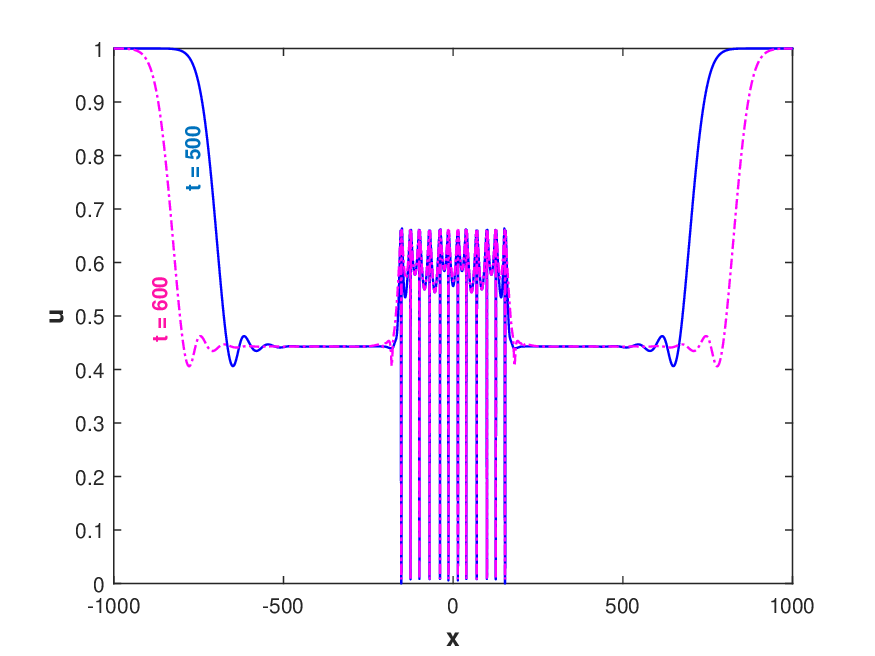}
         \caption{}\label{fig:14a}
     \end{subfigure}
     \begin{subfigure}{0.4\textwidth}
         \centering
         \includegraphics[width=7cm,height=4.5cm]{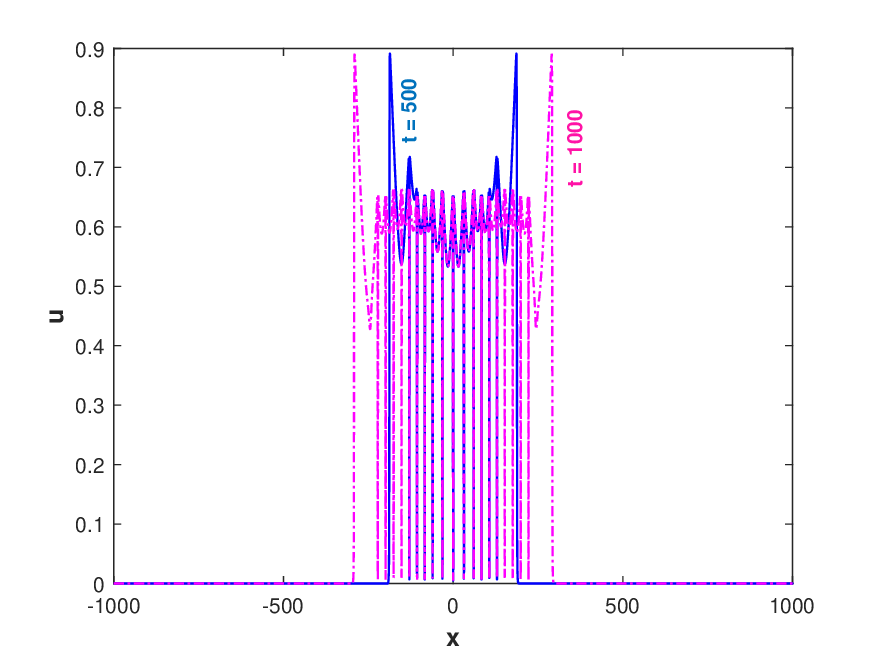}
         \caption{}\label{fig:14b}
     \end{subfigure}
\caption{Travelling wave solution corresponds to the prey species $(u)$ for the nonlocal model (\ref{eq:loc1}) for $\delta=50,\ \alpha=0.04$ and $d_{1}=0.01$ connecting (a) $E_{1}$ and $E^{*}$; and (b) $E_{2}$ and $E^{*}$. } \label{fig:14}
\end{figure}

Till now, we have studied the travelling wave solutions when $d_{1}$ lies in the HSS domain. Now, we choose the value of $d_{1}$ for which it satisfies the Turing instability conditions. Figure \ref{fig:13} portrayed the travelling wave solution for the local model for $\alpha=0.04(<\alpha_{[H]})$ with $d_{1} = 0.01$. The travelling wave solution is shown for two different time instances where it is observed that the waves move toward the boundary with time. It is shown when the prey diffusion coefficient $(d_{1})$ is chosen small enough (i.e., $d_{1}<d_{1c}$), the species do not move smoothly with the wave. It is the psychology of the species that they try to stay and form patches due to low diffusion parameter values, whereas the predator species try to form a patch with mutual agreement so that they can go for their targeted prey. Furthermore, their interactions produce an ordered structure out of random movement, which is guaranteed by the Turing instability. The same dynamics can be shown when the nonlocal interaction is taken into consideration, but the formation of patches takes a longer time to form in the presence of $\delta$ [see Fig. \ref{fig:14}]. This implies that when a predator cooperates with neighbourhood predators while attacking the prey, the colonization of species takes a longer time.

\section{Conclusions} \label{sec:6}

Psychological effects play a crucial role in shaping the dynamics of individual and collective behaviour. It impacts emotional regulation, group dynamics, stress and health, social skills, altruism and cooperation, etc. In addition, collective decision-making also plays a vital role in achieving a fruitful outcome for a whole group/herd, and cooperative hunting is counted as an example of that. Now, bio-social dynamics connects the interaction between biological and social factors to outline the population's behaviour and psychological processes. So, it is good to explore the population dynamics and the dynamic nature of a system under the influence of certain psychological effects. For example, mutual cooperation in a group or even in dyads can be considered as an outcome of a process where people's psychology flows in the same direction or they come to a mutual agreement to make the outcome fruitful. It not only benefits a single individual but becomes worthwhile for each member. In population dynamics, different psychological effects work directly or indirectly, leading the species to apply different strategies for survival. To make our analysis simpler, we have considered a predator-prey model as an exemplification where the growth of the prey population is affected by the cooperative hunting of the predator species.

The sustainability of predator species in an ecological system depends on their consumption process and the search strategy for prey. This consumption rate of predators depends on the adequacy of the targeted prey species, the level of mutual cooperation they impose, etc. So, the growth of prey becomes affected by the frequent attacks of their predator. Earlier research on the topic has already revealed that there are many species, such as lions, wolves, wild dogs, etc., that act as pack hunters with collective decision-making while attacking their prey \cite{stander1992cooperative, creel1995communal, schmidt1997wolf}.

In this work, we have formulated a predator-prey model emphasizing the psychological phenomena induced by the cooperative hunting strategy of the predators. The main intention here has been to elucidate the importance of this factor in the dynamic behaviour of the model. To achieve this, we have chosen the predator to be a specialist one to exclude the possibility of species extinction because of prey shortage. Also, we have looked upon the influence of predators' intra-specific competition, which causes psychological stress among them due to the scarcity of their food. The numerical results suggest that the cooperation among predators $(\alpha)$ as well as their intra-specific competition $(\gamma)$- both can be chosen as controlling factors as these parameters regulate the stability of the system. Figure \ref{fig:1a} shows that the population coexists in a stable state if the cooperation rate of predators crosses a threshold value through Hopf bifurcation. Furthermore, $\gamma$ has a dual role (stabilizing and destabilizing) in the model [see Fig. \ref{fig:1b}]. A stable coexisting state has been found for a very small as well as large value of competition $(\gamma<\gamma_{[H_1]}\ \mbox{and}\ \gamma>\gamma_{[H_2]})$, but oscillations have been observed when it lies within a range $\gamma\in(\gamma_{[H_1]},\ \gamma_{[H_2]})$. The figure reveals that $E^{*}$ is stable when $\gamma<\gamma_{[H_1]}$, but an unstable limit cycle starts to form around stable $E^{*}$ when $\gamma$ crosses the threshold through supercritical Hopf bifurcation. Around this unstable limit cycle, there will be a stable limit cycle too. The amplitude of the unstable limit cycle starts to decrease for increasing $\gamma$, and gradually $E^{*}$ becomes unstable. Only a stable limit cycle is observed around the unstable $E^{*}$ in this case, and once $\gamma$ crosses the supercritical Hopf bifurcation, the coexistence state becomes stable again. So, the existence of multiple limit cycles is observed in the system when $\gamma$ lies between two Hopf thresholds. Not only that, but the consumption rate also controls the dynamic nature of the system as the coexistence state switches its behaviour between stability and oscillation by regulating the parameter and ultimately tends to a prey-free situation. Now, in this model, species extinction is not an inherently stable situation. Constant growth of prey and predators, or other influences, disrupt the equilibrium easily. Here, we have shown that there will not be any heteroclinic orbit directly joining the population extinction state to the population persistence state. Nevertheless, we have shown that population extinction is connected to a predator-free state, and the coexistence state will be connected to a predator-free state. 

We have studied Turing and non-Turing patterns for the local and nonlocal models, along with travelling wave solutions. Figure \ref{fig:5} depicts that when the coexisting state acts as a stable equilibrium, the increase in $\alpha$ expands the region of the Turing domain, increasing the chances of non-homogeneous pattern formation in the spatio-temporal model. As the species are not always homogeneously distributed over a domain, this expansion will help in long-term species persistence. It indicates that the prey will move at a higher rate in the mentioned direction when the predators make a pack with mutual cooperation to hunt the prey. Furthermore, when the nonlocal terms are incorporated into the system, the chances of species colonization remain always higher when nonlocal interaction is considered in the system, but eventually, the chances of patch formation are reduced by increasing the range of nonlocal interaction. It indicates that when the predators start to cooperate with a distant neighbourhood, the chances of pattern formation reduce but remain higher than the case when no nonlocal interaction is considered. As the species can be colonized more with the nonlocal interaction, this will ultimately benefit the survival of both species in a longer time. Additionally, travelling wave solutions are observed when $d_{1}$ lies in the homogeneous solution and Turing domains. Also, the travelling waves take a longer time to distribute over the domain at a higher range of nonlocal interactions. It indicates that the species take more time to start colonizing when the predators from a neighbourhood region cooperate with each other while hunting their prey through a nonlocal model. 

There are many more things that can be incorporated as an extension of the present work to make it more realistic for biology. In the natural environment, the prey species may adopt different defence mechanisms as a counteraction to the predator's attack. Therefore, considering the contribution of group defence in the prey's growth will move the situation closer to reality. In addition, carrying out an analysis with other psychological effects, including herd/schooling behaviour, a fear effect, prey refuge, etc., with a nonlocal approach will lead to diversification in overall dynamics. Moreover, in ecological systems, the carryover effect may take place in any predator-prey interaction where species' past experiences and backgrounds are used to explain their present behaviour. Mutual cooperation among predators is an adaptive strategy that gets better with time, which can be carried over to the next generations and incorporated into the prey species of the model (\ref{eq:det1}). In this work, the boundary of the domain for the spatio-temporal model is chosen as periodic, but in future, the analysis can be performed with no-flux boundary conditions so that there will be no migration of any of the populations across the boundary of their habitat representing a scenario where the boundaries act as barriers or reflectors for the species. Also, it can be instructive to incorporate two-dimensional diffusion in the model instead of one to observe the existence of several homogeneous and non-homogeneous stationary patterns. Nonetheless, Gaussian noise can be incorporated into the system in order to depict the scenario where a predator's cooperation depends on the stochasticity of the environment. The environmental noise will help to observe uncertainty in parameter estimates and model structure, explore the stability and resilience of ecological systems, and simulate the unpredictable fluctuations that influence the persistence or extinction of populations. In addition to ecological modelling, epidemiological models can be explored with nonlocal models accounting for psychological effects as disease transmission may be affected because of people's behavioural changes. The methodology developed here, with more realistic aspects incorporated into the analysis, can also be helpful in such situations.

 \bibliography{P_References}

\end{document}